\documentclass[dvips,12pt]{article}

\usepackage{amsmath,amssymb,amsfonts,amscd,pstricks,pstricks-add,epsfig
,soul%
}
\renewcommand{\arraystretch}{1.5}

\def\bt{\begin{tabular}}
\def\te{\end{tabular}}

\def\BM{\begin{pmatrix}}
\def\EM{\end{pmatrix}}

\def\txt{\textstyle}
\def\ds{\displaystyle}

\def\cit{\text{\it I\hskip -6ptC\/}}
\def\ptcit{\hbox{${\scriptstyle I\hskip -4ptC\/}$}}

\def\square{\hfill\hbox{\vrule height .9ex width .8ex depth -.1ex}}
\def\rit{\text{\it I\hskip -2pt  R}}

\def\zit{\text{\it Z\hskip -4pt  Z}}
\def\qit{\text{\it I\hskip -5.5pt  Q}}

\def\nit{\text{\it I\hskip -2pt  N}}

\def\Bd{{\text B}}

\def\Ed{{\text E}}

\def\Ms{{\cal M}}

\def\be{\begin{equation}}
\def\ee{\end{equation}}
\def\beqn{\begin{eqnarray}}
\def\eeqn{\end{eqnarray}}
\def\nobeqn{\begin{eqnarray*}}
\def\noeeqn{\end{eqnarray*}}
\def\ba{\left(\begin{array}}
\def\ea{\end{array} \right) }

\def\bpr{\paragraph{Proof.}}

\def\epr{\square\vskip 6pt}
\def\eop{\hbox{\vrule height .9ex width .8ex depth -.1ex}}

\def\u{\underline}
\def\o{\overline}
\def\and{\; \mbox{and} \;}

\newcommand{\half}{\frac{1}{2}}

\def\Ker{\mathop{\rm Ker}\nolimits}

\def\mod{\mathop{\rm mod}\nolimits}

\def\Be{\begin{enumerate}}
\def\Ee{\end{enumerate}}

\def\Bena{\begin{enumerate}
\def\labelenumi{\theenumi)}
\def\theenumi{\arabic{enumi}}
\def\labelenumii{\theenumii)}
\def\theenumii{\alph{enumii}}}

\def\Bean{\begin{enumerate}
\def\labelenumii{\theenumii)}
\def\theenumii{\arabic{enumii}}
\def\labelenumi{\theenumi)}
\def\theenumi{\alph{enumi}}}

\def\Bero{\begin{enumerate}
\def\labelenumii{\theenumii)}
\def\theenumii{\arabic{enumii}}
\def\labelenumi{(\theenumi)}
\def\theenumi{\roman{enumi}}}

\def\BeRo{\begin{enumerate}
\def\labelenumii{\theenumii)}
\def\theenumii{\arabic{enumii}}
\def\labelenumi{(\theenumi)}
\def\theenumi{\Roman{enumi}}}

\def\Bi{\vskip 11pt\begin{itemize}\itemsep=18pt}
\def\Ei{\end{itemize}}

\def\Bd{\begin{description}}
\def\Ed{\end{description}}

\def\R{\right}
\def\L{\left}
\def\F{\frac}

\def\bigoplus{\mathop{\oplus}\limits}

\def\prod{\mathop{\Pi}\limits}
\def\sum{\mathop{\Sigma}\limits}

\def\BsM{\begin{smallmatrix}}
\def\EsM{\end{smallmatrix}}

\def\resp#1{(resp. #1)}
\def\rresp#1{\qquad \mbox{(resp.} \quad #1\ )}

\def\bbf{\boldmath\bf}
\def\o{\overline}
\def\wt{\widetilde}
\def\wh{\widehat}

\renewcommand{\arraystretch}{2}
\def\Bi{\begin{itemize}}
\def\Ei{\end{itemize}}

\newcommand{\MM}{\mathbb{M}\,}

\newcommand{\QQ}{\mathbb{Q}\,}

\def\tr{\operatorname{tr}}
\def\det{\operatorname{det}}
\def\Gal{\operatorname{Gal}}

\def\End{\operatorname{End}}
\def\mod{\operatorname{mod}}

\def\Hom{\operatorname{Hom}}
\def\Phom{\operatorname{Phom}}

\def\Re{\operatorname{Re}}

\def\FREPSP{\operatorname{FREPSP}}

\def\Aut{\operatorname{Aut}}

\def\ELLIP{\operatorname{ELLIP}}
\def\ellip{\operatorname{ellip}}

\def\GL{\operatorname{GL}}

\def\Aut{{\rm Aut}}

\def\bM{\begin{matrix}}
\def\eM{\end{matrix}}
\def\Bm{\L(\begin{smallmatrix}}
\def\Em{\end{smallmatrix}\R)}

\def\lr{left (resp. right) }

\def\bpr{\noindent{\bf{Proof\/}}:\;\;}
\def\To{\longrightarrow }

\def\RL{_{R\times L}}

\begin{document}

\setcounter{page}{0}
{\pagestyle{empty}
\null\vfill
\begin{center}
{\LARGE Generalized modular forms including\\[15pt]
 the weak Maass forms,\\[15pt] the Ramanujan's Theta functions\\[15pt] and the Tau function\\[15pt]}
\vfill
{\sc C. Pierre\/}
\vskip 11pt


\vfill
\vfill

\begin{abstract}
The modular forms are revisited from a geometric and an algebraic point of view leading to a geometric interpretation of the weak Maass forms connecting them to the Ramanujan Mock Theta functions and to the cusp forms generated from the Langlands global program.
\end{abstract}
\vfill

\end{center}

\eject

\tableofcontents
\vfill\eject}

\setcounter{page}{1}
\def\thepage{\arabic{page}}
{\parindent=0pt 
\setcounter{section}{0}
\section*{Introduction}
\addcontentsline{toc}{section}{Introduction}

The theory of automorphic forms was initiated by H. Poincare with the Fuchsian functions, pursued by F. Klein in connection with the elliptic functions and developed in the frame of analytic number theory by many mathematicians among which E. Hecke \cite{Hec1}, G. Shimura, C.L. Siegel, A. Ogg and J.P. Serre. But, according to the opinion of R. Langlands \cite{Lan}, ``{\bbf the theory of automorphic forms remains\/} in 1997 as it was in 1967: {\bbf a diffuse, disordered subject driven as much by the availability of techniques as by any high esthetic purpose\/}''.

On the other hand, F. Dyson, referring to Mock theta functions, claimed in 1987: ``{\bbf The Mock theta functions give us a tantalizing hints of a grand synthesis still to be discovered\/}.  Somehow, it should be possible to build them into a coherent group-theoretical structure.  This remains a challenge for the future''.
\vskip 11pt

Taking into account the great challenge of the program of Langlands, {\bbf the two-dimen\-sional modular forms \cite{Gel}, \cite{Ogg}, are revisited here from a geometric and an algebraic point of view in the sense that\/}:
\Bean
\item {\bbf the weight $k$ is proved to refer to $k-$dimensional modular forms\/};
\item {\bbf the level $N$ is the Galois  (or transcendence) extension degree of algebraic (or transcendental) quanta\/} \cite{Pie5}.
\Ee

The precedent advances in this field, developed in \cite{Pie1}, are:
\Bena
\item the Fourier series development of a two-dimensional cusp form decomposes into {\bbf a tower of increasing embedded semitori: geometric interpretation\/}.

\item a cusp form is a function into $\cit$ from a set of {\bbf $\rit^2/(\zit\big/ N \zit)^2$-lattices of transcendental quanta having extension degree $N$: algebraic interpretation\/}.

\item {\bbf a cusp form is covered by a global elliptic semimodule, which is a real analytic cusp form\/}, in such a way that every two-dimensional semitorus is covered by a set of semicircles: {\bbf cuspidal representation\/}.

\item {\bbf a cosemialgebra of dual cusp forms exists\/} with respect to the semialgebra of cusp forms {\bbf leading to a bisemialgebra of cusp biforms\/} (sections 1.11 and 1.12).
\Ee

All that constitutes the pieces of a theory of two-dimensional modular forms developed in chapters~1 and 2.

But, in order to enlarge this framework, I have also considered in chapter~3:
\Bena
\item {\bbf the $L$-functions as nonperiodic transforms of modular forms, the nonperiodicity reflected by the different unit frequencies or quanta of energy of their nontrivial zeros\/};

\item {\bbf the theta series as real analytic modular forms\/} in one-to-one correspondence with the respective complex modular forms;

\item {\bbf a geometric interpretation of the weak Maass forms\/} in the sense that {\bbf their series expansion consists of a sum of increasing embedded two-dimensional semitori with elliptic cross sections\/} in such a way that their nonholomorphic part corresponds to a shadow;

\item {\bbf the connections of the Ramanujan theta functions \cite{Ram} with weak Maass forms, modular forms and global elliptic semimodules\/} in such a way that:
\Be
\item the nonhomomorphic parts of the weak Maass forms correspond to the shadow of the Mock modular forms;
\item the partitions in Ramanujan's theta functions are partitions of quanta.
\Ee

\item {\bbf the (2-dimensional) cusp form of weight 12 and level 1 associated with the tau function and corresponding to a universal orthogonal cusp form\/}.
\Ee

More concretely, {\bbf the theory of modular forms of weight $k$ and level $N$ is reviewed in chapter~1\/} with a special emphasis on:
\Bean
\item the decomposition of the Poincare upper half plane $H$ into $\rit^2(\zit\big/N \zit)^2$-lattices partitioning it into coset representatives $t^{(2)}_{\zit_N}[n]$ of the quotient (semi)group 
\[
T^2_{\zit_N}(H)= T_2(\rit) \big / T_2(\zit \big/ N \zit)\;, \qquad 1\le n\le \infty \;, \quad n\in\nit\;, 
\]
$T_2(\rit)$ being the upper triangular group of matrices;

\item the action of $T_2(\rit)$  on {\bbf the fundamental domain $D_{nN}^{(2)}$, being a unitary measurable set, i.e. a period\/}, of $t^{(2)}_{\zit_N}[n]$:

\item the action of the Hecke operator $T_k(n)$ on {\bbf a two-dimensional cusp form $f_k(z_N)$ of weight $k$ and level $N$ decomposing it into Fourier series of which $n$-th term\/}, corresponding to the coset representative
$t^{(2)}_{\zit_N}[n]$, {\bbf is a two-dimensional semitorus $T^2_{nN}$ generated by two orthogonal semicircles at $n$ transcendental quanta\/}.
\Ee

{\bbf In order to be able to vary the weight of the two-dimensional cusp form $f_k(z_N)$ of weight $k$ and level $N$, the origin of this weight $k$ is assumed to be related\/}
\Bi
\item {\bbf to a $k$-dimensional cusp form $\phi ^{(k)}(z_{Nk})$\/}, having a decomposition in Fourier series
\[ \phi ^{(k)}(z_{Nk})=\sum_{n=1}^\infty  \lambda ^{(k)}_n\ e^{2\pi inz_{Nk}}\]
where $z_{Nk}$ is a $k$-tuple of complex numbers of order $N$ as developed in the global program of Langlands over number fields \cite{Pie2},

\item {\bbf by a projective map\/}:
\[ C\ P_{(k)\to 2}: \qquad \phi ^{(k)}(z_{N_k})\To f_k(z_N)\]
under the conditions of proposition~1.14,

\item {\bbf in such a way that the fundamental domain $D^{(2)}_{nNk}$ of $f_k(z_N)$ be equal to the fundamental domain $D^{(k)}_{nN}$ of $\phi ^{(k)}(z_{Nk})$\/}.
\Ei

The morphism of projection from a two-dimensional cusp biform $(f_k(z^*_N)\times f_k(z_N))$ of weight $k$ and level $N$ to a cusp biform
$(f_h(z^*_N)\times f_h(z_N))$ of weight $h$ and some level $N$, $k>h$, $k$ and $h$ being even integers, depends on :
\Bean
\item the morphism of projection $\Phom_{k\to h}$ from a $k$-dimensional cusp biform to a $h$-dimensional cusp biform;

\item the Langlands functoriality conjecture \cite{Pie10} allowing to decompose $2$-dimensional cusp biforms of weight $k$ and $h$ into $2$-dimensional cusp biforms of weight $2$;
\Ee
as it was developed in proposition~1.17.
\vskip 11pt

{\bbf The global program of Langlands over number fields is used in chapter~2 to generate cusp (bi)forms from the Weil (bisemi)group of finite symmetric algebraic extensions\/} characterized by degrees being integers module $N$ of zero-th class, i.e. multiples of quanta.

{\bbf Langlands global correspondences are then associated with the covariant (bi)\-functor\/}
\begin{align*}
 {\rm FLGC} : \qquad & \begin{array}[t]{ccc}
{\rm CABG} & \To & {\rm CBCF}\\
\GL_1(\wt  F_{\o\omega }\times \wt F_{\omega }) & \To &
\Pi (\GL_1(\wt  F_{\o\omega } \times \wt F_{\omega } )\end{array}
\\
 {\rm FLGC}(\Phom_{k\to h}) : \qquad & \begin{array}[t]{ccc}
\Pi_k (\GL_1(\wt  F_{\o\omega }\times \wt F_{\omega } ) & \To &
\Pi_h (\GL_1(\wt  F_{\o\omega }\times \wt F_{\omega } )\end{array}
\end{align*}
which is a  (bi)function assigning:
\Bi
\item to each algebraic bilinear semigroup
$\GL_1(\wt F_{\o\omega }\times \wt F_{\omega } )$ over the product, right by left, of symmetric complex finite algebraic extensions
$\wt F_{\o\omega }$ and $ \wt F_{\omega } $ of the bisemigroup category \cite{Pie3} ${\rm CABG}$;

\item its cuspical representation
$\Pi (\GL_1(\wt  F_{\o\omega }\times \wt F_{\omega } ))$ of the bisemigroup category ${\rm CBCF}$ of complex cuspidal representations
\Ei
in such a way that:
\Bi
\item to each map
\[
\Phom_{k\to h} : \qquad   
\GL_k(\wt  F_{\o\omega }\times \wt F_{\omega } )   \To  
\GL_h(\wt  F_{\o\omega }\times \wt F_{\omega } )\]
sending a bilinear semigroup of dimension $k$ to a bilinear semigroup of dimension $h$, $k>h$, in the sense of what was abovementioned;

\item corresponds a map
$ {\rm FLGC}(\Phom_{k\to h}) $ sending the cuspidal representation\linebreak
$\Pi_k (\GL_1(\wt  F_{\o\omega }\times \wt F_{\omega } )$ of 
$\GL_1(\wt  F_{\o\omega }\times \wt F_{\omega } )$, which is a two-dimensional cusp (bi)form of weight $k$, to the cuspidal representation
$\Pi_h(\GL_1(\wt  F_{\o\omega }\times \wt F_{\omega } ))$ which is a two-dimensional cusp biform of weight $h$.
\Ei

It then results that a (Weil) cusp form ${}^{(\omega )}f_k(z_N)$ of weight $k$ and level $N$, generated from a Langlands global correspondence, can be identified with a (classical) cusp form $f_k(z_N)$ of weight $k$ and level $N$ if the Weil algebraic unitary fundamental domain $D^{(2)}_{F_{\omega _n};N;k}$ of ${}^{(\omega )}f_k(z_N)$ covers the fundamental classical domain
 $D^{(2)}_{nN;k}$ of $f_k(z_N)$.
 
 {\bbf Now, results on the local and global curvatures of two-dimensional tori are presented in proposition~2.9 in the light of the covering of the cusp form
 $f_k(z_N)$ by the global elliptic semimodule
 $\ELLIP_L(2,n,m_n)$\/}: this leads to a new dynamical transition from global euclidean geometry to local hyperbolic and spherical geometries.
 \vskip 11pt
 
 {\bbf In chapter~3, a generalization of two-dimensional cusp forms towards the weak Maass forms, the theta series, the Mock modular forms and the tau function is proposed\/}.
 
 But, first, {\bbf the (Mellin) transform\/} \cite{Bom} (which is a linear continuous map):
 \begin{align*}
 \phi _L: \qquad f_k ( z_{N-k} ) &\To L(f_k,s_{N-k_+})=\sum_{n=1}^\infty c_{nk}\ n^{-s_{N-k_+}}\\[11pt]
\rresp{ \phi _R: \qquad f_k ( z^*_{N-k} ) &\To L(f^*_k,s_{N-k_+})=\sum_{n= 1}^\infty c^*_{nk}\ n^{-s_{N-k_-}}}
\end{align*}
{\bbf from the two-dimensional cusp form of weight $k$ and level $N$, $f_k(z_{N-k})$
\resp{the dual cusp form  $f^*_k(z^*_{N-k})$} to the corresponding $L$-function 
$L(f_k,s_{N-k_+})$
\resp{$L(f^*_k,s_{N-k_-})$} is pointed out to be nonperiodic\/} because of the factors
$n^{-s_{N-k_+}}$
\resp{$n^{-s_{N-k_-}}$} where 
{\bbf $s_{N-k_+}$\/}
\resp{$s_{N-k_-}$} {\bbf is the complex variable conjugate to $z_{N-k}$\/}
\resp{$z^*_{N-k}$} being a complex point of order $(N\times k)$.

{\bbf The nonperiodicity of $\phi _L$ \resp{$\phi _R$} results from the map of the unique period $T$ of $f_k(z_{N-k})$ \resp{$f^*_k(z^*_{N-k})$} to the set
$\L\{\dfrac1{T^{(n)}_{N-k}}\R\}_{n=1}^\infty $ of inverse periods, i.e. unit complex frequencies, associated with the energies of one space quantum on the different levels ``$n$'' which can be evaluated from the consecutive spacings
$\delta \gamma _n=\gamma _{n+1}-\gamma _n$ between the nontrivial zeros of
$\zeta  (s)$ \cite{Pie8}\/}.
\vskip 11pt

In sections~3.3 to 3.8, {\bbf the theta series\/}
\[ \theta _{n^2k/2}(z_{N-k}) =\sum_n d_{nk/2}\ e^{2\pi in^2z_{N-k}}\;, \]
introduced from the quadratic form
\[ Q(n) = \sum_{i=1}^k  \sum_n a^{(i)}_{nn}\ n^2_i\;, \]
{\bbf are proved to be a real analytic modular form of weight $k/2$ and level $N$ in one-to-one correspondence with the two-dimensional modular form
\[
f_k(z_{N-k})=\sum_nc_{nk}\ q^n_{N-k}\]
of weight $k$ and level $N$\/} in such a way that the two-dimensional semitori
$c_{nk}\ q^n_{N-k}$, of which generators are two semicircles at $n$ transcendental quanta, are sent into semicircles $d_{nk/2}\ e^{2\pi in^2x_{N-k/2}}$ at $n^2$ transcendental quanta.

This approach corresponds to the classical one developed by G. Shimura and J.P. Serre who proved that every modular form of weight $1/2$ on $\Gamma _1(N)$ is a linear combination of theta series with characters \cite{Shi}, \cite{Ser}, \cite{S-S}.
\vskip 11pt

{\bbf In sections 3.9 to 3.15, a geometric interpretation of weak Maass forms is proposed in terms of modular curves\/}.  As in the Fourier series expansion
\[ f^{\omega M}_k(z) =\sum_{n=n_0}^\infty \gamma (f,n;y)\ q^{-n}+
\sum_{n=n_1}^\infty  a (f,n)\ q^{n}\;, \qquad q=e^{2\pi iz}\]
of a weak form $f^{\omega M}_k(z)$ of weight $k$, the imaginary dimension ``$y$'' is manifestly lowered with respect to the dimension ``$x$'' in $z=x+iy$, the weak
Maass form $f^{\omega M}_k(z)$ will be decomposed in {\bbf series expansion\/}
\[
f^{\omega M}_k(z_{N-k})=\sum_nT_n^{2,(e\ell)}(S^1_{a_{n_{N-k/2}}},e\ell^1_{ef,n_{N-k/2}})
\]
{\bbf consisting in the sum of two-dimensional semitori $T^{2,e\ell}_n(-)$ with semielliptic cross sections $e\ell^1_{ef,n_{N-k/2}}$\/}
in such a way that:
\Bean
\item the semicircular sections of the two-dimensional semitori $T^2_n(a_{nm|d^2},d)\in f_k(z_{N-k})\approx f_k^{\omega M}(z_{N-k})$ are transformed bijectively into semielliptic sections;

\item {\bbf $f^{\omega M}_k(z_{N-k})$ be periodic, holomorphic and weakly modular (i.e. deviated from circularity)\/};

\item the holomorphic part $\sum\limits_na(f,n)\ q^n$ of the Fourier series of $f^{\omega M}_k(z)$ corresponds to the sum of the products of two orthogonal semicircles of which the one at ``imaginary'' semicircular section is the equation of a semicircle inscribed in the ellipse $e\ell^1_{ef,n_{N-k/2}}$;

\item {\bbf the nonholomorphic part $\sum\limits_n\gamma (f,n;y)\ q^{-n}$ of $f^{\omega M}_k(z)$ corresponds to its shadow\/} as developed in proposition~3.11.
\Ee

{\bbf In section~3.14, generalized weak Maass forms are introduced as elliptic forms $f^{E\ell}_k(z_{N-k})$ of weight $k$ and level $N$\/} having a decomposition into the sum of surfaces of revolution of (semi)ellipses rotating around ellipses instead of circles as for the weak Maass forms.

The commutative diagramm
\[ \begin{psmatrix}[colsep=.5cm,rowsep=1cm]
f_k(z_{N-k}) & &f^{\omega M}_k(z_{N-k})\\
& f^{E\ell}_k(z_{N-k})
\psset{arrows=->,nodesep=3pt}
\ncline{1,1}{1,3}^{\sim}
\ncline{1,3}{2,2}>{\raisebox{-2mm}{\rotatebox{45}{$\sim$}}}
\ncline{2,2}{1,1}<{\rotatebox{-45}{$\sim$}}
\end{psmatrix}
\]
indicates {\bbf the possible transformation of cups forms $f_k(z_{N-k})$ of weight $k$ and level $N$ into the respective weak Maass forms $f^{\omega M}_k(z_{N-k})$ and elliptic forms $f^{E\ell}_k(z_{N-k})$ which are elliptic functions, i.e. doubly periodic meromorphic functions.

In sections~3.16 to 3.20, the Ramanujan theta functions are analyzed in the light of the new geometric interpretation of weak Maass forms\/}.

Each Ramanujan theta function \cite{Ono1} is a Mock theta function \cite{Wat} given by the $q$-series $H(q)=\sum\limits_n a_n\ q^n$ in such a way that $q^\lambda \ H(q)$, $\lambda \in \QQ$, be a Mock modular form of weight $1/2$ whose shadow is a unary thera series of weight $3/2$.  A Mock theta function is thus a Mock modular form \cite{Fol} of the space $\MM_k$ of such forms extending the space $M_k$ of classical modular forms of weight $k$ and characterized by a shadow $g=S[h]$ which is a modular form of weight $(2-k)$.

{\bbf Proposition~3.17 introduces a geometric interpretation of Ramanujan theta functions transformed into Mock theta functions related to weak Maass forms in such a way that\/}:
\Bean
\item {\bbf the nonholomorphic part of the weak Maass form $f^{\omega M}_2(\tau _{1-2})$ of weight $2$ and level $1$ corresponds to the shadow $g^*(\tau )$ of the Mock modular form $\hat h(\tau )$ of weight $1/2$ resulting from the Ramanujan theta function $H(q)$\/} after the three step sequence transformation \cite{Zag} recalled in section~3.16.

\item the space $\MM_1$ of Mock theta functions $\hat h(\tau )$ and $\widehat \MM_2$ of weak Maass forms of weight $2$ are isomorphic.
\Ee

Let 
\[
R(\omega ;q)= \sum_{n=1}^\infty  \sum_{m=-\infty }^\infty N(n,m)\ \omega ^m\ q^n\]
be the partition function specializing the 17 Ramanujan's Mock theta functions $H(q)$ \cite{B-O2}.

As $N(n,m)$ is assumed to denote the number of partitions of $n^{(2)}$ transcendental quanta, {\bbf the Dyson's rank $m$ of a partition of $n$ must be the order of the maximal Galois group associated with the considered transcendental extension minus the number of Galois subgroups\/}: this result is the Galois interpretation of the rank of a partition introduced by Dyson as being its largest part minus the number of its parts.
\vskip 11pt

{\bbf In sections~3.21 to 3.23, the Ramanujan tau function\/}
\[ \sum_{n=1}^\infty \tau (n)\ q^n=q\prod_{n=1}^\infty (1-q^n)^{24}=\Delta (z)\;, \]
which is a modular form of weight $12$ and level $1$, 
{\bbf is associated with the two-dimensional cusp form $f_2(z_{12-1})$ of weight $12$ and level $1$ which is proved to be a universal orthogonal cusp form corresponding throughout Langlands global correspondences to the sum of the cuspidal representations of six bilinear algebraic semigroups\/} generating three two-dimensional embedded toric bismisheaves as well as their orthogonal equivalents according to:

\[ \Delta (z) \To f_2(z_{12-1}) \To \Pi^{(12)}(\GL_6(\wt F_{\o \omega }\times_D\wt F_{\omega }))=
\bigoplus_{i=1}^6 \Pi^{(2_i)}(\GL_{1_i}(\wt F_{\o \omega }\times_D\wt F_{\omega }))\]
where:
\Bi
\item $\Pi^{(12)}(\GL_6(\wt F_{\o \omega }\times_D\wt F_{\omega }))$ is the twelve-dimensional cuspidal representation of\linebreak 
$\GL_6(\wt F_{\o \omega }\times_D\wt F_{\omega })$;

\item $\Pi^{(2_i)}(\GL_{1_i}(\wt F_{\o \omega }\times_D\wt F_{\omega }))$ is the two-dimensional cuspidal representation of\linebreak
$\GL_{1_i}(\wt F_{\o \omega }\times_D\wt F_{\omega })$.
\Ei

This can be finally generalized to {\bbf universal nonorthogonal cuspidal representations including bilinear crossed cuspidal representations of interaction in such a way that $\Pi^{(12)}(\GL_6(\wt F_{\o \omega }\times_D\wt F_{\omega }))$ decomposes then nonorthogonally by means of the Langlands bilinear functoriality conjecture\/} according to:
\begin{align*}
\Pi^{(12)}_{re\ell}(\GL_6(\wt F_{\o \omega }\times_D\wt F_{\omega }))
= &
\bigoplus_{i=1}^6 \Pi^{(2_i)}(\GL_{1_i}(\wt F_{\o \omega }\times_D\wt F_{\omega }))\\
& \qquad \bigoplus_{i\neq j=1}^6 \Pi^{(2_i)}(\GL_{1_i}(\wt F_{\o \omega }))\otimes_{OD}\Pi^{(2_j)}(\GL_{1_j}(\wt F_{\omega }))
\end{align*}
where the second sum on the right hand side refers to the six relevant offdiagonal crossed cuspidal representations of interaction as developed in section~3.22.

This universal nonorthogonal cuspidal representation
$\Pi^{(12)}_{re\ell}(\GL_6(\wt F_{\o \omega }\times_D\wt F_{\omega }))$ is then\linebreak mapped injectively:
\[
\Ms_{\Pi ^{(12)}\to\mathop{\u \times}\limits^6}: \qquad
\Pi^{(12)}_{re\ell}(\GL_6(\wt F_{\o \omega }\times_D\wt F_{\omega }))\To
\mathop{\times}\limits_{\o{i=1}}^6 (f_{2_i}(z^*_{2-1})\times f_{2_i}(z_{2-1}))
\]
into the cross binary product
\[
\mathop{\times}\limits_{\o{i=1}}^6 (f_{2_i}(z^*_{2-1})\times f_{2_i}(z_{2-1}))
=\L( \sum_{i=1}^6 f_{2_i}(z^*_{2-1})\R)\times\L( \sum_{i=1}^6 f_{2_i}(z_{2-1})\R)
\]
between the six cusp biforms
$f_{2_i}(z^*_{2-1})\times f_{2_i}(z_{2-1})$ of dimension $2$, weight $2$ and level $1$.

\section{Cusp (bi)forms of weight $k$}

\subsection{Classical definitions of the modular forms}

Let $H$ be the Poincare upper half plane of complex numbers $z=x+iy$ with strictly positive imaginary parts $y$.  

Let $S_\omega $ denote the set of pairs $\omega =(\omega _1,\omega _2)$ of complex numbers $\omega _i=x_i+iy_i$, $y_i>0$, $i=1,2$.

Let $z\in H$ be given by $z=\dfrac{\omega _1}{\omega _2}$ in such a way that
$\omega =(\omega _1,\omega _2)$ be sent into $z$ by the $\cit$-linear map $g_{\ptcit}:\cit^2\to\cit$ \cite{Del}.

{\bbf The set $G_\zit^{(2)}(H)$ of lattices of $H$ is the quotient $H\big/ \GL(2,\zit)$ of $H$ by the group $\GL(2,\zit)$ \/} in such a way that $G_\zit^{(2)}(H)$ corresponds to $\Hom(\zit^2,\cit)$.
\vskip 11pt

Let $F(\omega )=\omega _2^{-k}f(\omega _1/\omega _2)$ define a holomorpic and homogeneous elliptic function of weight $k$ in the upper half plane \cite{God}, $k$ being fairly often an even integer.

{\bbf $F(\omega )$ is a modular (or ``automorphic'') form of weight $k$ if it is invariant under the substitution\/} $\omega \to g_{SL_2}\omega$ where $g_{SL_2}$ is the matrix $\Bm a & b \\ c & d\Em$, $a,b,c,d\in\zit$, {\bbf of the homogeneous modular group $SL(2,\zit)$\/} of homographic  transformations, verifying $ad-bc=1$. The modular form $F(\omega )$ of weight $k$ is then equivalent to a $\cit$-valued function $f(z)$ of moderate growth on $H$ verifying:
\[
f(z)=(cz+d)^{-k}\ f\L(\dfrac{az+b}{cz+d}\R)\;, \qquad \forall\ \Bm a & b\\ c & d\Em\in SL(2,\zit)\;,
\]
where $(cz+d)^{-k}$ is the weight factor.

{\bbf The action of $PSL(2,\zit)=SL(2,\zit)\big/\pm I$ on $z\in\cit$ then generates on $f(z)$\/}:
\Bean
\item {\bbf a ``periodic'' translation by $S:z\to z+1$\/};

\item {\bbf a ``periodic'' space inversion by $T:z\to-1/z$\/}.
\Ee

The modular form of $f(z)$ of  weight $k$ is of level $N$ if it is invariant under the congruence (sub)group $\Gamma _0(N)$ \resp{$\Gamma _1(N)$}:
\begin{align*}
\Gamma _0(N) &= \L\{ \Bm a & b \\ c & d\Em \in SL(2,\zit): c\equiv 0 \text{ (mod }N)\R\}\\[11pt]
\rresp{\Gamma _1(N) &= \L\{ \Bm a & b \\ c & d\Em \in SL(2,\zit): c\equiv 0 \text{ (mod }N), a,d= 1 \text{ (mod }N), \R\}}.
\end{align*}
\vskip 11pt

\subsection{$\rit^2/\zit^2$-lattices of the Poincare upper half plane}

Let $G^{(2)}_\zit(H)=H\big/\GL(2,\zit)$, or more precisely, $G^{(2)}_\zit(H)=\GL_2(\rit)\big/\GL_2(\zit)$ 
{\bbf be the set of $\rit^2\big/\zit^2$-lattices of the Poincare upper half plane $H$\/} generated by the group $\GL_2(\rit)$ isomorphic to the group $\GL(H)$ of automorphisms of $H$ \cite{Bor}.

Similarly, let $G^{(2)}_{\zit_N}(H)=\GL_2(\rit)\big/\GL_2(\zit/N \zit)$ be the set of $\rit^2\big/(\zit/N \zit)^2$-lattices of $H$.

As the subgroup $\GL_2(\zit)$ 
\resp{the congruence subgroup $\GL_2(\zit/N \zit)$} leads to the generation of the
set $\{\Lambda ^{(2)}_\zit[n]\}$
\resp{$\{\Lambda ^{(2)}_{\zit_N}[n]\}$} of $\zit^2$-
\resp{$(\zit/N \zit)^2$-} lattices in $H$, {\bbf the quotient group
$G^{(2)}_\zit)(H)$
\resp{$G^{(2)}_{\zit_N})(H)$} is given by the coset representatives\/}:
\[ \L\{ g^{(2)}_\zit[n]\R\}^t_{n=1}
\rresp{ \L\{ g^{(2)}_{\zit_N}[n]\R\}^t_{n=1} }, \qquad t\le \infty \;.\]
{\bbf The order of the $n$-th coset representative $g^{(2)}_\zit[n]$
\resp{$g^{(2)}_{\zit_N}[n]$} is the integer $n$ \resp{$n\ N$}, corresponding to the number of its elements or automorphisms\/}

On the other hand, as $\Lambda ^{(2)}_\zit[n]$ is included into
$\Lambda ^{(2)}_{\zit_N}[n]$, the quotient subgroup\linebreak
$\Lambda ^{(2)}_{\zit_N}[n] \big/\Lambda ^{(2)}_{\zit}[n] $ is finite and its order is $N$: 
$\Lambda ^{(2)}_{\zit_N}[n] $ is thus of index $N$ in
$\Lambda ^{(2)}_{\zit}[n] $ \cite{Ser}.

Let $H_n$ \resp{$H_{nN}$} denote the Poincare upper half plane restricted to the $n$-th coset representative
$g^{(2)}_\zit[n]$
\resp{$g^{(2)}_{\zit_N}[n]$}, then $H_n\big/\Lambda ^{(2)}_\zit[n]$ is included into 
$H_{nN}\big/\Lambda ^{(2)}_{\zit_N}[n]$.

For every integer $n$ labeling a coset representative $g^{(2)}_\zit[n]$
\resp{$g^{(2)}_{\zit_N}[n]$} of $G^{(2)}_\zit(H)$
\resp{$G^{(2)}_{\zit_N}(H)$}, we have the monomorphism
\[ h_{\Lambda \to g}:\quad \Lambda ^{(2)}_\zit[n]\to g^{(2)}_\zit[n]
\rresp{h_{\Lambda_N \to g_N}:\quad \Lambda ^{(2)}_{\zit_N}[n]\to g^{(2)}_{\zit_N}[n]}, \quad 1\le n\le t\le \infty \;.\]
\vskip 11pt

\subsection{Fundamental domains}

{\bbf A fundamental domain $D^{(2)}_n$ 
\resp{$D^{(2)}_{nN}$} of $H_n$ with respect to the sublattice
$\Lambda ^{(2)}_{\zit}[n]$
\resp{$\Lambda ^{(2)}_{\zit_N}[n]$} is a unitary measurable set (i.e. a period) of $H_n$\/} in such a way that its translations by the vectors of 
$\Lambda ^{(2)}_{\zit}[n]$
\resp{$\Lambda ^{(2)}_{\zit_N}[n]$}  are a partition of $H_n$.

{\bbf As a modular form is periodic, all the fundamental domains
$D^{(2)}_n$ 
\resp{$D^{(2)}_{nN}$}, $\forall\ n$, are equal\/}:
\begin{align*}
 D^{(2)}_1=D^{(2)}_2= \dots &= D^{(2)}_n= \dots = D^{(2)}_t\\[11pt]
\rresp{D^{(2)}_{1N}=D^{(2)}_{2N}= \dots &= D^{(2)}_{nN}= \dots = D^{(2)}_{tN})}.
\end{align*}

{\bbf The fundamental domain $D^{(2)}_n$
\resp{$D^{(2)}_{nN}$} corresponds then to the isotropy group
$I^{(2)}_\zit=\{i_\zit\in I^{(2)}_\zit \mid i_\zit\ z \ i_\zit^{-1}=z\}$
\resp{$I^{(2)}_{\zit_N}=\{i_{\zit_N}\in I^{(2)}_{\zit_N} \mid i_{\zit_N}\ z \ i_{\zit_N}^{-1}=z\}$} of $z$\/}.
\vskip 11pt

\subsection{Proposition}

{\em
The isotropy subgroup
$I^{(2)}_\zit$
\resp{$I^{(2)}_{\zit_N}$} acts on $\GL_2(\rit)$ by conjugation in such a way that the coset representatives 
$\{g^{(2)}_\zit[n]\}_n$
\resp{$\{g^{(2)}_{\zit_N}[n[\}_n$} 
of $G^{(2)}_\zit(H)$
\resp{$G^{(2)}_{\zit_N}(H)$} are the conjugacy classes of $\GL_2(\rit)$.
}
\vskip 11pt

\bpr
The isotropy subgroup $I^{(2)}_\zit$
\resp{$I^{(2)}_{\zit_N}$} has order $1$ \resp{order $N$} according to section~1.2.

As the isotropy subgroup $I^{(2)}_\zit$
\resp{$I^{(2)}_{\zit_N}$} is defined with respect to the set
$\{\Lambda ^{(2)}_\zit[n]\}^t_{n=1}$
\resp{$\{\Lambda ^{(2)}_{\zit_N}[n]\}^t_{n=1}$} of sublattices of $\GL_2(\rit)$ which is isomorphic to the group of automorphisms of $H$, it is clear that there are $t$ classes of automorphisms in $\GL_2(\rit)$ generated by the translation vectors of the sublattices 
$\{\Lambda ^{(2)}_\zit[n]\}_{n}$
\resp{$\{\Lambda ^{(2)}_{\zit_N}[n]\}_{n}$} on the fundamental domain
$D^{(2)}_n$ 
\resp{$D^{(2)}_{nN}$}.  As every automorphism of $H$ is induced by a conjugation of $\GL_2(\rit)$, there are ``$t$'' classes of conjugation of $\GL_2(\rit)$ given by their orders which are integers 
$1\le n\le t\le \infty $
\resp{$N\le n\cdot N\le t\cdot N\le \infty $}.

And, thus, the cosets of $G^{(2)}_\zit(H)$
\resp{$G^{(2)}_{\zit_N}(H)$} are in one-to-one correspondence with the conjugacy classes of $\GL_2(\rit)$.\epr
\vskip 11pt

\subsection{Corollary}
{\em
The action of $\GL_2(\rit)$ on the fundamental domain $D^{(2)}_n$
\resp{$D^{(2)}_{nN}$} induces a homomorphism \cite{Hun}
\[ h_{G\to P_{n}}: \quad \GL_2(\rit) \To P(D^{(2)}_n)
\rresp{h_{G\to P_{nN}}: \quad \GL_2(\rit) \To P(D^{(2)}_{nN})}\]
into the group $P(D^{(2)}_n)$
\resp{$P(D^{(2)}_{nN})$} of all permutations of $D^{(2)}_n$
\resp{$D^{(2)}_{nN}$}.
}
\vskip 11pt

\bpr The map $h_{G\to P_{n}}$
\resp{$h_{G\to P_{nN}}$} being a bijection implies that 
$P(D^{(2)}_n)$
\resp{$P(D^{(2)}_{nN})$} is a group of permutations.\epr
\vskip 11pt

\subsection{Actions of Hecke operators}

{\bbf The sum of the functions on the cosets of $T^{(2)}_\zit(H)=T_2(\rit)\big/T_2(\zit)$
\resp{$T^{(2)}_{\zit_N}(H) =T_2(\rit)\big/T_2(\zit/N \zit)$}\/}
($T_2(\rit)$ being the subgroup of upper triangular matrices of $\GL_2(\rit)$)
{\bbf can be reached by the action of the Hecke operator $T_k(n)$
\resp{$T_k(n;N)$} on the modular form $f_k(z)$ \resp{$f_k(z)_N$}\/} of weight $k$ and level 1 \resp{$N$} given by:
\begin{align*}
T_k(n)\ f_k(z) &= n^{k-1}\ T_k(n)\ F_k(G^{(2)}_\zit(H))\\[11pt]
&=n^{k-1}\sum_{\substack{ad=n\\ 0\le b< d}}d^{-k}\ f_k\L(\dfrac{az+b}d\R)\\[15pt]
\rresp{T_k(n;N)\ f_k(z)_N &= (n\ N)^{k-1}\ T_k(n;N)\ F_k(G^{(2)}_{\zit_N}(H))\\[11pt]
&=(n\ N)^{k-1}\sum_{ad=nN}d_N^{-k}\ f_k\L(\dfrac{az+b}d\R)_N}\end{align*}
according to a complete set of upper triangular representatives
$\Bm a & b\\ 0 & d\Em$
\resp{$\Bm a & b\\ 0 & d\Em_N$}, $ad=0$ modulo $N$, $b=*$ modulo $N$ of the group
$M^n$ \resp{$\Gamma _1(N)$, $\Gamma _0(N),\dots$)} of integral matrices with determinant $n$ \resp{$n\ N\dots$} with respect to $SL_2(\zit)$ \resp{$SL_2(\zit/N \zit)$}.
\vskip 11pt

\subsection{Proposition}

{\em
The sum $\sum\limits_{ad=n} d^{-k}\ f_k\L(\dfrac{az+b}d\R)$
\resp{$\sum\limits_{ad=nN} d_N^{-k}\ f_k\L(\dfrac{az+b}d\R)_N$}
of the functions on the coset representatives $t^{(2)}_\zit)[n]$
\resp{$t^{(2)}_{\zit_N})[n]$} of 
$T_2(\rit)/T_2(\zit)$  
\resp{$T_2(\rit)/T_2(\zit/N \zit)$} is in one-to-one correspondence with the action
\begin{align*}
T_k(n) : \qquad T^{(2)}_\zit (H)&\To \sum_{(T^{(2)}_\zit(H):T^{(2)}_\zit(H))=n} t^{(2)}_\zit(H)\\[11pt]
\rresp{T_k(n;N) : \qquad T^{(2)}_{\zit_N} (H)&\To \sum_{(T^{(2)}_{\zit_N}(H):T^{(2)}_{\zit_N}(H))=nN} t^{(2)}_{\zit_N}(H)}
\end{align*}
of the Hecke operator $T_k(n)$ \resp{$T_k(n;N)$} associating (by a correspondence) with the lattice
$T^{(2)}_\zit(H)$
\resp{$T^{(2)}_{\zit_N}(H)$} the sum of its sublattices
$t^{(2)}_\zit(H)$
\resp{$t^{(2)}_{\zit_N}(H)$} of index $n$ \resp{$n\ N$} in $T^{(2)}_\zit(H)$
\resp{$T^{(2)}_{\zit_N}(H)$}.
}
\vskip 11pt

\bpr
This is immediate if we take into account the monomorphism
\[ h_{\Lambda \to g} : \quad \Lambda ^{(2)}_\zit[n] \To t^{(2)}_\zit[n]
\rresp{h_{\Lambda_N \to g_N} : \quad \Lambda ^{(2)}_{\zit_N}[n] \To t^{(2)}_{\zit_N}[n]}
\]
of section~1.2 mapping the $n$-th $\zit^2$-sublattice $\Lambda ^{(2)}_\zit[n]$
\resp{$(\zit/N \zit)^2$-sublattice\linebreak   $\Lambda ^{(2)}_{\zit_N}[n]$}
into the $n$-th coset representative $t^{(2)}_\zit[n]$
\resp{$t^{(2)}_{\zit_N}[n]$} of 
$T_2(\rit)/T_2(\zit)$
\resp{$T_2(\rit)/ T_2(\zit/N \zit)$}.

Then, the functions on the cosets $t^{(2)}_\zit[n]$
\resp{$t^{(2)}_{\zit_N}[n]$} are in one-to-one correspondence with the functions
$f_k\L(\dfrac{az+b}d\R)$
\resp{$f_k\L(\dfrac{az+b}d\R)_N$} on the sublattices
$\Lambda ^{(2)}_{\zit}[n]$
\resp{$\Lambda ^{(2)}_{\zit_N}[n]$} described by the coset representatives
$\Bm a & b\\ 0 & d\Em$
\resp{$\Bm a & b\\ 0 & d\Em_N$} of $M^n$
\resp{$\Gamma _1(N),\Gamma _0(N),\dots$}.\epr
\vskip 11pt

\subsection{Fourier series development of cusp forms}

Let
\begin{align*}
T_k(n) : \qquad F_k(T^{(2)}_\zit(H)) &\To \sum_{ad=n}d^{-k}\ f_k\L(\dfrac{az+b}d\R)\\[11pt]
\rresp{T_k(n;N) : \qquad F_k(T^{(2)}_{\zit_N}(H)) &\To \sum_{ad=nN}d_N^{-k}\ f_k\L(\dfrac{az+b}d\R)_N}
\end{align*}
be the action of the Hecke operator
$T_k(n)$ \resp{$T_k(n;N)$} on the function
$F_k(T^{(2)}_\zit(H))$
\resp{$F_k(T^{(2)}_{\zit_N}(H))$} over the two-dimensional lattice
$T^{(2)}_{\zit}(H)$
\resp{$T^{(2)}_{\zit_N}(H))$} given by the Laurent series:
\[
F_k(T^{(2)}_\zit(H))=\sum_{n=0}^\infty a'_n\ h^n
\rresp{F_k(T^{(2)}_{\zit_N}(H))=\sum_{n=0}^\infty b'_n\ k^n}
\]
with $h=H_k(z)$
\resp{with $k=K_k(z_N)$} being a function of $z$ \resp{$z_N$} {\bbf which is a complex point of order $(1\times k)$
\resp{$(N\times k)$} to which a period ``$1\times k$''
\resp{$N\times k$} corresponds\/}.

Let 
\[
\tau ^{(k)}_{h\to q} : \quad h=H_k(z)\To q=e^{2\pi iz}
\rresp{\tau ^{(k)}_{k\to q_N} : \quad k=K_k(z_N)\To q_N=e^{2\pi iz_N}}
\]
be {\bbf the toroidal mapping sending $h$ \resp{$k$} into $q$ \resp{$q_N$}
in such a way that the Laurent series
$F_k(T^{(2)}_\zit(H))$
\resp{$F_k(T^{(2)}_{\zit_N}(H))$} be transformed into the Fourier series\/}:
\[ f_k(z) = \sum_{n=1}^\infty a_n\ q^n
\rresp{f_k(z_N) = \sum_{n=1}^\infty b_n\ q_N^n}.\]

Taking into account the main property
\[
T_k(m)\ T_k(n)=\sum_{d|m;n} T_k(mn/d^2)\ d^{k-1}\]
of the Hecke operators, we find that \cite{Lan}
\[
T_k(m)\ f_k(z)=\sum_{n=1}^\infty a_{nk}(m)\ q^n
\rresp{T_k(m;N)_\varepsilon \ f_k(z_N)=\sum_{n=1}^\infty c_{nk}\ q_N^n}
\]
with $\varepsilon :(\zit/N \zit)^*\to \cit^*$ a Dirichlet character mod $N$, where
\[
a_{nk}(m) = \sum_{d|m,n}d^{k-1}\ a_{nm|d^2}
\rresp{c_{nk} = \sum_{d|(m,n)}\varepsilon (d)\ d^{k-1}\ a_{nm|d^2}}.\]
If $f_k(z)$ 
\resp{$f_k(z_N)$} is an eigenfunction of all the Hecke operators 
$T_k(m)$ \resp{$T_k(m;N)$}, $1\le m\le\infty $, according to:
\[
T_k(n)\ f_k(z) = \lambda (n)\ f_k(z)
\rresp{T_k(m;N)_\varepsilon \ f_k(z_N) = \lambda (n)_N\ f_k(z_N)},\]
we have that
$a_{nk}=\lambda (n)$ for $f_k(z)$ normalized with $a_1=1$ \resp{$c_{nk}=\lambda (n)_N=\sum_{d|n}\varepsilon (d)\ d^{k-1}$}
where
{\bbf $\lambda (n)$ \resp{$\lambda (n)_N$} is a Hecke eigencharacter of $f_k(z)$ \resp{$f_k(z_N)$}\/}.

If $a_1=0$ \resp{$c_1=0$}, then $T_k(n)\ f_k(z)$
\resp{$T_k(n;N)\ f_k(z_N)$} is a cusp form of weight $k$ and level $1$ \resp{$N$}.
\vskip 11pt

\subsection{Proposition (Geometric interpretation of cusp forms)}

{\em Let
\[
T_2(n)\ f_2(z) = \sum_{n=1}^na_n(m)\ q^n
\rresp{T_2(n;N)\ f_2(z_N) = \sum_{n=1}^nc_n\ q^n_N}
\] be a weight $2$ ($k=2$) cusp form of level $1$ \resp{$N$}.

Then, we have that {\bbf the $n$-th term $a_n(m)\ q^n$ \resp{$c_n\ q^n_N$} given by
\begin{align*}
& \begin{aligned}[t]
&T_n^2(a_{nm|d^2},d) \\
& \qquad \simeq a_n(m)\ q^n 
=  r_{S^1_{a_n}}\ e^{2\pi inx}\times r_{S^1_{d_n}}\ e^{2\pi in(iy)}
\end{aligned}\\[11pt]
\text{(resp.} \quad & \begin{aligned}[t]
&T_{nN}^2(a_{nm|d^2},d\varepsilon (d)) \\
& \qquad \simeq c_n\ q^n_N 
 = r_{S^1_{a_{nN}}}\ e^{2\pi inx_N}\times r_{S^1_{d_{nN}}}\ e^{2\pi in(iy_N)})\end{aligned}
\end{align*}
is a two-dimensional semitorus
$T_n^2(a_{nm|d^2},d)$
\resp{$T_{nN}^2(a_{nm|d^2},d\varepsilon (d))$}\/}:
\Bean
\item {\bbf generated by two orthogonal semicircles
$S^1_{a_n}$ 
\resp{$S^1_{a_{nN}}$} and $S^1_{d_n}$ 
\resp{$S^1_{d_{nN}}$} respectively at $a_n$ and $d_n$   (transcendental) quanta (see proposition~2.8)\/} in such a way that the Hecke eigencharacter $\lambda (n)=a_{nk}$ \resp{$\lambda (n)_N=c_{nk}$} decomposes according to
$\lambda (n)=\lambda _a(n)\times \lambda _d(n)$
\resp{$\lambda (n)_N=\lambda _a(n)_N\times \lambda _d(n)_N$}
where $\lambda _a(n)$
\resp{$\lambda _a(n)_N$} refers to the radius $r_{S^1_{a_n}}$ \resp{$r_{S^1_{a_{nN}}}$} and
$\lambda _d(n)$
\resp{$\lambda _d(n)_N$} refers to the radius $r_{S^1_{d_n}}$ \resp{$r_{S^1_{d_{nN}}}$};

\item {\bbf whose area results from $n$ \resp{$nN$} permutations of the fundamental domain
$D^{(2)}_n$
\resp{$D^{(2)}_{nN}$} of $T_2(n)\ f_2(z)$
\resp{$T_2(n;N)\ f_2(z_N)$}\/}.
\Ee
}
\pagebreak

\bpr
\Bena
\item The $n$-th term $f^{(n)}_2(z)=a_n(m)\ q^n$
\resp{$f^{(n)}_2(z_N)$}
of the cusp form
$T_2(n)\ f_2(z)$
\resp{$T_2(n;N)\ f_2(z_N)$}corresponds to the $n$-th cuspidal sublattice on products
$(\omega '_1\times\omega '_2)$
\resp{$(\omega '_{1_N}\times\omega '_{2_N})$} of complex numbers with positive imaginary parts and results from the map:
\[
z=\tfrac{\omega _1}{\omega _2}\To f^{(n)}_2(z) = \omega '_1\times\omega '_2
\rresp{z_N=\tfrac{\omega _{1N}}{\omega _{2N}}\To f^{(n)}_2(z_N) = \omega '_{1N}\times\omega '_{2N}}
\]
in such a way that $\omega '_1$ and $\omega '_2$ \resp{$\omega '_{1N}$ and $\omega '_{2N}$}
are respectively the orthogonal semicircles
\begin{align*}
S^1_{a_n} &= r_{S^1_{a_n}}\ e^{2\pi inx}
&\rresp{S^1_{a_{nN}} &= r_{S^1_{a_{nN}}}\ e^{2\pi inx_n}}, && x>0\;, \\
\text{and}\quad
S^1_{d_n} &= r_{S^1_{d_n}}\ e^{2\pi in(iy)}
&\rresp{S^1_{d_{nN}} &= r_{S^1_{d_{nN}}}\ e^{2\pi in(iy_n)}}, && y>0\;, 
\end{align*}
of radii $r_{S^1_{a_{n}}}$ \resp{$r_{S^1_{a_{nN}}}$}
and $r_{S^1_{d_{n}}}$ \resp{$r_{S^1_{d_{nN}}}$}.

Indeed, as it was proved in \cite{Pie1} and \cite{Pie2},
{\bbf
\[ f^{(1)}_2(z) \simeq q=e^{2\pi iz}=e^{2\pi ix}\cdot e^{2\pi i(iy)}\]
corresponds to the product of two orthogonal unitary (semi)circles because
$e^{2\pi i(iy)}\neq e^{-2\pi y}$ where $(iy)$ is an ``imaginary'' angle which cannot be ``assimilated'' to the function
$\exp(i\theta )=\cos\theta +i\sin\theta $\/}.  And, thus, this  ``imaginary'' angle is the angle of a circle orthogonal to $e^{2\pi ix}$.

\item As the $n$-th cusp subform $f^{(n)}_2(z)$ 
\resp{$f^{(n)}_2(z_N)$} is a function on the $n$-th coset representative 
$t^{(2)}_\zit[n]$
\resp{$t^{(2)}_{\zit_N}[n]$} of 
$T^{(2)}_{\zit}(H)=T_2(\rit)/T_2(\zit)$
\resp{$T^{(2)}_{\zit_N}(H)=T_2(\rit)/T_2(\zit/N \zit)$}, the action of
$t^{(2)}_{\zit}[n]$
\resp{$t^{(2)}_{\zit_N}[n]$} on the fundamental domain
$D^{(2)}_n$
\resp{$D^{(2)}_{nN}$} induces the map:
\begin{align*}
h_{f^{(n)}_2(z)\to P_n} :   && f^{(n)}_2(z) \quad &\To \quad P_n(f^{(1)}_2(z))\\
\rresp{h_{f^{(n)}_2(z_N)\to P_{nN}} :   && f^{(n)}_2(z_N) \quad &\To   \quad P_{nN}(f^{(1)}_2(z_N))}
\end{align*}
of $f^{(n)}_2(z)$
\resp{$f^{(n)}_2(z_N)$} into the subgroup 
$P_n(f^{(1)}_2(z))$
\resp{$P_{nN}(f^{(1)}_2(z_N))$} of $n$ \resp{$nN$} permutations of the unitary cusp subform
$f^{(1)}_2(z)=a_1\ q^1$
\resp{$f^{(1)}_2(z_N)=c_1\ q^1_N$} on $D^{(2)}_n$
\resp{$D^{(2)}_{nN}$}
according to corollary~1.5.\epr
\Ee
\vskip 11pt

\subsection{Proposition}

{\em
{\bbf The weight $2$ cusp form of level $1$ \resp{$N$}
\[
T_2(n)\ f_2(z) = \sum_{n=1}^\infty a_n(m)\ q^n
\rresp{T_2(n;N)\ f_2(z_N) = \sum_{n=1}^\infty c_n\ q^n_N}
\]
is a function
\Bean
\item composed of the sum of increasing embedded two-dimensional semitori localized in the upper half plane:
\[
a_1\ q^1\subset a_2\ q^2\subset \dots \subset a_n\ q^n
\rresp{c_1\ q^1_N\subset c_2\ q^2_N\subset \dots \subset c_n\ q^n_N};
\]

\item which is periodic in the sense that each cusp subform $f^{(n)}_2(z)$
\resp{\linebreak $f^{(n)}_2(z_N)$} has an area being $n$ \resp{$nN$} times the area of the unitary cusp subform
$f^{(1)}_2(z)=a_1\ q^1$
\resp{$f^{(1)}_2(z_N)=c_1\ q^1_N$}.
\Ee
}}
\vskip 11pt

\bpr This results immediately from proposition~1.9.\epr
\vskip 11pt

\subsection{Cosemialgebra of dual cusp forms}

Let $S_L(k,1)$
\resp{$S_L(k,N)$} denote the (semi)algebra of cusp forms
$f_k(z)$
\resp{\linebreak $f_k(z_N)$} of weight $k$ and level $1$ \resp{$N$} which, expanded in Fourier series
\[
f_k(z)=\sum_{n=1}^\infty a_{nk}(m)\ q^n
\rresp{f_k(z_N)=\sum_{n=1}^\infty c_{nk}(m)\ q^n_N}\]
are eigenfunctions of Hecke operators $T_k(m)$ 
\resp{$T_k(m;N)$} and are holomorphic in the Poincare upper half plane $H$.

{\bbf The dual (semi)algebra of $S_L(k,1)$
\resp{$S_L(k,N)$} is the cosemialgebra $S_R(k,1)$
\resp{$S_R(k,N)$} of dual cusp forms $f_k(z^*)$
\resp{$f_k(z^*_N)$}\/} of weight $k$ and level $1$ \resp{$N$} which:
\Bi
\item are holomorphic in the Poincare lower half plane $H^*$ of complex numbers $z^*=x-iy$ with negative imaginary parts $-y$;

\item are eigenfunctions of Hecke operators $T^*_k(m)$
\resp{$T^*_k(m;N)$};

\item are {\bbf expanded in Fourier series
\begin{align*}
f_k(z^*) &= \sum_{n=1}^\infty a^*_{nk}(m)\ q^{*n}\\[15pt]
\rresp{f_k(z^*_N) &= \sum_{n=1}^\infty c^*_{nk}(m)\ q^{*n}_N}, \qquad q^*=e^{2\pi iz^*}\;, \end{align*}
according to the complete set of lower triangular coset representatives\/}
$\Bm a&0\\ b&d\Em$ 
\resp{$\Bm a&0\\ b&d\Em_N$ $ad=0\mod N$, $b=*\mod N$} of the group
$M^n$ 
\resp{the congruence subgroup $\Gamma _1(N),\Gamma _0(N),\dots$} of integral matrices with determinant $n$ \resp{$nN$}
with respect to $SL_2(\zit)$
\resp{$SL_2(\zit/N \zit)$}.
\Ei

We thus have that:
\Bean
\item the $n$-th dual cusp subform $f_k^{(n)}(z^*)$
\resp{$f_k^{(n)}(z^*_N)$} of weight $k$ with respect to
$f_k^{(n)}(z)$
\resp{$f_k^{(n)}(z_N)$} is a cofunction on the $n$-th coset representative
$t^{*(2)}_{k,\zit}[n]$
\resp{$t^{*(2)}_{k,\zit_N}[n]$} of
\[
T^{(2)}_{k,\zit}(H^*)=T^t_2(\rit)/T^t_2(\zit)
\rresp{ T^{2)}_{k,\zit_N}(H^*)= \linebreak T^t_2(\rit)/T^t_2(\zit\big/N \zit) }\]
where:
\Bi
\item $T^{(2)}_{k,\zit}(H^*)$
\resp{$T^{(2)}_{k,\zit_N}(H^*)$} is the set of $\rit^2/\zit^2$
\resp{$\rit^2/(\zit/N \zit)^2$}-lattices in the Poincare lower half plane $H^*$;

\item $T_2(\rit)$
\resp{$T_2^t(\rit)$} is the two-dimensional subgroup of upper \resp{lower} triangular matrices.
\Ei

\item $D^{(k)}_n$
\resp{$D^{(k)}_{nN}$} is a fundamental domain of  $H^*_n$ with respect to the $n$-th 
$\zit^2$-sublattice$\Lambda ^{*(2)}_{k,\zit}[n]$
\resp{$(\zit/N \zit)^2$-sublattice$\Lambda ^{*(2)}_{k,\zit_N}[n]$}.
\Ee
\vskip 11pt

From now on, we shall only consider the semialgebra
$S_L(k,N)$ \resp{the cosemialgebra $S_R(k,N)$} of weight $k$ cusp forms $f_k(z_N)$ \resp{dual cusp forms $f_k(z^*_N)$} of level $N$, the level $N=1$ being a particular case of it.
\vskip 11pt

\subsection{Bisemialgebra of cusp biforms}

{\bbf Let
\[ R_{f_{k,N}}=\{f^{(n)}_k(z_N)\}_n
\rresp{R_{f^*_{k,N}}=\{f^{(n)}_k(z^*_N)\}_n}
\]
denote the \lr semiring of cusp subforms $f^{(n)}_k(z_N)$
\resp{dual cusp subforms $f^{(n)}_k(z^*_N)$}\/}
as introduced in section~1.11.

Then, $S_L(k,N)$
\resp{$S_R(k,N)$} is a left $R_{f_{k,N}}$-semialgebra
\resp{right $R_{f^*_{k,N}}$-cosemi\- algebra}.

Let $R_{f^*_{k,N}}\times R_{f_{k,N}}$ be the associated bisemiring \cite{Pie3}.

Then, {\bbf the $R_{f^*_{k,N}}\times R_{f_{k,N}}$-bisemialgebra
$S\RL(k,N)$ of cusp biforms $f_k(z^*_N)\times f_k(z_N)$ of weight $k$ and level $N$ is a bisemiring such that\/}:
\Bean
\item $S\RL(k,N)$ is a unitary $R_{f^*_{k,N}}\times R_{f_{k,N}}$-bisemimodule:
\item $\mu \RL:S\RL(k,N)\u\times S\RL(k,N)\to S\RL(k,N)$ is a bilinear homomorphism where {\bbf $\u\times$ is the cross binary operation acting on the cusp biforms\/} as follows:
\begin{multline*}
\L[ \L({}^1f_k(z^*_N)\times{}^1f_k(z_N) \R)\u\times
 \L({}^2f_k(z^*_N)\times{}^2f_k(z_N) \R)
 \R]\\[11pt]
 =
 \L[ \L({}^1f_k(z^*_N)+{}^2f_k(z^*_N) \R)\times
 \L({}^1f_k(z_N)+{}^2f_k(z_N) \R)
 \R]
 \end{multline*}
 {\bbf allowing cross products
$\L({}^1f_k(z^*_N)\times{}^2f_k(z_N) \R)$ and
$ \L({}^2f_k(z^*_N)\times{}^1f_k(z_N) \R)
$ between left and right cusp forms $1$ and $2$\/}.

\item $\eta\RL:R_{f^*_{k,N}}\times R_{f_{k,N}}\to S\RL(k,N)$ is an injective homomorphism \cite{Pie1}.
\Ee

The existence of the bisemialgebra $S\RL(k,N)$ of cusp biforms is especially important due to the homomorphism:
\[ \psi :\quad S\RL(k,N) \To \End(S_L(k,N))\]
allowing to compute the endomorphism $\End(S_L(k,N))$ of the semialgebra
$S_L(k,N)$ of cusp forms $f_k(z_N)$.

As the semialgebra $S_L(k,N)$ 
\resp{the cosemialgebra $S_R(k,N)$} of cusp forms $f_k(z_N)$
\resp{dual cusp forms $f_k(z^*_N)$} is defined on the quotient semigroup
\[
T^{(2)}_{\zit_H}(H) = T_2(\rit)\big/ T_2(\zit/N \zit)
\rresp{T^{t(2)}_{\zit_H}(H^*) = T^t_2(\rit)\big/ T^t_2(\zit/N \zit)}
\]
of lattices of $H$ \resp{$H^*$},
the bisemialgebra $S\RL(k,N)$ of cusp biforms $f_k(z^*_N)\times f_k(z_N)$ will be defined on the quotient bisemigroup
\[ G^{(2)}_{\zit_H}(H^*\times H) =\GL_2(\rit\times\rit)\big/\GL_2(\zit/N \zit)^2\]
where
\begin{align*}
G^{(2)}_{\zit_H}(H^*\times H) &= T^{t(2)}_{\zit_H}(H^*)\times T^{(2)}_{\zit_H}(H)\;, \\[6pt]
\GL_2(\rit\times\rit) &= T^t_2(\rit)\times T_2(\rit)\;, \\[6pt]
\GL_2(\zit/N \zit)^2 &= T^t_2(\zit/N \zit)\times T_2(\zit/N \zit)\;.
\end{align*}
\vskip 11pt

\subsection{Varying the weight of cusp (bi)forms}

The bisemigroup category CAF of cusp biforms $f_k(z^*_N)\times f_k(z_N)$ is then given by:
\Bena
\item the cusp biforms $(f_k(z^*_N)\times f_k(z_N))$, $k$ and $N$ varying;

\item {\bbf the (bi)morphisms of projection
\[ \Hom \L[
(f_k(z^*_N)\times f_k(z_N)),(f_h(z^*_N)\times f_h(z_N)) \R]\;, \qquad k>h\;, \]
sending cusp biforms
$(f_k(z^*_N)\times f_k(z_N))$ of weight $k$ into corresponding biforms
$(f_h(z^*_N)\times f_h(z_N))$ of weight $h$\/}.
\Ee

But, in order to understand the (bi)morphism of projection:
\[
\Phom^{(2)}_{k\to h}: \qquad f_k(z^*_N)\times f_k(z_N)
\quad \To \quad f_h(z^*_N)\times f_h(z_N)\]
from a cusp biform of weight $k$ onto a cusp biform of weight $h$, we have to give
{\bbf a new meaning to the weight of cusp forms\/}: this will be done in the next sections.
\vskip 11pt

\subsection{Proposition (Origin of the weight of cusp forms)}

{\em Let
\[ \phi ^{(2)}_k(z_N) \equiv T_k(n;N)_\varepsilon \ f_k(z_N)=\sum_{n=1}^\infty 
c_{nk}\ q^n_N\]
denote a weight $k$ two-dimensional cusp form of level $N$ as described in section~1.8,

where
\[ c_{nk}=\sum_{d|(m,n)} \varepsilon (d)\ d^{k-1}\ a_{mn|d^2} \qquad \text{and}
\qquad q^n_N=e^{2\pi inz_N}\;, \]

Let
\[ \phi ^{(k)}(z_{Nk})=\sum_{n=1}^\infty \lambda ^{(k)}_n\ e^{2\pi inz_{Nk}}\]
be the Fourier series development of a cusp form of real dimension $k$ ($k$ being an even integer) \cite{Pie2}
and level $N$ on $\cit^{k/2}$ above the weight $k$ two-dimensional cusp form
\[ \phi _k^{(2)}(z_N)\equiv f_k(z_N)\]
where:
\Bi
\item $\lambda ^{(k)}_n=\prod\limits_{d=1}^{k/2}\lambda ^{(2)}_{nd}$ with $\lambda ^{(2)}_{nd}=r_{S^1_{d_1}}\times r_{S^1_{d^2}}$ the product of the radii of two circles
$S^1_{d_1}$ and $S^1_{d^2}$;

\item $z_{Nk}=z_N^{(1)}+z_N^{(2)}+\dots+z_N^{(k)}$ is the sum of the elements of a $k$-tuple of complex numbers $z^{(\cdot)}_N$ of order $N$.
\Ei

Then, {\bbf there exists a map\/}:
\[ CP_{(k)\to 2}: \qquad \phi ^{(k)}(z_{Nk}) \quad \To \quad \phi ^{(2)}_k(z_N)\]
{\bbf in such a way that\/}, $\forall\ n$, $1\le n\le \infty $:
\Bena
\item {\bbf every two-dimensional semitorus
$T^{(2)}_{nN}(a,d\varepsilon (d)_k)=c_{nk}\ q^n_N$ of $\phi ^{(2)}_k(z_N)$ be the inverse image of a  $k$-dimensional real semitorus
$T^{(k)}_{nN}(\lambda ^{(k)}_n)=\linebreak \lambda ^{(k)}_n\ e^{2\pi inz_{Nk}}$ of $\phi ^{(k)}(z_{Nk})$\/};

\item the product $\lambda ^{(k)}_n=\prod\limits_{d=1}^{k/2}\lambda ^{(2)}_d$, representing the product of the radii of the generators of
$T^{(k)}_{nN}(\lambda ^{(k)}_n)$, be sent by the map
\[ CP_{(k)\to2}(n) : \qquad \lambda ^{(k)}_n \quad \To \quad c_{nk}\]
into the coefficient $c_{nk}$ of $\phi ^{(2)}_k(z_N)$ implying that $\|\lambda ^{(k)}_n\|=\|c_{nk}\|$.
\Ee
}
\vskip 11pt

\bpr
The weight $k$ of the two-dimensional cusp form of level $N$, $\phi ^{(2)}_k(z_N)$ then originates from the map
$CP_{(k)\to2}$ of a $k$-dimensional cusp form of level $N$, $\phi ^{(k)}(z_{Nk})$ into
$\phi ^{(2)}_k(z_N)$ at the conditions of this proposition.

The Fourier series development of $\phi^{(k)}(z_{Nk})$ was introduced in \cite{Pie2} from which it results that it is the sum of increasing embedded $k$-dimensional real semitori $T^{(k)}_{nN}(\lambda ^{(k)}_n)$.

{\bbf The conditions given by the map
\[ CP_{(k)\to2}(n): \qquad \lambda ^{(k)}_n \quad \to \quad c_{nk}\]
on the coefficients imply the inflation of the coefficients $c_n$ of a weight $2$ cusp form to the corresponding coefficients $c_{nk}$ of a weight $k$ cusp form\/} (according to section~1.8 and proposition~1.9) associated with the inflation  map
\[
I_{T_2\to k}: \qquad T_2(n,N)\ f_2(z_N) \quad \To \quad T_k(m,N)_\varepsilon \ f_k(z_N)\]
of a weight $2$ cusp form into a weight $k$ cusp form.

The origin of this inflation map $I_{T_2\to k}$ proceeds mainly from the factor
$d^{(k-1)}$ of $c_{nk}$ and, similarly, from the   ``normalization factor'' $(cz+d)^{-k}$ in the automorphic condition
\be
f(z)=(cz+d)^{-k}\ f\L(\dfrac{az+b}{cz+d}\R)\;.\tag*{\eop}
\ee
\vskip 11pt

\subsection{Proposition}

{\em
The fundamental domain $D^{(2)}_{nNk}$ of the weight $k$ two-dimensional cusp form
$\phi ^{(2)}_k(z_N)\equiv  f_k(z_N)$ is equal to the fundamental domain $D^{(k)}_{nN}$ of the corresponding $k$-dimensional cusp form $\phi ^{(k)}(z_{Nk})$ implying that
\[ T^{(2)}_{nN}(a,d,\varepsilon (d)_k)\simeq T^{(k)}_{nN}(\lambda ^{(k)}_n)\;, \qquad \forall\ n\;.\]
}
\vskip 11pt

\bpr
\Bena
\item The fundamental domain $D^{(2)}_{nNk}$ of the weight $k$ two-dimensional cusp form
$\phi ^{(2)}_k(z_N)$ is given by $q^1_N=e^{2\pi i1z_N}$ where $z_N$ is a complex point of order $(k\times N)$ according to propositions~1.4 and 1.14 while the fundamental domain
$D^{(k)}_{nN}$ of the corresponding $k$-dimensional cusp form
$\phi ^{(k)}(z_{Nk})$ is given by $q^1_{Nk}=e^{2\pi i1z_{Nk}}$
where $z_{Nk}$ is also a complex point of order $(k\times N)$ according to proposition~1.14.

Consequently, we have that:
\[ D^{(2)}_{nNk}=D^{(k)}_{nN}\;.\]
\item The domain $A^{(2)}_{nNk}$ of the $n$-th periodic cusp subform
$T^{(2)}_{nN}(a,d\varepsilon (d)_k)=f^{(n)}_k(z_N)$ is evidently given by:
\[ A^{(2)}_{nNk}=n\ D^{(2)}_{nNk}\]
while the domain $A^{(k)}_{nN}$ of the $n$-th periodic cusp subform $T^{(k)}_{nN}(\lambda ^{(k)}_n)$ of $\phi ^{(k)}(z_{Nk})$ is:
\[ A^{(k)}_{nN}=n\ D^{(k)}_{nN}\;.\]
Consequently, we have that:
\[ A^{(2)}_{nNk}=A^{(k)}_{nN}\]
and, then, that
\be
T^{(2)}_{nN}(a,d\varepsilon (d)_k)=T^{(k)}_{nN}(\lambda ^{(k)}_n)\;.\tag*{\eop}
\ee
\Ee
\vskip 11pt

\subsection{$k$-dimensional cusp biforms}

The $k$-dimensional real cusp biform of level $N$ associated with the $k$-dimensional real cusp form $\phi ^{(k)}(z_{Nk})$ of level $N$ is:
\[ \phi ^{(k)}(z^*_{Nk})\times\phi ^{(k)}(z_{Nk})\;.\]
As $ \phi ^{(k)}(z_{Nk})$ \resp{$\phi ^{(k)}(z^*_{Nk})$} is a function on the sum of the cosets of $T_k(\rit)/T_k(\zit/N \zit)$
\resp{$T^t_k(\rit)/T^t_k(\zit/N \zit)$} according to section~1.6,
$\phi ^{(k)}(z^*_{Nk})\times  \phi ^{(k)}(z_{Nk})$ will be a bifunction on the sum of the cosets of
\[
\GL_k(\rit\times\rit)\big/\GL_k(\zit/N \zit)^2
\equiv T^t_k(\rit)\times T_k(\rit)\big/T^t_k(\zit/N \zit)\times T_k(\zit/N \zit)\]
and will be written as follows:
\[
 \phi ^{(k)}_{z^*_{Nk}}
\times \phi ^{(k)}_{z_{Nk}}(\GL_k(\rit\times\rit)\big/\GL_k(\zit/N \zit)^2)\;.\]
\vskip 11pt

\subsection{Proposition (Changing the weight of cusp (bi)forms)}

{\em
Let 
\[
 \phi ^{(h)}_{z^*_{Nh}}
\times \phi ^{(h)}_{z_{Nh}}(\GL_h(\rit\times\rit)\big/\GL_h(\zit/N \zit)^2)\]
be the $h$-dimensional real cusp biform of level $N$ above the $2$-dimensional cusp biform
$(f_h(z^*_N)\times f_h(z_N))$ (also written
$(\phi ^{(2)}_h(z^*_N)\times \phi ^{(2)}_h(z_N))$) of weight $h$ and level $N$.

Then, the bimorphism (of projection) between $2$-dimensional cusp biforms of weights $k$ and $h$ (introduced in section~1.13):
\[
\Phom^{(2)}_{k\to h}: \qquad\quad
f_k(z^*_N)\times f_k(z_N)\quad \To \quad 
f_h(z^*_N)\times f_h(z_N) \;, \quad k>h\;, \]
directly depends on:
\Bean
\item the bimorphism (of projection)
\begin{multline*}
\Phom_{k\to h}: \qquad 
 \phi ^{(k)}_{z^*_{Nk}}
\times \phi ^{(k)}_{z_{Nk}}(\GL_k(\rit\times\rit)\big/\GL_k(\zit/N \zit)^2)\\
\To \quad
 \phi ^{(h)}_{z^*_{Nh}}
\times \phi ^{(h)}_{z_{Nh}}(\GL_h(\rit\times\rit)\big/\GL_h(\zit/N \zit)^2)
\end{multline*}
between $k$ and $h$-dimensional cusp biforms;

\item the decomposition
\begin{align*}
f_k(z^*_N)\times f_k(z_N) &= \bigoplus_{\ell=1}^{k/2} (f_{2_\ell}(z^*_N)\times f_{2_\ell}(z_N))\;,\\[11pt]
f_h(z^*_N)\times f_h(z_N) &= \bigoplus_{\ell=1}^{h/2} (f_{2_\ell}(z^*_N)\times f_{2_\ell}(z_N))\;,
\end{align*}
$h$ andf $k$ being even integers,

of weight $k$ and $h$ $2$-dimensional cusp biforms
$(f_{k}(z^*_N)\times f_{k}(z_N))$ and
$(f_{h}(z^*_N)\times f_{h}(z_N))$ 
into $2$-dimensional cusp biforms of weight $2$ according to the Langlands functoriality conjecture.
\Ee
\vskip 11pt

{\bbf This can be summarized in the following diagram\/}:
\[ \begin{psmatrix}[colsep=.5cm,rowsep=1cm]
{\renewcommand{\arraystretch}{1}
\begin{array}[b]{c}
\GL_2(\rit\times\rit)\qquad \\
\big/\GL_2(\zit/N \zit)^2 \\
= {\mathbb{GL}}_2 \end{array}} & & &
{\renewcommand{\arraystretch}{1}
\begin{array}[b]{c}
\GL_k(\rit\times\rit)\qquad \\
\big/\GL_k(\zit/N \zit)^2 \\
= {\mathbb{GL}}_k \end{array}} & & &
{\renewcommand{\arraystretch}{1}
\begin{array}[b]{c}
\GL_h(\rit\times\rit)\qquad \\
\big/\GL_h(\zit/N \zit)^2 \\
= {\mathbb{GL}}_h \end{array} }\\
\phi ^{(2)}_{z^*_{N}}
\times \phi ^{(2)}_{z_{N}}({\mathbb{GL}}_2) & & &
\phi ^{(k)}_{z^*_{Nk}}
\times \phi ^{(k)}_{z_{Nk}}({\mathbb{GL}}_k) & & &
\phi ^{(h)}_{z^*_{Nh}}
\times \phi ^{(h)}_{z_{Nh}}({\mathbb{GL}}_h)\\
 f_2(z^*_N)\times f_2(z_N)  & & &
f_k(z^*_N)\times f_k(z_N) & & &
 f_h(z^*_N)\times f_h(z_N) 
\psset{arrows=->,nodesep=3pt}
\everypsbox{\scriptstyle}
\ncline{1,1}{1,4}
\ncline{1,4}{1,7}
\ncline{1,1}{2,1}>{\Pi ^{(2)}}
\ncline{1,4}{2,4}>{\Pi ^{(k)}}
\ncline{1,7}{2,7}>{\Pi ^{(h)}}
\ncline{2,1}{2,4}^{\Phom_{2\to k}}
\ncline{2,4}{2,7}^{\Phom_{k\to h}}
\ncline{2,1}{3,1}>{{\rm Id}}
\ncline{2,4}{3,4}>{CP_{(k)\to2\RL}}
\ncline{2,7}{3,7}>{CP_{(h)\to2\RL}}
\ncline{3,1}{3,4}^{\Phom^{(2)}_{2\to k}}
\ncline{3,4}{3,7}^{\Phom^{(2)}_{k\to h}}
\end{psmatrix}\]
where:
\Bi
\item $\phi ^{(2)}_{z^*_{N}}\times \phi ^{(2)}_{z_{N}}(\GL_2(\rit\times\rit)
/\GL_2(\zit/N \zit)^2)\equiv f_2(z^*_N)\times f_2(z_N)$ is a two-dimensional cusp biform of weight $2$;
\item the map $\Pi ^{(k)}$ is a cuspidal representation of the quotient bisemigroup
$(\GL_k(\rit\times\rit)
/\GL_k(\zit/N \zit)^2)$ given by the corresponding cusp biform and depending on a toroidal compactification \cite{Pie2}.
\Ei
}
\vskip 11pt

\bpr
{\bbf The bimorphism (of projection) $\Phom_{k\to h}$ is directly related to the Langlands functoriality conjecture\/} transposed from cuspidal representations of bilinear algebraic semigroups \cite{Pie4} to cuspidal representation of quotient bisemigroups 
$\GL_k(\rit\times\rit)
/\GL_k(\zit/N \zit)^2$ and asserting that the cuspidal representation $\Pi ^{(k)}(\GL_k(\rit\times\rit)
/\GL_k(\zit/N \zit)^2)$ of the quotient bisemigroup
$(\GL_k(\rit\times\rit)
/\GL_k(\zit/N \zit)^2)$ is orthogonally completely reducible if it decomposes diagonally according to \cite{Pie4}:
\[
\Pi ^{(k)}(\GL_k(\rit\times\rit)
\big/\GL_k(\zit/N \zit)^2)
= \bigoplus_{\ell=1}^{k/2} \Pi ^{(2_\ell)}(\GL_{2_\ell}(\rit\times\rit)
\big/\GL_{2_{\ell}}(\zit/N \zit)^2)\;.
\]
This corresponds to the injective map $\Phom_{2\to k}$ which is above the injective map
 $\Phom_{2\to k}$ leading to {\bbf the decomposition
\[
 f_k(z^*_N)\times f_k(z_N)= \bigoplus_{\ell=1}^{k/2} (f_{2_\ell}(z^*_N)\times f_{2\ell}(z_N))\]
 of weight $k$ cusp biforms into weight $2$ cusp biforms\/}.
\vskip 11pt

The surjective morphism $\Phom^{(2)}_{k\to h}$ is then directly associated with:
\Bean
\item the deflation of the product $(c^*_{nh}\times c_{nh})$ of the coefficients of the two-dimensional cusp biforms $f_h(z^*_N)\times f_h(z_N)$ of weight $h$ from the product
$(c^*_{nk}\times c_{nk})$ of the coefficients of the cusp biform $f_k(z^*_N)\times f_k(z^*_N)$ of weight $k$, $1\le n\le t\le \infty $, according to proposition~1.14 and giving the increasing sequence
\[ c^*_{n2}\times c_{n2}\To c^*_{nh}\times c_{nh}\To c^*_{nk}\times c_{nk}\;;\]

\item the increasing sequence
\[
(z^*_N\times z_N)_2 \To (z^*_N\times z_N)_h \To  (z^*_N\times z_N)_k\]
of cusp bipoints, where $(z^*_N\times z_N)$ is a complex bipoint of order $(N\times 2)$, $(z^*_N\times z_N)_h$ is a complex bipoint of order $(N\times h)$ and
$(z^*_N\times z_N)_k$ is a complex bipoint of order $(N\times k)$, with $k>h$.\epr
\Ee

\section{Cusp (bi)forms from the Langlands global program}

Cusp (bi)forms are characterized by the two integer numbers $k$ and $N$ where $k$ refers to the geometric dimension of a cusp (bi)form above a cusp (bi)form of weight $2$ while $N$ will be recalled to refer to an algebraic dimension.

Indeed, we shall at first prove that the quotient bisemigroup $\GL_k(\rit\times\rit)/\GL_k(\zit/N \zit)^2$, on which a cusp biform can be defined, is covered by the (algebraic) bilinear semigroup $\GL_k(F_{\o v}\times F_v)$ over the product $(F_{\o v}\times F_v)$ of sets $F_{\o v}$ and $F_v$ of pseudoramified transcendental extensions referring respectively to the lower and upper half space.
\vskip 11pt

\subsection{Algebraic and transcendental quanta}

\Bean
\item Let then $\wt F$ denote a set of finite algebraic extensions of a number field $k$ of characteristic $0$: $\wt F$ is assumed to be a set of {\bbf symmetric splitting fields composed of the left and right algebraic extension semifields $\wt F_L$ and $\wt F_R$\/} being respectively the set of complex and conjugate complex simple roots of the polynomial ring $k[x]$.

In the real case, the symmetric splitting fields are the left and right symmetric splitting semifields $\wt F_L^+$ and $\wt F_R^+$ composed of the set of positive and symmetric negative simple real roots.

\item The left and right equivalence classes of infinite Archimedean completions of $\wt F_L$ and $\wt F_R$ are the left and right symmetric infinite complex places
\[ 
\omega =\{\omega _1,\dots,\omega _n,\dots,\omega _t\} \quad \and \quad
{\o \omega} =\{{\o \omega }_1,\dots,{\o \omega} _n,\dots,{\o \omega} _t\} \;, \qquad t\le \infty \;.\]
In the real case, the infinite places are similarly
\[ 
v =\{v _1,\dots,v _n,\dots,v _t\} \quad \and \quad
{\o v} =\{{\o v }_1,\dots,{\o v} _n,\dots,{\o v} _t\} \;,\]
covering the respective infinite complex places as described subsequently.

\item All these (pseudoramified) completions, corresponding to the transcendental extensions, proceed from the associated algebraic extensions by a suitable isomorphism of compactification \cite{Pie5} and are built from irreducible real subcompletions
$F_{v^1_n}$
\resp{$F_{\o v^1_n}$} characterized by a transcendence degree
\[
\tr\cdot d\cdot F_{v^1_n}\big / k=\tr\cdot d\cdot F_{\o v^1_n}\big/ k=N
\quad \text{equal to} \quad
[\wt F_{v^1_n}:k]= [\wt F_{\o v^1_n}:k]= N\;,\]
which is the Galois extension degree of the associated irreducible algebraic closed subsets
$\wt F_{v^1_n}$ and 
$\wt F_{\o v^1_n}$.

{\bbf All these irreducible (unitary)  subcompletions \resp{subextensions} are assured  to be transcendental \resp{algebraic} quanta\/} \cite{Pie6}.

\item The pseudoramified real extensions are characterized by degrees:
\[
[\wt F_{v_n}:k]=[\wt F_{\o v_n}:k]=*+n\ N\;, \qquad 1\le n\le t\le\infty \]
which are integers modulo $N$, $\zit/N \zit$, where:
\Bi
\item $\wt F_{v_n}$ and 
$\wt F_{\o v_n}$ are algebraic extensions corresponding to the infinite completions
$F_{v_n}$ and
$F_{\o v_n}$, respectively at the $v_n$-th and
 $\o v_n$-th symmetric real pseudoramified places;
 \item $*$ denotes an integer inferior to $N$, taken generally to be ``$0$''.
 \Ei
 
Similarly, the pseudoramified complex extensions $\wt F_{\omega _n}$ and $\wt F_{\o\omega _n}$, corresponding to the completions 
$  F_{\omega _n}$ and $  F_{\o\omega _n}$ at the infinite places $\omega _n$ and $\o\omega _n$,
are characterized by extension degrees:
\[
[\wt F_{\omega _n}:k]=[\wt F_{\o \omega _n}:k]=(*+n\ N)\ m_n\;,  \]
 where $m_n$ is the multiplicity of the $n$-th real extension $F_{v_n}$ and $F_{\o v_n}$ covering its $n$-th complex equivalent.
 
 Let then 
$\wt F_{v_{n,m_n}}$ 
\resp{$\wt F_{\o v_{n,m_n}}$} denote a \lr pseudoramified real extension equivalent to 
$\wt F_{v_n}$ 
\resp{$\wt F_{\o v_n}$}.

The corresponding pseudounramified real extensions
$\wt F^{nr}_{v_{n,m_n}}$ and
$\wt F^{nr}_{\o v_{n;m_n}}$ are characterized by their global residue degree (case $N=1$):
\[
f_{v_n}=[\wt F^{nr}_{v_{n,m_n}}:k]=j \quad \text{and} \quad
f_{\o v_n}=[\wt F^{nr}_{\o v_{n,m_n}}:k]=j \;.\]

\item If the orders of  the Galois sub(semi)groups correspond to the class zero of the integers modulo $N$, then these Galois sub(semi)groups $\Gal(\wt F_{v_n}/k)$
\resp{$\Gal(\wt F_{\o v_n}/k)$} are global Weil sub(semi)groups of extensions
$\wt F_{v_n}$
\resp{$\wt F_{\o v_n}$} constructed from $n$ algebraic quanta.

{\bbf By an isomorphism of compactification, these $n$ noncompact algebraic quanta of 
$\wt F_{v_n}$
\resp{$\wt F_{\o v_n}$} are sent into the corresponding compactified $n$ transcendental compact quanta $F_{v_n}$
\resp{$F_{\o v_n}$}\/}.

As in the Galois case, there is a one-to-one correspondence between all transcendence extension subfields:
\[
F_{v_1} \subset \dots \subset F_{v_n} \dots  \subset  F_{v_t} 
\rresp{F_{\o v_1} \subset \dots \subset  F_{\o v_n} \dots  \subset  F_{\o v_t} }
\]
and the set of all (normal) sub(semi)groups of these:
\begin{align*}
\Aut_k(F_{v_1}) &\supset \dots \supset 
\Aut_k(F_{v_n}) \supset  \dots \supset 
\Aut_k(F_{v_t}) \\
\rresp{\Aut_k(F_{\o v_1}) &\supset  \dots \supset 
\Aut_k(F_{\o v_n}) \supset \dots \supset 
\Aut_k(F_{\o v_t})}
\end{align*}
taking into account that 
$\Gal(\wt F_{v_n}/k)\simeq \Aut_k(F_{v_n})$
\resp{$\Gal(\wt F_{\o v_n/k})\simeq \Aut_k(F_{\o v_n})$}.
 \Ee
 \vskip 11pt
 
\subsection{Abstract bisemivarieties}

\Bean
\item Let 
\[
\GL_k(\wt F_{\o v}\times\wt F_{v})\equiv T^t_k(\wt F_{\o v})\times T_k(\wt F_{v})\]
 be {\bbf the algebraic bilinear semigroup of matrices\/} over the product of sets
$\wt F_v=\{\wt F_{v_1},\dots,\wt F_{v_n},\dots,\wt F_{v_t}\}$ and
$\wt F_{\o v}=\{\wt F_{\o v_1},\dots,\wt F_{\o v_n}, \dots,\wt F_{\o v_t}\}$ of pseudoramified algebraic extensions in such a way that its representation bisemispace is given by the tensor product
$\wt M_{v_R}\otimes \wt M_{v_L}$ of a right 
$T^t_k(\wt F_{\o v})$-semimodule $\wt M_{v_R}$ by a left
$T_k(\wt F_{v})$-semimodule $\wt M_{v_L}$.

Considering the monomorphism:
\[ \wt \sigma _{v_R}\times \wt \sigma _{v_L}: \qquad
W ^{ab}_{\wt F_{\o v}}\times W ^{ab}_{\wt F_{v}} \quad \To \quad \GL_k(\wt F_{\o v}\times\wt F_v)\]
from the product of the global Weil semigroup
$W ^{ab}_{\wt F_{\o v}}\times W ^{ab}_{\wt F_{v}}$ into 
$\GL_k(\wt F_{\o v}\times\wt F_v)$ and the isomorphism:
\[
W ^{ab}_{\wt F_{\o v}}\times W ^{ab}_{\wt F_{v}}\quad  \overset{\sim}{\To}\quad  \Aut_k(F_{\o v})\times \Aut_k(F_v)\]
of $(W ^{ab}_{\wt F_{\o v}}\times W ^{ab}_{\wt F_{v}})$ with respect to the product
$\Aut_k(F_{\o v}) \times \Aut_k(F_v)$ of the automorphisms (semi)groups of the sets $F_{\o v}$ and $F_v$ of increasing transcendental extensions, we have the commutative diagram:
\[ \begin{psmatrix}[colsep=.5cm,rowsep=1cm]
W ^{ab}_{\wt F_{\o v}}\times W ^{ab}_{\wt F_{v}} & & & & & & &G^{(k)}(\wt F_{\o v}\times\wt F_{v})\\
\Aut_k(F_{\o v}) \times \Aut_k(F_v) & & & & & & &G^{(k)}( F_{\o v}\times F_{v})
\psset{arrows=->,nodesep=3pt}
\ncline{1,1}{1,8}^{\wt \sigma _{v_R}\times \wt \sigma _{v_L}}
\ncline{2,1}{2,8}^{\sigma _{v_R}\times  \sigma _{v_L}}
\ncline{1,1}{2,1}>{\wr}
\ncline{1,8}{2,8}>{\wr}
\end{psmatrix}\]
where $G^{(k)}(\wt F_{\o v}\times\wt F_{v})\equiv \wt M_{v_R}\otimes\wt M_{v_L}$ and
$G^{(k)}(  F_{\o v}\times  F_{v})$ is an abstract bisemivariety covered by the affine bisemigroup
$G^{(k)}(\wt F_{\o v}\times\wt F_{v})$.

\item Let
{\bbf $G^{(k)}( F^{nr}_{\o v}\times F^{nr}_{v})$ be the abstract bisemivariety of dimension $k$\/} over the product
$(  F^{nr}_{\o v}\times  F^{nr}_{v})$ of the sets 
$ F^{nr}_{\o v}=\{ F^{nr}_{\o v_1},\dots,F^{nr}_{\o v_n}\}$ and
$ F^{nr}_{  v}=\{ F^{nr}_{ v_1},\dots,F^{nr}_{ v_n}\}$ of increasing pseudounramified (case $N=1$) extensions.

Then, the kernel
$\Ker (G^{(k)}_{F\to F^{nr}})$ of the map:
\[ G^{(k)}_{F\to F^{nr}} : \qquad G^{(k)}(F_{\o v}\times F_v)\quad \To \quad G^{(k)}(
F^{nr}_{\o v}\times F^{nr}_v)\]
is the minimal (unitary) bilinear parabolic subsemigroup $P^{(k)}(F_{\o v^1}\times F_{v^1})$ over the product of the sets
$F_{\o v^1}=\{ F_{\o v^1_1},\dots,F_{\o v^1_n},\dots\}$ and 
$F_{v^1}=\} F_{v^1_1},\dots,F_{v^1_n},\dots\}$ of unitary Archimedean pseudoramified completions.

$G^{(k)}(F_{\o v}\times F_v)$ acts on the unitary bilinear parabolic subsemigroup
$P^{(k)}(F_{\o v^1}\times F_{v^1})$ by conjugation \cite{Pie5}.

\item At every infinite Archimedean biplace $(\o v_n\times v_n)$ corresponds
{\bbf a conjugacy class $cg^{(k)}_{v\RL}[n]$ of the abstract bisemivariety $G^{(k)}(F_{\o v}\times F_v)$\/} whose number of representatives corresponds to the  number of equivalent transcendental extensions of $F_{\o v_n}\times F_{v_n}$.

{\bbf The $n$-th conjugacy class representative $g^{(k)}_{v\RL}[n]$ is defined over $n$ transcendental biquanta (a biquantum being the product of a right quantum by its symmetric left equivalent)\/} in such a way that the number of biquanta in 
 $g^{(k)}_{v\RL}[n]$ is 
\[ n_{\rm biq}(g^{(k)}_{v\RL}[n])=n^k\;.\]

\item Let $G^{(k)}( F_{\o \omega }\times F_{\omega })$ denote the complex abstract bisemivariety which is a $\GL_{k/2}(F_{\o \omega }\times F_\omega )$-bisemimodule
$ M_{\omega _R}\otimes M_{\omega _L}$ and the representation space of the bilinear semigroup of matrices $\GL_{k/2}(F_{\o \omega }\times F_\omega )$ over the product
$F_{\o \omega}\times F_\omega $ of sets of complex pseudoramified increasing transcendental extensions.

{\bbf Assume that each conjugacy class representative
$g^{(k)}_{\omega \RL}[n]$ of $G^{(k)}(F_{\o\omega }\times F_\omega )$ is unique and can be covered by $m_n$ real conjugacy class representatives
$\{ g^{(k/2)}_{v\RL}[n,m_n]\}$\/} of $G^{(k/2)}
(F_{\o v}\times F_v)$, $1\le n\le t\le\infty $.

So, the complex bipoints of $G^{k}(F_{\o \omega }\times F_\omega )$ are in one-to-one correspondence with the real bipoints of
 $G^{k}(F_{\o v }\times F_v )$ and we have the inclusion:
 \[G^{(k)}(F_{\o v}\times F_v )\big/ G^{(k/2)}(F_{\o  v}\times F_v )\quad \hookrightarrow \quad 
G^{(k)}(F_{\o \omega }\times F_\omega  )\;.\]
\Ee
\vskip 11pt

\subsection{Cuspidal representation of complex (algebraic) bilinear semigroups}

\Bean
\item {\bbf Providing a cuspidal representation of the complex abstract bisemivariety
$G^{(k)}(F_{\o\omega }\times F_\omega )$ consists in finding a cusp biform of dimension $k$ on
$G^{(k)}(F_{\o\omega }\times F_\omega )$ by summing the cuspidal subrepresentations of its conjugacy class representatives\/} \cite{Pie2}.

Let then
\[
\gamma ^T_{F_{\omega _n}}: \quad F_{\omega _n}\To F^T_{\omega _n}
\rresp{\gamma ^T_{F_{\o\omega _n}}: \quad F_{\o\omega _n}\To F^T_{\o\omega _n}}\;, \quad\forall\ n\;, \]
be {\bbf the toroidal isomorphism\/} mapping each \lr complex transcendental extension
$F_{\omega _n}$ \resp{$F_{\o\omega _n}$} into its 
toroidal equivalent $F^T_{\omega _n}$ \resp{$F^T_{\o \omega _n}$} which is a complex
one-dimensional semitorus localized in the upper \resp{lower} half space.

\item {\bbf Each \lr function on the conjugacy class representative
$g^{(k)}_{\omega ^T_L}[n]\in T^{(k)}(F^T_\omega )\subset G^{(k)}(F^T_{\o\omega }\times F^T_\omega )$\/}
\resp{$g^{(k)}_{\omega ^T_R}[n]\in T^{(k)}(F^T_{\o\omega })$} {\bbf is a function\/} \resp{cofunction}
$\phi _L(T^k_L[n])$ \resp{$\phi _R(T^k_R[n])$}, $g^{(k)}_{\omega ^T_R}[n]\equiv T^k_L[n]$, 
{\bbf on the even $k$-dimensional real semitorus $T^k_L[n]$\/} \resp{$T^k_R[n]$} having the analytic development:
\[ \phi _L(T^k_L[n])=\lambda (k,n)\ e^{2\pi in z_{Nk}}\qquad 
\rresp{\phi _R(T^k_R[n])=\lambda^*(k,n)\ e^{2\pi in z^*_{Nk}}}\]
where:
\Bi
\item $z_{Nk}=z^{(1)}_N+z^{(2)}_N+\dots+z^{(k)}_N$
\resp{$z^*_{Nk}=z^{*(1)}_N+z^{*(2)}_N+\dots+z^{*(k)}_N$} is a complex \resp{conjugate complex} point of order $(k\times N)$ according to proposition~1.14;

\item $\lambda (k,n)\times \lambda ^*(k,n)$ is the product of the eigenvalues of the $n$-th coset representative of the product, right by left, of Hecke operators \cite{Pie2}.
\Ei

\item {\bbf This \lr function $\phi _L(T^k_L[n])$ \resp{$\phi _R(T^k_R[n])$} constitutes the cuspidal representation\/} $\Pi ^{(n)}(g^{(k)}_{\omega _L}[n])$ \resp{$\Pi ^{(n)}(g^{(k)}_{\omega _R}[n])$} {\bbf of the $n$-th conjugacy class representative\/} of $G^{(k)}(F_\omega )$ \resp{$G^{(k)}(F_{\o\omega} )$} in such a way that {\bbf the cusp biform of
$\GL_{k/2}(F_{\o\omega }\times_D F_\omega )$ is given by the Fourier biseries\/}:
\begin{align*}
\Pi (\GL_{k/2}(F_{\o\omega }\times_D F_\omega )
&= \bigoplus^t_{n=1} \Pi ^{(n)}(g^{(k)}_{\omega \RL}[n])\\[15pt]
&=\L( \sum_n \lambda ^*(k,n)\ e^{2\pi inz^*_{Nk}}\R)\times_D
\L( \sum_n \lambda (k,n)\ e^{2\pi inz_{Nk}}\R)\end{align*}
where $\times_D$ is the diagonal product.
\Ee
\vskip 11pt

\subsection[Cuspidal representation of real (algebraic) bilinear semigroups]{Cuspidal representation of real (algebraic) bilinear\linebreak semigroups}

A real cuspidal representation, covering the complex cuspidal representation\linebreak
$\Pi (\GL_{k/2}(F_{\o\omega }\times_D F_\omega )$, can be obtained for the real diagonal bilinear semigroup
$G^{k}(F_{\o v }\times_D F_v )$ by summing the cuspidal subrepresentations
$\Pi ^{(n)}(g^{(k)}_{v\RL}[n,m_n])$ of its conjugacy class representatives
$g^{(k)}_{v\RL}[n,m_n]$.

Every \lr function on the set of conjugacy class representatives\linebreak
$\{g^{(k)}_{v^T_L}[n,m_n]\}_{m_n}$
\resp{$\{g^{(k)}_{v^T_R}[n,m_n]\}_{m_n}$} is a function \resp{cofunction}
$\psi _L(T^k_L[n,m_n])$
\resp{$\psi _R(T^k_R[n,m_n])$} on the $k$-dimensional real semitorus
$T^k_L[n,m_n]$
\resp{$T^k_R[n,m_n]$} localized in the upper \resp{lower} half space, covered by $m_n$ semitori of dimension $k/2$ and having the analytic development:
\begin{align*}
\psi _L(T^k_L[n,m_n]) &= \sum_{m_n} \lambda (k/2,n,m_n)\ e^{2\pi in_{m_n}x_{Nk/2}}\\[11pt]
\rresp{\psi _R(T^k_R[n,m_n]) &= \sum_{m_n} \lambda^* (k/2,n,m_n)\ e^{-2\pi in_{m_n}x_{Nk/2}}}
\end{align*}
where $x_{Nk}$ is a real point of order $N\times (k/2)$.

{\bbf The real cuspidal representation
$\Pi (\GL_k(F_{\o v}\times_D F_v)$ is then given by the $k$-\linebreak dimensional global elliptic 
$\Pi (\GL_k(F_{\o v}\times_D F_v)$-bisemimodule\/}
\[
\ELLIP\RL(k,n,m_n)=\sum_n\sum_{m_n}
(\lambda ^*(k/2,n,m_n)\ e^{-2\pi in_{m_n}x_{Nk/2}}\times_D
(\lambda (k/2,n,m_n)\ e^{2\pi in_{m_n}x_{Nk/2}}\;, \]
{\bbf covering the $k$-dimensional cusp biform
$\Pi (\GL_{k/2}(F_{\o\omega} \times_D F_\omega ))$\/}
in the sense that
\[ \ELLIP\RL(k,n,m_n) \quad \hookrightarrow \quad \Pi (\GL_{k/2}(F_{\o\omega} \times_D F_\omega ))\;.\]
\vskip 11pt

\subsection{(Bi)functor FLGC of the Langlands global correspondence(s)}

Let {\bbf CABG denote the bisemigroup category of algebraic bilinear semigroups
$\GL_{k/2}(\wt F_{\o\omega }\times \wt F_\omega )$\/} given by
\Bean
\item the algebraic bilinear semigroups $\GL_{k/2}(\wt F_{\o\omega }\times \wt F_\omega )$, the geometric dimension $k$ varying;

\item the (bi)morphisms of projection
\[ \Hom\L(\GL_{k/2}(\wt F_{\o\omega }\times \wt F_\omega ),\GL_{h/2}(\wt F_{\o\omega }\times \wt F_\omega )\R)\;, \qquad k>h\;, \]
sending $\GL_{k/2}(\wt F_{\o\omega }\times \wt F_\omega )$ into the algebraic bilinear semigroup
$\GL_{h/2}(\wt F_{\o\omega }\times \wt F_\omega )$.
\Ee

Similarly, let {\bbf CBCF denote the bisemigroup category of (complex) cuspidal representations
$\Pi (\GL_{k/2}(\wt F_{\o\omega }\times \wt F_\omega ))$\/} of the algebraic bilinear semigroup
$\GL_{k/2}(\wt F_{\o\omega }\times \wt F_\omega )$ whose (bi)morphisms of projection are:
\[\Hom\L(\Pi (\GL_{k/2}(\wt F_{\o\omega }\times \wt F_\omega )),\Pi (\GL_{h/2}(\wt F_{\o\omega }\times \wt F_\omega )\R)\;.\]

Then, {\bbf there exists a covariant (bi)functor FLGC associated with the Langlands global correspondences\/} \cite{Pie2}:
\begin{align*}
{\rm FLGC} :  &&{\rm CABG}  \To& \quad {\rm CBCF}\\[6pt]
&&\GL_{k/2}(\wt F_{\o\omega }\times\wt F_\omega )   \To &\quad \Pi (\GL_{k/2}(\wt F_{\o\omega }\times \wt F_\omega )\\[6pt]
{\rm FLGC}(\Phom_{k\to h}):  &&\Pi (\GL_{k/2}(\wt F_{\o\omega }\times \wt F_\omega ))
  \To &\quad 
\Pi (\GL_{h/2}(\wt F_{\o\omega }\times \wt F_\omega ))\end{align*}
which is a (bi)function assigning:
\Bean
\item to each algebraic bilinear semigroup
$\GL_{k/2}(\wt F_{\o\omega }\times \wt F_\omega ) \in {\rm CABG}$ its cuspidal representation
$\Pi (\GL_{k/2}(\wt F_{\o\omega }\times \wt F_\omega ))\in{\rm CBCF}$;

\item to each morphism $\Phom^{\rm (alg)}_{k\to h}: \GL_{k/2}(\wt F_{\o\omega }\times \wt F_\omega )\to \GL_{h/2}(\wt F_{\o\omega }\times \wt F_\omega )$ a morphism
${\rm FLGC}(\Phom^{\rm (alg)}_{k\to h}): \Pi (\GL_{k/2}(\wt F_{\o\omega }\times \wt F_\omega ))\to \Pi (\GL_{h/2}(\wt F_{\o\omega }\times \wt F_\omega ))$.
\Ee
\vskip 11pt

\subsection{Hecke eigenbivalues as Galois representation}

{\bbf Consider now the two-dimensional case $\GL_2(F_{\o v}\times F_v)$.

The ring of the endomorphisms of the $\GL_2(F_{\o v}\times F_v)$-bisemimodule
$(M_{v_R}\otimes M_{v_L})$\/}, decomposing it into the set of subbisemimodules
$\{M_{F_{\o v_{n,m_n}}}\otimes M_{F_{v_{n,m_n}}}\}_{n,m_n}$ or conjugacy class representatives
$\{g^{(2)}_{v\RL}[n]\}_{n,m_n}$ (see section~2.2), according to the $(\zit\big/N \zit)^2$-bisemilattices
$(\Lambda ^{(2)}_{\o v_{n,m_n}}\otimes\Lambda ^{(2)}_{v_{n,m_n}})$ \cite{Pie1}, {\bbf is generated over $(\zit\big/N \zit)^2$ by the product $(T_{n_R}\otimes T_{n_L})$ of Hecke operators\/}
$T_{n_R}$ over $T_{n_L}$ for $n\nmid N$ and by the product $(U_{n_R}\times U_{n_L})$ of Hecke operators $U_{n_R}$ and $U_{n_L}$ for $n\mid N$: it is noted $T_H(N)_R\otimes T_H(N)_L$ (weight two case).

The coset representative of $U_{n_L}$ \resp{$U_{n_R}$}, referring to the upper \resp{lower} half plane, can be chosen to be upper \resp{lower} triangular and given by the integral matrix
$\L(\BsM 1 & b_N \\ 0 & n_N\EsM\R)$
\resp{$\L(\BsM 1 & 0 \\ b_N & n_N\EsM\R)$} of the congruence subgroup
$\Gamma _L(N)$
\resp{$\Gamma _R(N)$} in $\GL_2(\zit)$ in such a way that:
\Bi
\item $n_N=*+n\ N$ where $n_N$ is the number of transcendental subfields in the $n$-th conjugacy class representative $g^{(2)}_N[n]$;

\item $b_N=*+b\cdot N$ refers to a phase shift;
\Ei
where $*$ denotes an integer inferior to $N$.

The coset representative of $(U_{n_R}\times U_{n_L})$ will be given by
\[
U_{n_R}\times U_{n_L}=\L[ \BM 1 & b_N \\ 0&1\EM 
\BM 1&0\\ b_N & 1\EM \R]\BM 1&0\\ 0 & n^2_N\EM\]
where $\L(\BsM 1 & b_N\\ 0&1\EsM \R) \L(\BsM 1 & 0 \\ b_N&1\EsM \R) $ is an element of the nilpotent group of matrices $u_2(b_N)\cdot u_2(b_N)^t$ and of the unimodular decomposition group
$D_{n^2_N}$ acting on the split Cartan subgroup element
$\alpha _{n^2_N}=\L(\BsM 1&0\\ 0&n^2_N\EsM \R)$.

$(U_{n_R}\times U_{n_L})$ then factorizes according to the bilinear Gauss decomposition into nilpotent and diagonal parts \cite{Pie1}.

Let {\bbf $\lambda ^2_+(n^2_N,b^2_N)$ and 
$\lambda ^2_-(n^2_N,b^2_N)$ be the eigenvalues of
$(D_{n^2_N;b^2_N}\times\alpha _{n^2_N})\equiv \GL_2(T_H(N)_R\otimes T_H(N)_L)$ corresponding to the irreducible semisimple (pseudo)\-ramified representation:
\[ \rho _{\lambda ^2_\pm}: \qquad \Gal ( \wt F_{\o v_n}/k)\times \Gal ( \wt F_{v_n}/k)
\quad \To \quad \GL_2(T_H(N)_R)\otimes T_H(N)_L)\]
associated with a weight two cusp form\/} in such a way that $\lambda ^2_+$ and $\lambda ^2_-$ verify \cite{Pie1}:
\begin{align*}
{\rm trace}\ \rho _{\lambda ^2_\pm} &= 1+b^2_N+n^2_N\;, \\[6pt]
\det \rho _{\lambda ^2_\pm} &= \lambda ^2_+(n^2_N,b^2_N)\times \lambda ^2_-(n^2_N,b^2_N) = n^2_N\;.\end{align*}
Indeed, the eigenvalues $\lambda ^2_\pm$ are:
\[
\lambda ^2_\pm = \F{(1+b^2_N+n^2_N)\pm[(1+b^2_N+n^2_N)-4n^2_N]^{1/2}}{2}
\]
and the characteristic polynomial of $\rho _{\lambda ^2_\pm}$ has the form:
\[
X^2-{\rm trace}\ \rho _{\lambda ^2_\pm X}  +\det\rho _{\lambda ^2_\pm}=0\;.\]

All that can be summarized in the commutative diagram:
\[ \begin{psmatrix}[colsep=.5cm,rowsep=1cm]
\Gal (\wt F_{\o v}/k)\times \Gal (\wt F_{v}/k) && 
\GL_2 (  F_{\o v}\times   F_{v}) \\
\GL_2(T_H(N)_R\otimes T_H(N)_L) &&
\GL_2(F^+_R\times F^+_L)\big/\GL_2(\zit\big/ N \zit)^2)
\psset{arrows=->,nodesep=3pt}
\ncline{1,1}{1,3}
\ncline{2,1}{2,3}
\ncline{2,3}{1,3}
\ncline{1,1}{2,1}
\end{psmatrix}
\]
Note that, in the weight $k$ case, the Cartan subgroup element
$\alpha _{n^2_N}=\L(\BsM 1&0\\ 0&n^2_N\EsM \R)$ is then given by
$\alpha _{n^2_{N,k}}=\L(\BsM a^2_N&0\\ 0&d^{2k-2}_N\EsM \R)$ where
$n^2_{N,k}=\sum\limits_{d\mid n} a^2_N\cdot d^{2k-2}_N\cdot \epsilon ^2(d)$ denotes the number of transcendental quanta, i.e. $((k/2)\times n)^2$ referring to sections 1.8 and 2.2.
\vskip 11pt

\subsection{Two-dimensional cusp biforms and global elliptic bisemimodules}

The  $k$-dimensional bilinear cuspidal representation
$\Pi (\GL_{k/2}(F_{\o\omega }\times_{(D)} F_\omega )$ of the complex bilinear semigroup
$\GL_{k/2}(F_{\o\omega }\times_{D} F_\omega )$ is a $k$-dimensional cusp biform whose Fourier biseries are given by:
\[
\Pi (\GL_{k/2}(F_{\o\omega }\times_{(D)} F_\omega )=
\L( \sum_n\lambda ^*(k,n)\ e^{2\pi inz^*_{Nk}}\R)
\times _D\L(\sum_n\lambda (k,n)\ e^{2\pi inz_{Nk}}\R)
\]
according to section~2.3.

{\bbf The corresponding $2$-dimensional (``Weil'') cusp biform of weight $k$ can be reached by considering the projective (bi)map\/}:
\[
CP_{(k)\to(2)} : \qquad \Pi (\GL_{k/2}(F_{\o\omega }\times_{D} F_\omega ) \quad \To \quad 
{}^{(\omega )}f_k(z^*_N) \times_D{}^{(\omega )}f_k(z_N)\;, \]
introduced in proposition~1.14,

where ${}^{(\omega )}f_k(z^*_N) \times_D{}^{(\omega )}f_k(z_N)$ are given by
\[
{}^{(\omega )}f_k(z^*_N) \times_D{}^{(\omega )}f_k(z_N)
= \L(\sum_n\lambda ^*_k(2,n)\ e^{2\pi inz^*_N}\R)\times _D
\L(\sum_n\lambda _k(2,n)\ e^{2\pi inz_N}\R)
\] where:
\Bi
\item $\lambda _k(2,n)\equiv [\lambda ^2_{k_\pm}(n^2_N,b^2_N)]^{1/2}$ is the product of the radii of two orthogonal circles at $(k/2)\times n$ transcendental quanta referring to proposition~1.9 and section~2.6;

\item $z_N$ is a complex point of order $(N\times k)$ written now $z_{N-k}$.
\Ei
\vskip 11pt

Similarly, the real $k$-dimensional bilinear cuspidal representation
$\Pi (\GL_k(F_{\o v} \times_{(D)} F_v))$ of the bilinear semigroup
$\GL_k(F_{\o v}\times_{(D)} F_v)$ is a $k$-dimensional elliptic
$\Pi (\GL_k(F_{\o v}\times_{(D)} F_v))$-bisemimodule:
\[
\ELLIP\RL (k,n,m_n)  =\sum_n\sum_{m_n}(\lambda ^*(k/2,n,m_n)\ e^{-2\pi in_{m_n}x_{Nk/2}}\times_D
(\lambda (k/2,n,m_n)\ e^{2\pi in_{m_n}x_{Nk/2}}
\]
which is a real analytic $k$-dimensional cusp biform.

{\bbf The corresponding $2$-dimensional global elliptic bisemimodule of weight $k$ will be reached by the projective (bi)map\/}:
\[
CP^{\ELLIP}_{(k)\to (2)} : \qquad
\ELLIP\RL(k,n,m_n) \quad \To \quad \ELLIP_{k\RL}(2,n,m_n)\]
where:
\begin{multline*}
 \ELLIP_{k\RL}(2,n,m_n)\\
=\sum_n\sum_{m_n} (\lambda ^*_{k/2}(1,n,m_n)\ e^{-2\pi in_{m_n}x_{N-k/2}} 
\times_D \lambda _{k/2}(1,n,m_n)\ e^{2\pi in_{m_n}x_{N-k/2}})
\end{multline*}
with:
\Bi
\item $ \lambda ^2_{k/2}(1,n,m_n)$ being eigenvalues of the coset representative
$(U_{n_{Rk}}\times U_{n_{Lk}})$ of the product of Hecke operators where $n^2_{N;k}$, denoting the number of transcendental quanta, is now written according to:
\[ n^2_{N;k}/d^2=a^2_{N;k}=a^{2k-2}_N\ \epsilon ^2(a^{2k-2}_N)\] 
referring to section~2.6.

So, $\lambda _{k/2}(1,n,m_n)$ corresponds to the radius of a semicircle equal to the respective coefficient
$\lambda (k/2,n,m_n)$ of $\ELLIP_L(k,n,m_n)$ by considering proposition~1.14.

\item $x_N$ being a real point of order $(N\times k/2)$ written now $x_{N-k/2}$.
\Ei
\vskip 11pt

\subsection{Proposition}

{\em The {\bbf $2$-dimensional elliptic bisemimodule 
$\ELLIP_{k\RL}(2,n,m_n)$\/} of weight $k$,\linebreak which is a real analytic cusp biform, 
{\bbf covers the corresponding cusp biform\linebreak
${}^{(\omega )}f_k(z^*_{N-k}) \times{}^{(\omega )}f_k(z_{N-k})$\/} according to:
\[
\ELLIP_{k\RL}(2,n,m_n) \quad \overset{\sim}{\hookrightarrow } \quad 
{}^{(\omega )}f_k(z^*_{N-k}) \times_D{}^{(\omega )}f_k(z_{N-k})\]
in the sense that:
\Bean
\item {\bbf every $2$-dimensional semitorus
\[
T^2_L[n]=\lambda _k(2,n)\ e^{2\pi inz_{N-k}}\quad 
\rresp{T^2_R[n]=\lambda^*_k(2,n)\ e^{2\pi inz^*_{N-k}}}
\]
 of class ``$n$'' of
$
{}^{(\omega )}f_k(z_{N-k})$
\resp{${}^{(\omega )}f_k(z^*_{N-k})$}
 is covered by $m_n=\sum\limits_{d\mid n} d\cdot N\cdot nu$ semicircles \begin{align*}
T^1_L[n,m_n]&=\lambda _{k/2}(1,n,m_n)\ e^{2\pi in_{m_n}x_{N-k/2}}\\[6pt]
\rresp{T^1_R[n,m_n]&=\lambda^* _{k/2}(1,n,m_n)\ e^{-2\pi in_{m_n}x_{N-k/2}}}
\end{align*}
at $(\tfrac k2\times a_{N;k})$ transcendental quanta of\/}
$\ELLIP_{k_L}(2,n,m_n)$
\resp{\linebreak $\ELLIP_{k_R}(2,n,m_n)$} localized in the upper \resp{lower} half plane, where $nu$ is the number of nonunits of Galois extensions;

\item {\bbf the parameter $b_N$ of the nilpotent group of matrices
$u_2(b_N)\cdot u_2(b_N)^t$ is a ``phase shift''\/} in the first dimension of the considered $(\zit\big/N \zit)^2$-bisemilattice.

When this phase shift ``$b_N$'' increases, the difference between the two eigenvalues
$\lambda _{k/2_+}(1,n,m_n)$ and
$\lambda _{k/2_-} (1,n,m_n)$ of $(U_{n_{R_k}}\times U_{n_{L_k}})$ also increases but verifies\linebreak
$\lambda ^2_{k/2_+}(1,n,m_n)\times \lambda ^2_{k/2_-}(1,n,m_n)=n^2_{N;k}$.
\Ee
}
\vskip 11pt

\bpr \Bean
\item As the $n$-th semitorus $T^2_L[n]=\lambda _k(2,n)\ e^{2\pi inz_{N-k}}$ of the cusp form ${}^{(\omega )}f_k(z_{N-k})$ is generated by the product of two orthogonal semicircles at $(k/2\times n)$ transcendental quanta whose product of radii is
$\lambda _{k/2}(2,n)$ which is an eigenvalue of 
$(U_{n_{R_k}}\times U_{n_{L_k}})$ with
$\alpha_{n^2_N}=\L( \BsM 1&0\\ 0&n^2_N\EsM\R)$ and
$n^2_{N;k}=\sum\limits_{d\mid n}a^2_N\cdot d^{2k-2}_N\ \epsilon ^2(d)$, we have that the number ``$m_n$'' of semicircles $T^1_L[n,m_n]=\lambda _{k/2}(1,n,m_n)\ e^{2\pi in_{m_n}x_{N-k/2}}$ covering $T^2_L[n[$
 is the integer 
 $d_{N\cdot nu}=\sum_{d\mid n} d^{(k-1)}\cdot N\cdot nu$ since, in this real case,
 $a^2_{N;k}=a^{2k-2}_N\cdot \epsilon ^2(a^{2k-2}_N)$.
 
 \item It is easy to calculate that when the phase shift
 $b_N=b\cdot N$ increases,
 $\L|\lambda _{k/2_+}(1,n,b_N)\R.\linebreak -\L.\lambda _{k/2_-}(1,n,b_N)\R|$ also increases, reflecting a deformation of one radius
 $\lambda _{k/2_+}(1,n,b_N)$ with respect to the order
 $\lambda _{k/2_-}(1,n,b_N)$.
 \epr
 \Ee
 \vskip 11pt
 
 \subsection{Proposition (Local curvature of a torus)}
 
 {\em 
Let \[
\{T^1_L(n,m_n)=\lambda _{k/2}(1,n,m_n)\ e^{2\pi in_{m_n}x_{N-k/2}}\}_{m_n=1}^{d\cdot N\cdot nu}\]
 be the set of semicircles covering the $n$-th two-dimensional semitorus
$T^2_L[n]\in {}^{(\omega )}f_k(z_{N-k})$.

Then, {\bbf the increase of the length of the semicircles with respect to $m_n$ depends on the twisting of ``opposite'' semicircles in function of their curvature leading to degenerate singularities\/} of fold type (or possible of cusp type) on them.

The {\bbf blowups of the versal deformations of these singularities in codimension $1$ and $2
$ are then sent on opposite semicircles increasing then their lengths.
}}
\vskip 11pt

\bpr Assume that the semitorus $T^2_L[n]$ is generated from a cylinder covered by $d\cdot N\cdot nu$ line segments having the same length, i.e. by ``$d$'' transcendental quanta.

Then, by bending this cylinder in order to get a semitorus, degenerate singularities of fold type are generated on the most bent line segments.

Let $f(x)=x^3$ be one of these fold singularities.

Its versal unfolding in codimension $1$ is
\[ F(x,a_1)=x^3+a_1\ x^1\;.\]
Possible degenerate singularities of cusp type $f(x)=x^4$ can be generated of which versal unfolding in codimension $2$ are
\[ F(x,a_1,a_2) = x^4+a_1\ x^1+a_2\ x^2\;.\]
Blowups of these versal deformations, consisting in the extensions of their quotient algebras, are introduced in \cite{Pie7}.  They are smooth endomorphisms based on Galois antiautomorphisms disconnecting the monomials $x^1\in F(x,a_1)$ and $x^1,x^2\in F(x,a_1,a_2)$ which are then sent on the opposite least bent line segments.\epr
\vskip 11pt

{\bbf This dynamical process is then responsible for the transition from the global euclidean geometry of the semitorus to local hyperbolic and spherical geometries\/} \cite{Pie12}.
\vskip 11pt

\subsection{Proposition}
{\em
\Bena
\item Let \[
T^{(2)}_{\zit_N}(H)=T_2(\rit)\big/T_2(\zit,N \zit)\]
 be the set of $\rit^2\big/(\zit/N \zit)^2$-sublattices of the Poincare upper half plane $H$.

Let $g^{(2)}_{\zit_N}[n]$ denote the $n$-th coset representative of $T^{(2)}_{\zit_N}(H)$ whose fundamental domain is $D^2_{nN,k}$ with respect to the sublattices $\Lambda ^{(2)}_{\zit_N}[n]$.

Let $f_k(z_{N-k})$ be a cusp form of weight $k$ and level $N$ which is:
\Be
\item invariant under a congruence subgroup $\Gamma  (N)$ of $SL_2(\zit)$;

\item periodic and expanded in Fourier series $f_k(z_{N-k})=\sum\limits_{n=0}^\infty c_{n_k}\ q^n_{N-k}$;

\item eigenfunction of Hecke operators $T_k(n;N)$ acting on the function
$F_k(T^{(2)}_{\zit_N}(H))$ on $T^{(2)}_{\zit_N}(H)$ according to a complete set of upper triangular coset representatives of $\Gamma (N)$.
\Ee

\item Let
\[ (T_1(F_\omega )=T_1(F_L)\Big / T_1(\zit\big/N \zit))\simeq 
(T_2(F_v )=T_2(F^+_L)\Big / T_2(\zit\big/N \zit))\]
be the set of $F_\omega \big/(\zit\big/N \zit)^2\simeq (F_v)\big/(\zit\big/N \zit)^2$-(sub)lattices of Galois on the\linebreak Poincare upper half plane $H$.

Let $D^{(2)}_{F_{\omega _n},N;k}$
\resp{$D^{(2)}_{F_{v _n},N;k}$} be the fundamental unitary domain of the $n$-th conjugacy class representative
$g^{(2)}_{\omega _L}[n]\in T_1(F_\omega )$
\resp{$g^{(2)}_{v _L}[n]\in T_2(F_v )$} with respect to the $(\zit\big/N \zit)^2$-sub(semi)lattice $\Lambda ^{(2)}_{\omega _n}$
\resp{$\Lambda ^{(2)}_{v _n}$}.

Let ${}^{(\omega )}f_k(z_{N-k})$ be a ``Weil'' cusp form of weight $k$ and level $N$ on $T_1(F_\omega )$ or on $T_2(F_v)$ which is:
\Be
\item invariant under a suitable congruence subgroup $\Gamma (N)$;

\item periodic and expanded in Fourier series
${}^{(\omega )}f_k(z_{N-k})=\sum\limits_n\lambda _k(2,n)\ e^{2\pi inz_{N-k}}$;

\item eigenfunction of the Hecke operator $T_{H-k}(N)$ according to section~2.6.
\Ee

{\bbf If the Weil unitary fundamental domain $D^{(2)}_{F_{\omega _n};N;k}$\/} (covered by
 $D^{(2)}_{F_{v _n};N;k}$)
 {\bbf covers the fundamental classical domain $D^{(2)}_{nN,k}$, then the Weil cusp form ${}^{(\omega )}f_k(z_{N-k})$ of weight $k$ and level $N$ can be identified with the classical cusp form
 $f_k(z_{N-k})$\/} of weight $k$ and level $N$ in such a way that:
 \[ {}^{(\omega )}f_k(z_{N-k}) \approx f_k(z_{N-k})\;, \]
 the weight $k$ referring to a geometric dimension while the level $N$ is an algebraic dimension related to the degree of extension of transcendental quanta: this constitutes the main objective of the Langlands program.
 \Ee
 }
 \vskip 11pt
 
 \bpr
 As the cusp form $f_k(z_{N-k})$ of weight $k$ originates from a $k$-dimensional cusp form
 $\phi ^{(k)}(z_{N_k})$ according to proposition~1.14, the two-dimensional fundamental domain $D^{(2)}_{nN;k}$ of weight $k$ is given by:
 \[
 D^{(2)}_{nN;k}=\L( D^{(2)}_{nN;k=2}\R)^{k/2}\]
 where $D^{(2)}_{nN;k=2}$ is the corresponding two-dimensional fundamental domain of weight $2$ expanded by a power $k/2$.
 \vskip 11pt
 
 Similarly, let $D^{(2)}_{F_{\omega _n};N;k=2}$ be the two-dimensional Weil unitary fundamental domain of weight $2$: its surface
 $S(D^{(2)}_{F_{\omega _n};N;k})$ is the square of the one transcendental quantum of level $N$;
 \[
S(D^{(2)}_{F_{\omega _n};N;k=2})=N^2\times (nu)^2\]
where $nu$ is the number of nonunits.

So, the two-dimensional unitary domain
$D^{(2)}_{F_{\omega _n};N;k}$ of weight $k$ is:
\[
D^{(2)}_{F_{\omega _n};N;k}= \L(
S(D^{(2)}_{F_{\omega _n};N;k=2})\R)^{k/2}= \L(N^2\times (nu)^2\R)^{k/2}\;.\]
\vskip 11pt

Now, assume by hypothesis that:
\[
D^{(2)}_{F_{\omega _n};N;k}\simeq D^{(2)}_{\omega _n;N;k}\;.\]
By periodicity of cusp forms, the fundamental unitary domains are the same in all coset representatives of $T^{(2)}_{\zit_N}(H)$ which is covered by $T_1(F_\omega )\simeq T_2(F_v)$.
\vskip 11pt

The unitary fundamental classical \resp{Galois or Weil} domain(s)
$D^{(2)}_{nN;k}$ \resp{\linebreak $D^{(2)}_{F_{\omega _n};N;k}$} generate(s), under the action of the Hecke operator having a representation in the upper triangular group of matrices
$T_2(\zit\big/N \zit))$ (or $\Gamma (N)$), the classes of automorphisms in
$T_2(\rit)$ \resp{$T_1(F_\omega )$ or in $T_2(F_v)$} according to the
$(\zit\big/N \zit)^2$-sublattices in the Poincare upper half plane $H$.

As the $(\zit\big/N \zit)^2$-lattice is the same in the classical and Weil (or Galois) case, the classical and Galois two-dimensional cusp forms of weight $k$ and level $N$ correspond, and, thus, we have:
\be
{}^{(\omega )}f_k(z_{N-k})\approx f_k(z_{N-k})\;.\tag*{\eop}
\ee
\vskip 11pt

\section[Zeta functions, Theta series, weak Maass forms, Mock Theta functions and the Tau function]{Zeta functions, Theta series, weak Maass forms,\linebreak Mock Theta functions and the Tau function}

This chapter will be devoted to a generalization of classical (or Weil) cusp forms in order to take into account Maass forms and Mock modular forms.

But, first, the generation of
\Bean
\item $L$-series from cusp forms \cite{Hec1} by a globally nonperiodic transform map,
\item theta series from modular forms with integer weights
\Ee
 will be introduced.
\vskip 11pt

\subsection{Proposition ($L$-functions as nonperiodic transforms of cusp forms)}

{\em Let 
\[
f_k(z_{N-k})=\sum\limits_{n=1}^\infty  c_{nk}\ q^n_{N-k}
\rresp{f_k(z^*_{N-k})=\sum\limits_{n=1}^\infty  c^*_{nk}\ q^{*n}_{N-k}}, \quad q^n_{N-k}=e^{2\pi inz_{N-k}}\;, \]
be a left (resp.  right dual) cusp form of weight $k$ and level $N$ as developed in sections~1.11 to 1.14 and 2.6 to 2.8.

Let
\[
s_+=\sigma +i\tau 
\rresp{s_-=\sigma -i\tau }\]
be a complex variable conjugate to $z=x+iy$
\resp{$z^*=x-iy$}.

Then, there exists a linear continuous map (Mellin transform):
\[
\Phi _L:f_k(z_{N-k})\To L(f_k,s_{N-k_+})
\rresp{\Phi _R:f^*_k(z^*_{N-k})\To L(f^*_k,s_{N-k_-})}
\]
in such a way that the {\bbf $L$-function $L(f_k,s_{N-k_+})$
\resp{$L(f^*_k,s_{N-k_-})$} be the nonperiodic transform of
$f_k(z_{N-k})$
\resp{$f^*_k(z_{N-k})$} given by\/}:
\begin{align*}
L(f_k,s_{N-k_+})
&= \sum_{n=1}^\infty \int f^{(n)}_k(z_{N-k})\ e^{-2\pi inz_{N-k}}\cdot n^{-s_{N-k_+}}\ dz_{N-k}\\[11pt]
&= \sum_{n=1}^\infty c_{n_k}\ n^{-s_{N-k_+}}\\[15pt]
\rresp{
L(f^*_k,s_{N-k_-})
&= \sum_{n=1}^\infty \int f^{(n)}_k(z^*_{N-k})\ e^{-2\pi inz^*_{N-k}}\cdot n^{-s_{N-k_-}}\ dz_{N-k}\\[11pt]
&= \sum_{n=1}^\infty c^*_{n_k}\ n^{-s_{N-k_-}}}
\end{align*}
and corresponding to the Mellin transform
\[
L(f_k,s_{N-k_+})=(2\Pi )^{s_+}\ \Gamma (s_+)^{-1}\int^\infty _0(-i\ z_{N-k})^{s_+}\ f_k(z_{N-k}\ \F{dz_{N-k}}{z_{N-k}}\]
where:
\Bi
\item $\begin{array}[t]{rl}
f^{(n)}_k(z_{N-k})&=c_{nk}\ q^n_{N-k}=c_{nk}\ e^{2\pi inz_{N-k}}\\
\rresp{f^{(n)}_k(z^*_{N-k})&=c^*_{nk}\ q^{*n}_{N-k}=c^*_{nk}\ e^{2\pi inz^*_{N-k}}}
\end{array}$

is the $n$-th term of $f_k(z_{N-k})$
\resp{$f_k(z^*_{N-k})$};

\item $z_{N-k}$
\resp{$z^*_{N-k}$} is a point of order $(N\times k)$ conjugate to the point 
$s_{N-k_+}$
\resp{$s_{N-k_-}$} of order $(N\times k)$.
\Ei
}
\vskip 11pt

\bpr
Referring to the literature of $L$-functions, we must admit that
$s_+$
\resp{$s_-$} is conjugate to $z$ \resp{$z^*$} in the sense that
$s_+\approx \tfrac 1z$
\resp{$s_-\approx \tfrac 1{z^*}$} is a  complex inverse space variable restricted to the upper \resp{lower} half plane.

$s_{N-k_+}$ 
\resp{$s_{N-k_-}$} is a complex point of order $(N\times k)$, usually called period $(N\times k)$.

$\phi _L$ \resp{$\phi _R$} is the nonperiodic equivalent  of the Fourier transform.  Indeed, if
$\phi _L$ \resp{$\phi _R$} was globally periodic,
$L(f_k,s_{N-k_+})$
\resp{$L(f^*_k,s_{N-k_-})$} would be given by:
\begin{align*}
L(f_k,s_{N-k_+}) &= \sum_{n=1}^\infty \int f^{(n)}_k(z_{N-k})\ e^{-2\pi inz_{N-k}}\ dz_{N-k}=\sum\limits_n c_{nk}\\[15pt]
\rresp{L(f^*_k,s_{N-k_-}) &= \sum_{n=1}^\infty \int f^{(n)}_k(z^*_{N-k})\ e^{-2\pi inz^*_{N-k}}\ dz_{N-k}=\sum\limits_n c^*_{nk}}\;.\end{align*}

{\bbf The factor $n^{-s_{N-k_+}}$
\resp{$n^{-s_{N-k_-}}$} is thus responsible for the global nonperiodicity of 
$\phi _L$ 
\resp{$\phi _R$}\/} which can be analyzed by the Poisson formula \cite{E-K}.

$f_k(z_N)$ is a periodic function of period
$T=D^{(2)}_{nk}$
\resp{$T_N=D^{(2)}_{nNk}$}, being a unitary fundamental domain according to section~1.3 and proposition~1.15, which can be expressed generally as \cite{Rod}:
\begin{align*}
f_k(z_{N-k}) &= f_k(z_{N-k})\star \sum_{n=1}^\infty \delta (z_{N-k}-nNkT)\\[11pt]
&= \sum_nf_k(z_{N-k}-nNkT)\;, \end{align*}
the period $T$ being explicitly written.

Its Fourier transform is:
\[
\dfrac1T \wh f_k(\nu _{N-k})\sum_n\delta \L( \nu _{N-k}-\dfrac{nNk}T\R) = \dfrac1T\sum_n\wh f\L(\dfrac{nNk}T\R)\delta \L(\nu _{N-k}-\dfrac{nNk}T\R)\;, \]
$\nu _{N-k}$ having the ``periodic'' variable conjugate to $z_{N-k}$, i.e.  a series of Dirac equidistant distributions with interval $\F1T$. 

By taking into account the inverse Fourier transform, we get:
\[
\sum_nf_k(z_N-nNkT)=\dfrac1T\sum_n\wh f_k\L(\dfrac{nNk}T\R)\ e^{2\pi inz_{N-k}}\;.\]
That can be summarized by the Fourier transform
\[ \sum_n\delta \L(z_{N-k}-nNkT\R) \quad
\overset{\rm FT}{\To} \quad \sum_n\delta \L(T\nu _{N-k}-nNk\R)\simeq \dfrac 1T \sum_n\L(\nu _{N-k}-\dfrac{nNk}T\R)\]
of Dirac distributions which is periodic \cite{Rod} while the nonperiodic transform $\Phi _L$ would lead to:
\[ 
\sum_n\delta \L(\nu _{N-k}-nNkT\R) \quad \overset{\Phi _L}{\To} \quad
\sum_n\dfrac1{T^{(n)}}\ \beta \L(s_{N-k_+}-\dfrac{nNk}{T^{(n)}}\R)
\]
where $\beta \L(s_{N-k_+}-\dfrac{nNk}{T^{(n)}}\R)$ would be the nonperiodic equivalent of the Dirac distribution
$\delta \L(\nu _{N-k}-\dfrac{nNk}T\R)$.

{\bbf The nonperiodicity of $\Phi _L$
\resp{$\Phi _R$}
can then be evaluated from the map of the unique period $T$  of $f_k(z_{N-k})$
\resp{$f_k(z^*_{N-k})$}, or more exactly $T_N$ associated with one space quantum, to the set
$\L\{1\big/T^{(n)}\R\}_{n=1}^\infty $ of inverse periods, i.e. unit periods, or more exactly
$\L\{1\big/T^{(n)}_{N-k}\R\}_{n=1}^\infty $ associated with the energies of one space quantum on the different packets ``$n$'' \cite{Pie8}\/}.\epr
\vskip 11pt

\subsection{Corollary}

{\em
If $k=2$ and $N=1$ (pseudounramified case), let the weight $2$ cusp biform $f_2(z^*_{N=1})\otimes_Df_2(z_{N=1})$ be reduced to the double  global elliptic bisemimodule
\[
\ELLIP\RL(1,n)=\sum_n\L(2\lambda ^{nr}(1,n)\ e^{-2\pi inx}\otimes_D
2\lambda ^{nr}(1,n)\ e^{2\pi inx}\R)\]
where, according to sections 2.7 and \cite{Pie1}:
\[
\lambda ^{nr}(1,n)=4n^2\;, \qquad x\in\rit\;.\]
Then, the nonperiodic map
\[
\Phi _R\otimes_D\Phi _L : \qquad \ELLIP\RL(1,n) \quad \To \quad \zeta _R(s_-)\otimes_D\zeta _L(s_+)\]
sends  the global elliptic bisemimodule
$\ELLIP\RL(1,n)$ into the (diagonal) product $\zeta _R(s_-)\linebreak \otimes_D\zeta _L(s_+)$ of the classical zeta functions:
\[
\zeta _R(s_-)  = \sum_nn^{-s_-}\;, \quad s_-= \sigma -i\tau \;, \qquad
\and \qquad \zeta _L(s_+) = \sum_nn^{-s_+}\;, \quad s_+= \sigma +i\tau \;.
\]
{\bbf The nonperiodicity of $\Phi _L$ \resp{$\Phi _R$} can then be evaluated from the Gaussian distribution of the consecutive spacings
$\delta \gamma _n=\gamma _{n+1}-\gamma _n$ between the nontrivial zeros of 
$\zeta _R(s_-)$, 
$\zeta _L(s_+)$ and
$\zeta (s)$\/}.
}
\vskip 11pt

\bpr
Referring to \cite{Pie8}, the kernel
$\Ker[\Phi _R\otimes_D\Phi _L]$ of 
$(\Phi _R\otimes_D\Phi _L)$ maps the product, right by left, of degeneracies of
$(2\lambda ^{nr}(1,n)\ e^{-2\pi inx}\otimes_D2\lambda ^{nr}(1,n)\ e^{2\pi inx})$ into the product ``$4n^2$'' of the trivial zeros of $\zeta _R(s_-)$ and 
 $\zeta _L(s_+)$. And, these products of trivial zeros are mapped into the products of the corresponding pairs of nontrivial zeros under the action of
 $D_{4n^2;i^2}\ \epsilon _{4n^2}$ which is a coset representative of the Lie bisemialgebra of the decomposition (bisemi)group restricted to $4n^2$: this corresponds to the solution of the Riemann hypothesis \cite{Pie1}.
 
 So, {\bbf the consecutive spacings
 $\delta \gamma _n=\gamma _{n+1}-\gamma _n$ between the nontrivial zeros of $\zeta  (s)$ are the energies of one free transcendental quantum in subsemilattices of $(n+1)$ quanta as proved in \cite{Pie8}.
 
 As the energies $\delta \gamma _n$ vary from one level ``$n$'' to another, for example ``$n+j$'' the zeta functions $\zeta _R(s_-)$ and $\zeta  _R(s_+)$ are non periodic\/}, which is also the case for the map $\Phi _R\otimes_D\Phi _L$.\epr
 \vskip 11pt
 
 \subsection{Theta series}
 
 Let 
$
 Q(r)=\langle n,n\rangle=\sum\limits_{i=1}^k\sum_na_{nn}^{(i)}\ n_i^2$
 be a $\zit$-valued positive definite quadratic form on $\zit^k$ reduced to its diagonal form which is always possible.
 
 As envisaged classically \cite{Ser}, the diagonal terms $a_{nn}^{(i)}$ of the matrix $A$ of rank $(t\times k)$, $1\le n\le t\le \infty $, $1\le i\le k$, are equal to $1$ and the $n_i\in\zit$ are the generators of $\zit$-subsemilattices
 $\{\Lambda _\zit^{(2)}[n]\}$ in the Poincare upper half plane.
 
 From the quadratic form $Q(n)$, {\bbf we introduce the theta series $\Theta (z_{N-k})$ by\/}:
 \begin{align*}
 \Theta (z_N)
 &= \sum_i \sum_n (q_N)^{\sum\limits_in_i^2}\;, && 1\le n\le t\;, \quad 1\le i\le k\;, \\[11pt]
 &= \sum_{i=1}^k\sum_ne^{2\pi i\L(\sum\limits_in_i^2z_N\R)}\;, && z_N \text{\ being a point of period $N$}\;, \\[11pt]
 &=\sum_ne^{2\pi in^2z_{N-k}}\;, && z_{N-k} \text{\ being a point of period $(N\times k)$}\;,
 \end{align*}
 or, more generally as suggested by D. Zagier in \cite{Zag}:
 \[
 \Theta _{a,b}(z_{N-k}) = \sum_{n+a}e^{i\pi i(n+a)^2z_{N-k}}\sum_{i=1}^k
e^{2\pi ib(n_i^2+a^2)}\;, \qquad a,b\in \qit^k\;.\]
\vskip 11pt

\subsection{Theta series as real analytic modular forms}

{\em
Let 
$
\Theta _{n^2k/2}(z_{N-k})=\sum\limits_nd_{nk/2}\ e^{2\pi in^2z_{N-k}}$
be a theta series introduced from the quadratic form
\[
Q(n) = \sum_{i=1}^k\sum_na^{(i)}_{nn}\ n^2_i\;.\]
Then, $\Theta _{n^2k/2}(z_{N-k})$ is a real analytic form of weight $k/2$, in one-to-one correspondence with the weight $k$ modular form
$f_k(z_{N-k})=\sum\limits_nc_{nk}\ q^n_{N-k}$.
}
\vskip 11pt

\bpr
Let
\[ f_k(z_{N-k})=\sum_{n=1}^\infty c_{nk}\ q^n_{N-k}\]
be a (left) cusp form of weight $k$ and level $N$ such that, in $q_{N-k}=e^{2\pi iz_{N-k}}$, $z_{N-k}$ is a point of order $(N\times k)$ according to proposition~1.15.

Assume that the theta series $\Theta _{n^2k/2}(z_{N-k})$ originates from the cusp form $f_k(z_{N-k})$ by the map:
\[
f\Theta _{kN}: \qquad f_k(z_{N-k})\quad \To \quad \Theta _{n^2k/2}(z_{N-k})\]
in such a way that {\bbf the two-dimensional semitori $c_{nk}\ q^n_{N-k}$, of which generators are two semicircles at $n$ transcendental quanta\/} (see proposition~1.9), {\bbf are sent into semicircles $d_{nk/2}\ e^{2\pi in^2x_{N-k/2}}$ at $n^2$ transcendental quanta by the bijective maps:
\[
T-c_{nk}: \qquad c_{nk}\ e^{2\pi inz_{N-k}} \quad \To \quad d_{nk/2}\ e^{2\pi in^2x_{N-k/2}}\;, \qquad \forall\ n\;, \]
at the condition that the complex point of order $(k\times N)$ be sent into the real point $x_{N-k/2}$ of order $\L(\dfrac k2\times N\R)$, i.e. that, in $z_{N-k}=x_{N-k/2}+iy_{N-k/2}$, $iy_{N-k/2}\to 0$\/} which is suggested in \cite{D-S}.

As $f_k(z_{N-k})$ is a modular form of weight $k$, generated from the $k$-dimensional geometric cusp form
$\Pi (\GL_{k/2}(F_\omega ))$ according to sections~2.7 and 2.3, and as the theta series
$\Theta _{n^2k/2}(z_{N-k})$ results from $f_k(z_{N-k})$ by the isomorphic maps
$\{T-c_{nk}\}_n$ stretching two-dimensional semitori into semicircles, $\Theta _{n^2k/2}(z_{N-k})$ is a modular form of weight $k/2$ since the maps
$T-c_{nk}$ divide the geometric dimension by two.\epr
\vskip 11pt

\subsection{Proposition}

{\em
Let 
$
f_k(z_{N-k})=\sum_nc_{nk}\ q^n_{N-k}$
be a cusp form of weight $k$ and level $N$.

Let
\[
\Theta _{n^2k/2}(z_{N-k/2})=\sum_nd_{nk/2}\ e^{2\pi in^2z_{N-k/2}}\]
be a theta series generated from $f_k(z_{N-k})$ by the bijective maps $\{T-c_{nk}\}_n$.

Let
\[
\ELLIP_{k/2}(2,n,m_n)=\sum_n\sum_{m_n}\lambda _{k/2}(1,n,m_n)\ e^{2\pi in_{m_n}x_{N-k/2}}\]
be the global elliptic semimodule covering the cusp form $f_k(z_{N-k})$.

Then, {\bbf we have the following commutative diagram\/}:
\[ \begin{psmatrix}[colsep=.5cm,rowsep=.5cm]
f_k(z_{N-k}) && 
\ELLIP_{k/2}(2,n,m_n) \\
& \;\; \raisebox{2mm}{\rotatebox{-30}{$\sim$}} \hspace{1.4cm} \rotatebox{30}{$\sim$} &\\
 & \Theta _{n^2k/2}(z_{N-k/2})&
\psset{arrows=->,nodesep=3pt}
\everypsbox{\scriptstyle}
\ncline{1,1}{1,3}^{f{\rm EL}_{kN}}
\ncline{1,1}{3,2}<{f\Theta _{kN}}
\ncline{3,2}{1,3}>{\Theta {\rm EL} _{kN}}
\end{psmatrix}
\]
in such a way that $\Theta {\rm EL}_{kN}=f\Theta _{kN}\circ f{\rm EL}^{-1}_{kN}$ be a bijective map.
}
\vskip 11pt

\bpr
The map
\[
\Theta {\rm EL}_{kN} : \qquad \Theta _{n^2k/2}(z_{N-k/2}) \quad \To \quad
\ELLIP_{k/2}(2,n,m_n)\]
sends semicircles of $\Theta _{n^2k/2}(z_{N-k/2})$ at $\L(\frac k2\times n\R)^2$ transcendental quanta into semicircles of $\ELLIP_{k/2}(2,n,m_n)$ at $\L(n\times \frac k2\R)$ transcendental quanta.

As $f {{\rm EL}_{kN}}$ and $f\Theta _{kN}$ are bijective maps,
$\Theta {\rm EL}_{kN}$ is clearly bijective, reorganizing the transcendental quanta of the theta series
$\Theta _{n^2k/2}(z_{N-k/2})$ according to the underlying $\zit\big/N \zit$-subsemilattices of 
$\Theta _{n^2k/2}(z_{N-k/2})$ and $\ELLIP_{k/2}(2,n,m_n)$:
\be \begin{psmatrix}[colsep=.5cm,rowsep=1cm]
\{\Lambda ^{(k/2)}_{\zit_N}[n^2]\}_n && &&
\{\Lambda ^{(k)}_{\zit_N}[n]\}_n   \\
\Theta _{n^2k/2}(z_{N-k/2})&&&&
 \ELLIP_{k/2}(2,n,m_n)
\psset{arrows=->,nodesep=3pt}
\everypsbox{\scriptstyle}
\ncline{1,1}{1,5}
\ncline{2,1}{2,5}^{\Theta {\rm EL}_{kN}}
\ncline{1,1}{2,1}
\ncline{1,5}{2,5}
\end{psmatrix}
\tag*{\eop}
\ee
\vskip 11pt

\subsection{Locally compact global elliptic semimodule}

Assume that the compact two-dimensional global elliptic semimodule
$\ELLIP_{1_L}(2,n,m_n)$ of weight $1$ and level $N$ covering the weight $2$ cusp form
$f_2(z_{N-2})$ is transformed by the map:
\[
C-Lc_{\rm ELL_2}: \qquad \ELLIP_{1_L}(2,n,m_n) \quad \To \quad \ELLIP^{(Lc)}_{1_L}(1,n,m_n)\]
into the locally compact one-dimensional global elliptic semimodule $\ELLIP^{(Lc)}_{1_L}(1,n,m_n)$ in such a way that the $m_n$- semicircles
$\{\lambda _{1_L}(2,n,m_n)\ e^{2\pi in_{m_n}x_{N-2}}\}_{m_n}$ of
$\ELLIP_{1_L}(2,n,m_n)$, covering compactly the $n$-th semitorus of $f_2(z_{N-2})$, are no more connected.

Thus, $\ELLIP^{(Lc)}_{1_L}(1,n,m_n)$ becomes one-dimensional and its weight, related to the geometric dimension, can be considered as being $1\cdot \ELLIP^{(Lc)}_{1_L}(1,n,m_n)$ is then a modular form of weight $1$ and will be written according to
$\ELLIP^{(Lc)}_{1_L}(1,n,m_n)$.

So, we have the commutative diagram:
\[ \begin{psmatrix}[colsep=.5cm,rowsep=.5cm]
f_2(z_{N-2}) && &&  \ELLIP_{1_L}(2,n,m_n)&&  &&
 \ELLIP^{(Lc)}_{1_L}(1,n,m_n) \\
&& && \;\;\; \raisebox{2mm}{\rotatebox{-30}{$\sim$}} \hspace{3.5cm} \raisebox{1mm}{\rotatebox{30}{$\sim$}} && &&\\
 && && \Theta _{n^21/2}(z_{N-1/2})&& &&
\psset{arrows=->,nodesep=3pt}
\ncline{1,5}{1,1}^{\sim}
\ncline{1,5}{1,9}^{\sim}
\ncline{1,1}{3,5}
\ncline{3,5}{1,5}>{\wr}
\ncline{1,9}{3,5}
\end{psmatrix}
\]
\vskip 11pt

\subsection{Proposition (The Theta series of weight $1/2$)}

{\em
\[ \Theta _{n^21/2}(z_{N-1/2})=\sum_nd_{n1/2}\ e^{2\pi in^2x_{N-1/2}}\]
is generated either from the weight $2$ modular form $f_2(z_{N-2})$ or from the weight $1$ real analytic modular form $\ELLIP^{(Lc)}_{1_L}(1,n,m_n)$.
}
\vskip 11pt

\bpr
This results directly from the commutative diagram of section~3.6 by taking into account that the locally compact modular form of weight $1$:
$\ELLIP^{(Lc)}_{1_L}(1,n,m_n)$ covers the weight $2$ cusp form $f_2(z_{N-2})$.\epr
\vskip 11pt

\subsection{Remark (Classical approach of modular forms of weight $1/2$)}

This way of generating modular forms of weight $1/2$ generalizes the classical approach of G. Shimura \cite{Shi} and J.P. Serre \cite{Ser} in which it is proved that every modular form of weight $1/2$ on $\Gamma _1(N)$ is a linear combination of theta series with characters.
\vskip 11pt

\subsection{Weak Maass forms: a quick review}

Maass forms \cite{Maa} of weight $2$, now called generally weak Maass forms, are functions $f_2^{WM}(z)$, on the Poincare upper half plane, equipped with the Riemannian metric 
$ds^2=\frac{dx^2+dy^2}{y^2}$ and being modular under the action of $SL(2,\zit)$ implying that
\[
f_2^{WM} \L(\frac{az+b}{cz+d}\R) =f^{WM}_2(z)\;, \qquad \forall\ \L(\BsM a & b \\ c & d\EsM\R)\in SL(2,\zit)\;.\]
$F^{WM}_2(z)$ has a Fourier series expansion of the form \cite{Boo}:
\[ f^{WM}_2(z) = \sum_{n=1}^\infty \lambda (n)\ \sqrt y\ K_{ir}(2\pi ny)\ \cos(2\pi nx)\;, \qquad z=x+iy\;, \]
where $K_{ir}(2\pi ny)$ is the classical Bessel function, and is an eigenfunction of the hyperbolic Laplace operator $\Delta =y^2\L(\F \partial {\partial x^2}+\F \partial {\partial y^2}\R)$ having eigenvalue $\F14+r^2$.

The $L$-function associated with $f^{WM}_2(z)$ is the series $L(f^{WM}_2,s)=\sum\limits_{n=1}^\infty \lambda (n)\ n^{-s}$ converging for $\Re (s)>1$.
\vskip 11pt

A {\bbf weak Maass form of weight $k$ on a subgroup $\Gamma \subset \Gamma _0(4)$ \/} verifies \cite{Ono1}:
\[
f^{WM}_k(Az)=\L(\F cd\R)^{2k}\ \epsilon ^{-2k}_d(cz+d)^k\ f^{WM}_k(z)
\qquad \text{for\ } A=\L(\BsM a&  b \\ c & d\EsM\R)\in\Gamma \]
where $\L(\F cd\R)$ is the extended Legendre symbol

 and {\bbf has a Fourier expansion of the type\/}:
\[ f^{WM}_k(z) = \sum_{n=n_0}^\infty  \gamma (f,n;y)\ q^{-n}+\sum_{n=n_1}^\infty  a(f,n)\ q^n\;, \qquad q=e^{2\pi iz}\;, \]
where $\gamma (f,n;y)$ are functions in $y$ and $a(f,n)$ are complex numbers.

{\bbf $\sum\limits_{n=n_0}^\infty \gamma (f,n;y)\ q^n$ refers to the nonholomorphic part of $f^{WM}_k(z)$ while\linebreak $\sum\limits_{n=n_1}^\infty a(f,n)\ q^n$ is its holomorphic part\/}.
\vskip 11pt

As no explicit examples of the Maass forms are known at this day \cite{Boo}, {\bbf a geometric interpretation of weak Maass forms will be given to them in terms of modular curves\/}.

A first elementary fundamental remark is that, in the Fourier series expansion $f^{WM}_2(z)$ of a weak Maass form, {\bbf the ``imaginary'' dimension ``$y$'' is lowered with respect to the dimension ``$x$''\/} in $z=x+iy$.

Taking into account that the Fourier series expansion of a modular curve decomposes into two-dimensional semitori $\{T^2_n(a_{nm|d^2},d)\}_{n=1}^t$, $t\le \infty $, generated by two orthogonal semicircles referring to propositions~1.9, 1.10 and 2.8, we are led to the following proposition.
\vskip 11pt

\subsection{Proposition (Geometric interpretation of weak Maass forms)}

{\em
{\bbf The series expansion $f^{WM}_2(z)$ of a weak Maass form consists of a sum of two-dimensional semitori $\{T^{2,(e\ell)}_n(S^1_{a_n}\times e\ell^1_{ef,n})\}^t_{n=1}$, $t\le \infty $, with elliptic cross sections $e\ell^1_{ef,n}$\/} in such a way that:
\Bi
\item the elliptic section $e\ell^1_{ef,n}$ is a bijective deformation of the semicircle $S^1_{d_n}$ generator of the $n$-th semitorus $T^2_n(a_{nm|d^2},d)\simeq S^1_{a_n}\times S^1_{d_n}$ of the modular curve (cusp form);

\item the equation of the ellipse $e\ell^1_{ef,n}$ corresponding to the circle $S^1_{d_n}$ having radius $r_{S^1_{d_n}}$ is given by:
\[
z_{e\ell\ell^1_n}=r_{S^1_{d_n}}(f\cos 2\pi (iy)+ie\sin2\pi (iy))\]
where $e$ and $f$ are the half lengths of the axis of the ellipse $e\ell^1_{ef,n}$.
\Ei
}
\vskip 11pt

\bpr
As the imaginary dimension ``$y$'' is lowered with respect to the dimension ``$x$'', {\bbf the semicircular sections $S^1_{d_n}$ of the two-dimensional semitori $T^2_n(a_{nm|d^2},d)$ are transformed into semielliptic sections $e\ell^1_{df,n}$ according to the bijective map\/}:
\[ Se: \qquad S^1_{d_n}\quad \To \quad e\ell^1_{cf,n}\;, \qquad \forall\ n\in\nit\;,\]
in such a way that:
\Bena
\item the length (or area) of $S^1_{d_n}$ is equal to the length of $e\ell^1_{ef,n}$;

\item the number of transcendental quanta on $S^1_{d_n}$ is equal to that on $e\ell^1_{ef,n}$.
\Ee
\vskip 11pt

{\bbf The canonical equation
\[ f^2\ X^2+e^2\ Y^2=f^2\ e^2\]
of the ellipse $e\ell^1_{ef,n}$,\/} having the same length as the circle $S^1_{a_n}$, {\bbf is transformed under the bijective map ``$Se$'' into\/}:
\[f^2\ r^2_{S^1_{d_n}}\cos^22\pi iy+e^2\ r^2_{S^1_{d_n}}\sin^22\pi iy=e^2\ f^2\]
or
\[ r^2_{S^1_{d_n}}(f^2\cos^22\pi iy+e^2\sin^22\pi iy)=e^2\ f^2\]
or
{\bbf 
\[z_{e\ell\ell^1_n}=r_{S^1_{d_n}}(f\cos2\pi iy+i\ e\sin2\pi iy)\]
}
by the change of variables $\begin{cases}
X &= r_{S^1_{d_n}}\cos2\pi iy\;, \\
Y &= r_{S^1_{d_n}}\sin2\pi iy\;, \end{cases}$

where:
\Bi
\item the ``imaginary angle'' $iy$ refers to a circle localized in a plane orthogonal to the circle $S^1_{a_n}$ \cite{Pie2};

\item $e$ and $f$ are the half lengths of the axis of the ellipse $e\ell\ell^1_{ef,n}$.\epr
\Ei
\vskip 11pt

\subsection{Proposition (Holomorphic and nonholomorphic part of the new series expansion of a weak Maass form)}

{\em
Let
\begin{align*}
f^{WM}_2(z) &= \sum_n\L[ \L( r_{S^1_{d_n}} ( f\cos2\pi n(iy)+ie\sin2\pi n(iy))\R) \times
\L( r_{S^1_{a_n}}\  e^{2\pi inx} \R) \R] \\[11pt]
&= \sum_nT^{2,(e\ell)}_n ( S^1_{a_n},e\ell^1_{ef,n})\end{align*}
be {\bbf the new series expansion of the weak Maass form $f^{WM}_2(z)$ of weight $2$.

If
\[ f^{WM,hol}_2(z) = \sum_nc(n)\ e^{2\pi inz}\;, \]
where $c(n)=r_{S^1_{a_n}}\times r_{S^1_{d_n,inscr}}$ denotes the holomorphic part of the ``Fourier'' series expansion of $f^{WM}_2(z)$\/} in such a way that
\[ S^1_{d_n,inscr}=r_{S^1_{d_n,inscr}}\ e^{2\pi in(iy)}\]
be the equation of the circle inscribed in the elliptic $e\ell ^1_{ef,n}$,

then, {\bbf the nonholomorphic part $f^{WM,nonhol}_2(z)$ of the weak Maass form
$f^{WM}_2(z)$, is given by\/}:
\[ f^{WM,nonhol}_2(z) = \sum_nT^{2,(e\ell)}_n(S^1_{a_n},e\ell^1_{ef,n}) -
\sum_n \L( r_{S^1_{a_n}}\times r_{S^1_{d_n,inscr}}\R)\ e^{2\pi inz}\]
and {\bbf corresponds to the shadow of $f^{WM}_2(z)$.
}}
\vskip 11pt

\bpr
As the new series expansion of the weak Maass form $f^{WM}_2(z)$ is the sum, over the integers $n\in\nit$, of the products of the equation of the circle $S^1_{a_n}$ over ``$x$'' by the equation of the ellipse $e\ell^1_{ef,n}$ over ``$i\cdot y$'', which is a bijection with the equation of the circle having as radius $r_{S^1_{d_n}}$, it appears clearly that the holomorphic part 
$f^{WM,hol}_2(z)$ of $f^{WM}_2(z)$ must be given by
\[
f^{WM,hol}_2(z)
 = \sum_nr_{S^1_{a_n}}\ e^{2\pi inx}\times r_{S^1_{d_n,incr}}\ e^{2\pi i n(iy)}
 = \sum_nc(n)\ e^{2\pi inz}\]
where $r_{S^1_{d_n,incr}}\ e^{2\pi in(iy)}$ is the equation of the circle inscribed in the ellipse $e\ell^1_{ef,n}$.

So, the {\bbf nonholomorphic part $f^{WM,nonhol}_2(z)$ of $f^{WM}_2(z)$ will be given by the series expansion\/}:\pagebreak
\begin{align*}
& f^{WM,nonhol}_2(z)\\[11pt]
&\qquad= \sum_n T^{2,(e\ell)}_n(S^1_{a_n},e\ell^1_{ef,n})-\sum_nc(n)\ e^{2\pi inz}\\[11pt]
&\qquad= \sum_n\L[ S^1_{a_n}\times (e\ell^1_{ef,n}-S^1_{d_n,incr})\R]\\[11pt]
&\qquad= \sum_n\L\{
r_{S^1_{a_n}}\ e^{2\pi inx}\times\L[\L( r_{S^1_{d_n}}(f\cos2\pi iy+ie\sin2\pi iy)\R)
-\L( r_{S^1_{d_n,incr}}\ e^{2\pi in(iy)}\R)\R]\R\}
\end{align*}
which clearly corresponds to the nonholomorphic part
$ \sum\limits_{n=n_0}^\infty \gamma (f,n;y)\ q^{-n}$
of the Fourier expansion of the weak Maass form $f^{WM}_2(z)$ as developed in section~3.9.

This nonholomorphic part
$f^{WM,nonhol}_2(z)$ is then the ``shadow'' of the weak Maass form
$f^{WM}_2(z)$ with respect to its holomorphic part
$f^{WM,hol}_2(z)$ in the sense that the shadow introduced for describing the Mock modular forms \cite{Zag} is a unary theta series of weight $3/2$, i.e. a function of the form
$\sum\limits_n\varepsilon (n)\ n \ q^{\kappa n^2}$ with $\kappa \in\qit$ and $\varepsilon $ and odd periodic function.  Taking into account propositions~3.4 and 3.5, it can be shown that
$f^{WM,nonhol}_2(z)$ can be transformed in such a unary theta series of weight $3/2$.\epr
\vskip 11pt

\subsection{Weak Maass forms of weight $k$ and level $N$}

Referring to proposition~1.14 and section~2.7, we infer that {\bbf the weak Maass form
$f^{WM}_2(z)$ of weight $2$ and level $1$ can be extended to the weak Maas form
$f^{WM}_k(z_{N-k})$ of weight $k$ and level $N$ defined by the series\/}:
\begin{align*}
& f^{WM}_k(z_{N-k})\\[11pt]
&\qquad = \sum_nT^{2,(e\ell)}_n(S^1_{a_{n_{N-k/2}}},e\ell^1_{ef,n_{N-k/2}})\\[11pt]
&\qquad = \sum_n\L[ \L( r_{S^1_{d_{n-k/2}}}(f\cos2\pi n(iy)_{N-k/2})+ie\sin2\pi n(iy)_{N/2})\R)
\times r_{S^1_{a_{n-k/2}}}\ e^{2\pi inx_{N-k/2}} \R]\;,
\end{align*}
i.e. {\bbf by the sum over the integers $n$ of two-dimensional semitori
$T^{2,(e\ell)}_n(S^1_{a_{n_{N-k/2}}},\linebreak e\ell^1_{ef,n_{N-k/2}})$ with elliptic cross sections
$e\ell^1_{ef,n_{N-k/2}}$ of weight $k/2$ and level $N$\/} where:
\Bi
\item $S^1_{a_{n_{N-k/2}}}=r_{S^1_{a_{n-k/2}}}\ e^{2\pi inx_{N-k/2}}$ is the equation of a semicircle over ``$x$'' at $(a\times k/2)$ transcendental quanta with radius
$r_{S^1_{a_{n-k/2}}}$ and real points $x_{N-k/2}$ of order $(N\times k/2)$, $n=ad$;

\item $e\ell^1_{ef,n_{N-k/2}}=r_{S^1_{d_{n-k/2}}}(f\cos2\pi n(iy_{N-k/2})+ie\sin2\pi n(iy_{N-k/2}))$ is the equation of an ellipse over ``$i\cdot y$'' at $(d\times k/2)$ transcendental quanta with
$r_{S^1_{d_{n-k/2}}}$ being the radius of the circle in one-to-one correspondence with
$e\ell^1_{ef,n_{N-k/2}}$ and ``real'' points ``$iy_{N-k/2}$'' of order $(N\times k/2)$.
\Ei
\vskip 11pt

\subsection{Proposition}

{\em
{\bbf The weak Maass form $f^{WM}_k(z_{N-k})$ of weight $k$ and level $N$ defined by the series:
\[
f^{WM}_k(z_{N-k})=\sum_nT^{2,(e\ell)}_n(S^1_{a_{n_{N-k/2}}},e\ell^1_{ef,n_{N-k/}})
\]
is:
\Bean
\item periodic;
\item holomorphic and weakly modular.
\Ee
}}
\vskip 11pt

\bpr
\Bean
\item $f^{WM}_k(z_{N-k})$ is periodic because it verifies:
\[
f^{WM}_k(z_{N-k})=f^{WM}_k(z_{N-k}+1)\;, \]

\item {\bbf $f^{WM}_k(z_{N-k})$ is holomorphic but not modular with respect to the elliptic cross section $e\ell^1_{ef,n_{N-k/2}}$ deviation from circularity\/}. On the other hand, the circle factor $S^1_{a_{n_{N-k/2}}}$ of $T^{2,(e\ell)}_n$ is modular.

Consequently, $f^{WM}_k(z_{N-k})$ will be said to be weakly modular, the modularity given by:
\[ f^{WM}_k\L(\frac {az_{N-k}+b}{cz_{N-k}+d}\R)=
f^{WM}_k(z_{N-k})\qquad \text{for\ } \L(\BsM a & b \\ c & d\EsM\R)\in SL(2,\zit)\;.\]
Remark that the modularity can be checked by using the formula:
\[
\cos \L(a\cdot b\R)=\half \L(\exp (iab)+\exp (-iab)\R)
=\half \L( (\exp ia)^b+(\exp-ia)^b\R)\]
and
\be
\cos \L(\tfrac ab\R)={\txt\half}\ \L( (\exp (ia)^{1/b}+(\exp -ia)^{1/b}\R)\;.\tag*{\eop}
\ee
\Ee
\vskip 11pt

\subsection{Generalized weak Maass forms: Elliptic forms}

{\bbf
The weak Maass form $f^{WM}_k(z_{N-k})$ of weight $k$ and level $N$\/}, having a decomposition into the sum over the $n$ sublattices, $1\le n\le\infty $, of two-dimensional semitori
$T^{2,(e\ell)}_n$ with elliptic cross section, {\bbf can be generalized to elliptic forms of weight $k$ and level $N$ introduced as having a decomposition into the sum over $n$ sublattices of surfaces of revolution of (semi)ellipses rotating around ellipses instead of circles\/} as for the weak Maass forms.

Thus, {\bbf a two-dimensional elliptic form
$f^{EL}_k(z_{N-k})$ of weight $k$ and level $N$ will be defined by the series\/}:
\begin{align*}
f^{EL}_k(z_{N-k}) 
&  = \sum_n EL^{2,(e\ell)}_n(e\ell^1_{ab_{n_{N-k/2}}},e\ell^1_{ef,n_{N-k/2}})\\[11pt]
& = \sum_n \L[ \L( r_{S^1_{a_{n-k/2}}}(b\cos2\pi nx_{N-k/2}+ia\sin2\pi nx_{N-k/2})\R)  \R.\\[11pt]
& \qquad  \qquad \times\L.\L( r_{S^1_{d_{n-k/2}}}(f\cos2\pi n(iy_{N-k/2})+ie\sin2\pi n(iy_{N-k/2}))\R) \R]
\end{align*}
where:
\Bi
\item  $EL^{2,(e\ell)}_n(e\ell^1_{ab_{n_{N-k/2}}},e\ell^1_{ef,n_{N-k/2}})$ is the $n$-th surface of revolution of the ellipse\linebreak $e\ell^1_{ef,n_{N-k/2}}$ rotating around the ellipse $e\ell_{ab_{n_{N-k/2}}}$;
\item the equation of the ellipse $e\ell^1_{ab_{n_{N-k/2}}}$ is given by:
\[
e\ell_{ab_{n_{N-k/2}}}= r_{S^1_{a_{n-k/2}}}(b\cos2\pi nx_{N-k/2}+ia\sin2\pi nx_{N-k/2})\]
 with $r_{S^1_{a_{n-k/2}}}$ being the radius of the circle in bijection with
 $e\ell^1_{ab_{n_{N-k/2}}}$ referring to section~3.12 and $x_{N-k/2}$ being a point or order $(N\times k/2)$.
 \Ei
 
 {\bbf The elliptic form $f^{EL}_k(z_{N-k})$, defined over a doubly periodic lattice, is thus an elliptic function characterized by
 four parameters ``$a$, $b$, $e$, $f$'' which are the half lengths of the axis of the two ellipses
 $e\ell^1_{ab_{n_{N-k/2}}}$ and
 $e\ell^1_{ef,n_{N-k/2}}$\/}.
 \vskip 11pt
 
 \subsection{Proposition}
 
 {\em
 Let $f_k(z_{N-k})$, $f^{WM}_k(z_{N-k})$ and $f^{EL}_k(z_{N-k})$ be respectively a cusp form, a weak Maass form and an elliptic form of weight $k$ and level $N$.
 
 Then {\bbf the commutative diagram\/}
 \[ \begin{psmatrix}[colsep=.5cm,rowsep=.5cm]
f_k(z_{N-k}) && 
f^{WM}_k(z_{N-k})  \\
& \;\; \raisebox{2mm}{\rotatebox{-30}{$\sim$}} \hspace{1.4cm} \rotatebox{30}{$\sim$} &\\
 & f^{EL}_k(z_{N-k})&
\psset{arrows=->,nodesep=3pt}
\ncline{1,1}{1,3}^{\sim}
\ncline{3,2}{1,1}
\ncline{1,3}{3,2}
\end{psmatrix}
\]
{\bbf indicates the possible transformation of one of these forms into another\/}.
}
\vskip 11pt

\bpr The above commutative diagram results from the one-to-one correspondence between these forms referring to proposition~3.10.\epr
\vskip 11pt

\subsection{Mock Theta functions of Ramanujan: A brief synthesis}

Each of the 17 Ramanujan's theta functions belongs to one of the three families of functions:

\Bean
\item Lerch sums,
\item quotients of binary theta series by unary theta series,
\item Fourier coefficients of Jacobi forms,
\Ee
as discovered by S. Zwegers \cite{Zwe}.

These theta functions are $q$-hypergeometric series of the form $\sum\limits_{n=0}^\infty A_n(q)$, $A_n(q)\in \qit(q)$.  The first of these is:
\[ f(q) = \sum_{n=0}^\infty \frac {q^{n^2}}{(1+q^2)\dots (1+q^n)^2}\;.\]

{\bbf Each Ramanujan theta function is a Mock theta function given by the $q$-series
$ H(q) = \sum_{n>0} a_n \ q^n$ in such a way that $q^\lambda \ H(q)$, $\lambda \in \QQ$, be a Mock modular form of weight $1/2$ whose shadow is a unary theta series of weight $3/2$\/} given by a function of the form $\sum\limits_n\varepsilon (n)\ n\ q^{Kn^2}$, $K\in\QQ$, where $\varepsilon (n)$ is an odd periodic function.

{\bbf A Mock theta function is thus a Mock modular form of weight $k$ of the space $\MM_k$ of such forms   extending the space $M_k$ 
of classical modular forms\/} of weight $k$ and characterized {\bbf
by a shadow $g=S[h]$\/} which is a modular form of weight $(2-k)$.

Zagier then deduced that {\bbf the Ramanujan's Mock theta function $H(q)$ acquire modularity after carrying out the following three steps \cite{Zag}\/}:
\Bean
\item multiply $H(q)$ by a rational power $q^\lambda $ of $q$;
\item change the variable $q=e^{2\pi i\tau }$ by $\tau $, setting $h(\tau )=e^{2\pi i\lambda \tau }\ H(e^{2\pi i\tau })$;

\item add a nonholomorphic correction term $g^*(\tau )$ to $h(\tau )$ in such a way that $\widehat h(\tau )=h(\tau )+g^*(\tau )$ transforms like a modular form of weight $1/2$ with $g^*(\tau )$ associated to a theta series $g(\tau )=\sum\limits_{n\in\zit+d}n\ q^{Kn^2}$ being a modular form of weight $3/2$ and the shadow of $h(\tau )$.
\Ee

Finally, it appears that the space $\MM_k$ of Mock theta functions is in one-to-one correspondence with the space $\wt \MM_k$ of weak Maas forms.

Having proposed in the preceding sections a geometric and algebraic approach to theta series and weak Maass forms, {\bbf we will try to give in the next proposition a geometric interpretation of Ramanujan's theta functions transformed into Mock theta functions related to weak Maass forms\/}.
\vskip 11pt

\subsection[Proposition (Connection of Ramanujan Theta functions with weak Maass forms, modular forms and global elliptic semimodules)]{Proposition (Connection of Ramanujan Theta functions\linebreak  with weak Maass forms, modular forms and global elliptic semimodules)}

{\em
Let 
$\Ms_{H\to \wh h}:  H(q)  \To  \wh h(\tau )$
denote the map corresponding to the three step sequence recalled in section~3.16 and transforming each Ramanujan theta function $H(q)$ into a Mock modular form $\wh h(\tau )$ of weight $1/2$ characterized by a shadow $g^*(\tau )$ being a modular form of weight $3/2$.

Then, the map $\Ms_{H\to\wh h}$ corresponds to the composition of maps
$(\Ms_{f^{WM;hol}_2\to \wh h}\circ\linebreak \Ms_{H\to f_2^{WM;hol}})$ of {\bbf the commutative diagram\/}:
 \[ \scalebox{.9}{$\begin{psmatrix}[colsep=.5cm,rowsep=.5cm]
 && \ELLIP^{WM,hol}_{1_L}(2,n,m_n) &&   \ELLIP_{1_L}(2,n,m_n) &&&& f_2(\tau _{12})\\
& 
&&
\rotatebox{30}{$\begin{array}{l} \scriptstyle{\Ms_{f^{WM,hol}\to\ELLIP}}\\ {}\end{array}$}\\
H(q) && f_2^{WM,hol}(\tau_{1-2})
&& \wh h(\tau ) &&
&& f_2^M(\tau _{1-2})\\
&& &  \rotatebox{30}{${\renewcommand{\arraystretch}{1.2}\begin{array}{c} \sim\\ \Ms_{f_2^{WM}\to\wh h} \end{array}}$} \\
&& f^{WM}_2(\tau _{1-2})
\psset{arrows=->,nodesep=3pt}
\everypsbox{\scriptstyle}
\ncline{1,3}{1,5}
\ncline{1,9}{1,5}^{\phantom{\qquad\;\;}\Ms_{\ELLIP\to f_2}}_{\phantom{\qquad}\ds\sim}
\ncline{3,3}{1,3}<{\Ms_{f^{WM,hol}_2\to\ELLIP^{WM,hol}}}>{\ds\wr}
\ncline{3,5}{1,5}>{\Ms_{\wh h\to\ELLIP}}<{\ds\wr}
\ncline{3,9}{1,9}<{\Ms_{f^M_2\to f_2}}>{\ds\wr}
\ncline{3,1}{3,3}^{\!\!\Ms_{H\to f^{WM,hol}_2}\;\;}
\ncline{3,3}{3,5}^{\Ms_{f^{WM,hol}_2\to\wh h}}
\ncline{3,5}{3,9}^{\Ms_{\wh h\to f^M_2}}_{\ds\sim}
\ncline{3,3}{5,3}<{\Ms_{f_2^{WM,hol}\to f_2^{WM}}}
\ncline{5,3}{3,5}
\ncline{3,3}{1,5}
\end{psmatrix}$}
\]
where:
\Bi
\item $f^{WM,hol}_2(\tau _{1-2})$ is the holomorphic part of the weak Maass form $f^{WM}_2(\tau _{1-2})$ of the space $\wh \MM_2$ of weight $2$ and level $1$;

\item $f^M_2(\tau _{1-2})$ is the Mock modular form of weight $2$ and level $1$ of the space $\MM_2$ extending the space $SL(2,1)$ of classical (or Weil) modular forms $f_2(\tau _{1-2})$ of weight $2$ and level $1$;
\item $\ELLIP_{1_L}(2,n,m_n)$ and 
  $\ELLIP^{WM,hol}_{1_L}(2,n,m_n)$ are the global elliptic semimodules of weight $1$ and level $1$ covering respectively the modular forms $f_2(\tau _{1-2})$ and\linebreak $f_2^{WM,hol}(\tau _{1-2})$;
  
  \item $\Ms_{H\to f_2^{WM,hol}}$ and $\Ms_{f_2^{WM,hol}\to\wh h}$ are the maps respectively of the steps ((a) -- b) and c)) of the three-step sequence recalled in section~3.16;
  
  \item $\Ms_{\wh h\to f^M_2}$ is the covering map of the Mock modular form $\wh h(\tau )$ of weight $1/2$ and level $1$ by the Mock modular form $f^M_2(\tau _{1-2})$ of weight $2$ and level $1$,
  \Ei
  
  in such a way that:
 \Bean
 \item {\bbf the nonholomorphic part $f^{WM,nonhol}_2(\tau _{1-2})$ 
of the weak Maass form $f^{WM}_2(\tau _{1-2})$
corresponds to the shadow $g^*(\tau )$ of the Mock modular form $\wh h(\tau )$\/};
 
 \item the spaces $\MM_1$ of Mock theta functions $\wh h(\tau )$ and $|\wh M|_2$ of weak Maass forms $f_2^{WM}(\tau _{1-2})$ are isomorphic;
 
 \item the bijective map $\Ms_{\wh h\to\ELLIP}: \wh h(\tau )\overset{\sim}{\To}
 \ELLIP_{1_L}(2,n,m_n)$ is an extension of the bijective map:
 \[
 \Ms_{f^{WM,hol}_2\to\ELLIP}: \qquad f_2^{WM,hol}(\tau _{1-2}) \quad \overset{\sim}{\To} \quad
 \ELLIP_{1_L}^{WM,hol}(2,n,m_n)\]
 by the shadow $\wh h(\tau )$.
 \Ee
 }
 \vskip 11pt
 
 \bpr
 The map
 \[ \Ms_{H\to f_2^{WM,hol}} : \qquad H(q) \quad \To \quad f_2^{WM,hol}(\tau _{1-2})\;, \]
 corresponding to the steps a) and b) of section~3.16, transforms each Ramanujan's $q$-hypergeometric series $H(q)=\sum\limits_na_n\ q^n$ into a holomorphic function which may be the holomorphic part $f_2^{WM,hol}(\tau _{1-2})$ of a weak Maass form $f_2^{WM}(\tau _{1-2})$ of weight $2$ and level $1$,  referring to proposition~3.10.
 
 Then, the map
 \[ 
 \Ms_{f_2^{WM,hol}\to\wh h}: \qquad f_2^{WM,hol}(\tau _{1-2}) \quad \To \quad \wh h(\tau )\]
 corresponds to the map
 \[
 \Ms_{f_2^{WM,hol}\to f_2^{WM}}: \qquad f_2^{WM,hol}(\tau _{1-2}) \quad \To \quad f_2^{WM}(\tau _{1-2})\]
 extending $f_2^{WM,hol}(\tau _{1-2})$ by the shadow $g^*(\tau )$ which may be the nonholomorphic part
 $f_2^{WM,nonhol}(\tau _{1-2})$ of the weak Maass form $f_2^{WM}(\tau _{1-2})
 $ according to proposition~3.11.
 
 It then results that the global elliptic semimodule $\ELLIP_{1_L}(2,n,m_n)$, covering $\wh h(t\tau )$, is an extension by the shadow $g^*(\tau )$ of the global elliptic semimodule $\ELLIP^{WM,hol}_{1_L}(2,n,m_n)$ covering $f_2^{WM,hol}(\tau _{1-2})$.
 
 And, thus, the spaces $\MM_1$ of Mock theta functions $\wh h(\tau )$, which are real analytic modular forms \cite{Wat}, and $\wh \MM_2$ of weak Maass forms $f_2^{MW}(\tau _{1-2})$ are isomorphic (see proposition~3.7).\epr
 \vskip 11pt
 
\subsection{Corollary (Mock modular forms of weight $k$)}

{\em
The weight $2$ case of Mock modular forms studied in proposition~3.17 can be easily generalized to the weight $k$ by considering that the shadow of weight $(2-k)$ of a Mock modular form of weight $k$ corresponds to the nonholomorphic part of a weak Maass form $f_2^{MW}(\tau_{1-k})$ of weight $k$ and level $1$.
}
\vskip 11pt

\subsection{Partitions of $n$ and Dyson's rank}

In order to provide a combinatorial explanation of Ramanujan's congruence for the number of partitions $p(n)$ of an integer $n$ \cite{B-O1} of which generating function is 
\[ \sum_{n=0}^\infty p(n)\ q^n = \prod_{n=1}^\infty \frac 1{1-q^n}\;, \]
{\bbf F. Dyson \cite{Dys} introduced the rank of a partition defined to be its largest part minus the number of its parts\/}.

In this respect, let $N(n,m)$ denote the number of partitions of $n$ having rank $m$ congruent to $r\mod s$.

{\bbf The generating function giving the number of partitions of $n$ with rank $m$ is\/}:
\begin{align*}
R(\omega ,q)
&= \sum_{n=1}^\infty \sum_{m=-\infty }^{+\infty } N(n,m)\ \omega ^m\ q^n\\
&= 1+\sum_{n=1}^\infty \frac{q^{n^2}}{\prod\limits_{m=1}^n(1-\omega \ q^m)(1-\omega ^{-1}\ q^m)}
\end{align*}
where $\omega =e^{2\pi ix/s}$ is a $s$-th root of unity.

According to D. Zagier \cite{Zag}, knowing the functions $n\to N(n,m)$ for all $r\ (\mod s)$ is equivalent to knowing the specializations of $R(\omega ;q)$ to all $s$-th root of unity $\omega $.  For $\omega =-1$, $R(\omega ,q)$ specializes to the first Ramanujan's Mock theta function $f(q)$ given in section~3.16.  Bringmann and Ono \cite{B-O2} generalize this to other roots of unity.
\vskip 11pt

\subsection[Proposition (Partitions of quanta in Ramanujan's Theta functions)]{Proposition (Partitions of quanta in Ramanujan's Theta\linebreak functions)}

{\em
Let $R(\omega ,q)$ be the partition function specializing to the (17) Ramanujan's Mock theta functions $H(q)$ of weight $1/2$ and level $1$.

As $N(n,m)$ denotes the number of partitions of $n^{(2)}$ transcendental quanta, it follows that:
\Bena
\item {\bbf the Dyson's rank $m$ of a partition of $n$ must be the order of the maximal Galois sub(bisemi)group associated with the considered transcendental biextension minus the number of Galois sub(bisemi)groups\/};

\item {\bbf $\omega^m$ is a phase factor\/} related to $N(n,m)$;

\item {\bbf $N(n,m)\times nu$ may be the multiplicity of the $n$-th global elliptic subsemimodule $\ELLIP_{1_L}(1,n)$ covering the $n$-th term of the Ramanujan's Mock theta function $H(q)$\/} specialized from $R(\omega ;q)$ where $nu$ is the number of nonunits of Galois extensions;

\item there is a map:
\[
\Ms_{h\to f_{(2)}}: \qquad H(q) \quad \To \quad f_2(\tau _{1-2})\]
from the theta function $H(q)$ to the corresponding Hecke cusp form $f_2(\tau _{1-2})$ of weight $2$ and level $1$ in such a way that {\bbf $N(n,m)\ \omega ^m$ maps into the $n$-th coefficient $c_{n2}$\/} (being a global Hecke character) {\bbf of the cusp form $f_2(\tau _{1-2})=\sum\limits_n c_{n2}\cdot q^n_{1-2}$\/}.
\Ee
}
\vskip 11pt

\bpr
\Bena
\item As the integer $n$ in the Mock theta series refers to the $n$-th $\zit^2$-sublattice (see proposition~1.7) of the cusp form $f_2(\tau _{1-2})$ extending $H(q)$, as developed in proposition~3.17, and as $f_2(\tau _{1-2})$ is generated from the Weil (or Galois) group $W^{ab}_{F_\omega }$ according to section~2.2 by means of a Langlands global correspondence (see section~2.5), it (i.e. the integer $n$) must correspond to the order of a Weil (or Galois) subgroup.

Taking into account the definition of the rank $m$ of a partition of $n$, it is clear that the {\bbf largest part of partition is the order of the maximal Galois sub(bisemi)\-group associated with the considered partition and that the number of parts of the partition is the number of intermediate Galois (or Weil) sub(bisemi)\-groups\/}.

\item Referring to the diagram of proposition~3.17, we see that the global elliptic
semimodule $\ELLIP_{1_L}^{WM,hol}(2,n,m_n)$ covers the holomorphic part $f_2^{WM,hol}(\tau _{1-2})$ of the weak Maass form $f_2^{WM}(\tau _{1-2})$
 and thus also $H(q)$.
 
 Then, $(N(n,m)\times nu)$ must refer to the multiplicity ``$m_n$'' of the $n$-th global elliptic subsemimodule
 \[ \ellip_{1_L}(1,n)=\sum_{m_n}\lambda _1(1\text{\ or\ }2,n)\ e^{2\pi in_{m_n}x_{1-1}}\]
 according to section~2.7.
 \Ee
 
 \dots\ in contrast with the result of proposition~2.8 referring to the multiplicity with respect to Hecke cusp forms in the sense of:
 
 \Bena
 \item[3)] the map
 \[ \Ms_{H\to f(2)} : \qquad \begin{aligned}[t]
 H(q) \quad &\To \quad f_2(\tau _{1-2})\\
N(n,m)\ \omega ^m \quad & \To \quad c_{n2}\;, && \qquad \forall\ n\in\nit\;, \end{aligned}\]
is such that $N(n,m)\ \omega ^m$ be sent into the $n$-th coefficient of $f_2(\tau _{1-2})$ which is a product of radii of two orthogonal circles (see proposition~1.9).
\epr
\Ee
\vskip 11pt

\subsection{The Tau function: A brief summary}

The Ramanujan tau function is the function defined by:
\[
\sum_{n=1}^\infty \tau (n)\ q^n=q\ \prod_{n=1}^\infty (1-q^n)^{24}=\Delta (z)\;, \qquad q=e^{2\pi iz}\,;, \quad z=x+iy\;, \quad y>0\;.\]
As we have the equality
\[ \Delta (z) = \eta (z)^{24}\;, \]
where $\eta (z)=q^{1/24}\prod\limits_{n=1}^\infty (1-q_n)$ is the Dedekind's eta function of weight $1/2$ and level $1$ due to transformation laws:
\[
\eta (-1/z) = (-iz)^{1/2}\ \eta (z) \qquad \and \qquad n(z+1)=e^{\pi i/12}\ n(z)\;, \]
and as $\Delta (z)$ satisfies the symmetry condition:
\[ \Delta (-1/z) = z^{12}\ \Delta (z)\;, \]
{\bbf $\Delta (z)$ is a modular form of weight $12$ and level $1$\/}.

The coefficients $\tau (n)$ satisfy:
\Bean
\item $\tau (nm)=\tau (n)\ \tau (m)$ if $(n,m)=1$,

\item $\tau (p^{n+1})=\tau (p)\ \tau (p^n)-p^{11}\ \tau (p^{n-1})$, $n>1$,

\item $|\tau (p)|\le 2\ p^{11/2}$, $\forall\ $ prime 
\Ee
which is the Ramanujan's conjecture proved by P. Deligne.

{\bbf The $\tau (n)$ enjoy many congruence relations\/}, for example:
\[
\tau (n)=\sigma _{11}(n)\mod 2^{11}\qquad \text{for\ } n=1\mod 8\;, \]
and a very interesting one is:
\[
\tau  (n) =\sum_{0<d|n} d^{11}\mod 691\]
in such a way that the Eisenstein series
\[
E_{12}(q) = -B_{12}/24 + \sum_{n=1}^\infty \L(\sum_{0<d|n} d^{11}\R)\ q^n\;, \]
where
\[
B_{12}/24 = -691/65520 = 0\mod 691\;, \]
has the same Fourier expansion coefficients as $\Delta (z)$ as analyzed by B. Mazur in \cite{Maz}.

The Dirichlet series associated with the Ramanujan tau function is:
\[
L_\tau (s) = \prod_p \F1{1-\tau (p)\ p^{-s}+p^{11-2s}}\]
and the function $(2\pi )^{-s}\ \Gamma (s)\ L_\tau (s)$ is invariant under the substitution $s\to 12-s$.

The Ramanujan tau function $\Delta (z)=\sum\limits_n\tau (n)\ q^n=q\ \prod\limits_n(1-q^n)^{24}$ is a generating function on numbers of inverse partitions of quanta (corresponding perhaps to partitions in the dual (semi)space with respect to the generating function
\[
\sum_np(n)\ q^n =\prod _n \F1{1-q^n}\]
of numbers of partitions $p(n)$ of quanta as analyzed in section~3.19.

Referring to proposition~3.20, {\bbf the map
\[
\Ms_{\Delta \to f_{(12)}} : \qquad \Delta (z) \quad \To \quad f_2(z_{12-1})\]
from the tau function $\Delta (z)$ of weight $12$ to the associated two-dimensional (Hecke) Weil cusp form $f_2(z_{12-1})$ of weight $12$ and level $1$ will be studied in the next proposition and proved to correspond to a map from this tau function to a cusp form $\phi ^{(12)}(z_{12-1})$ in $12$ real dimensions and level $1$\/} (on $\cit^6$) {\bf which is an orthogonal universal structure\/}.  This may be a clue to know how to compute the congruence relations of $\tau (n)$.
\vskip 11pt

\subsection{Proposition (Universal (orthogonal) cusp form in dimension $2$ and weight $12$)}

{\em
Let $f_2(z_{12-1})$ be a two-dimensional Weil cusp form of weight $12$ and level $1$.

Let
\[ CP_{(12)\to2}: \qquad \phi ^{(12)}(z_{12-1})\quad \To \quad f_2(z_{12-1})\]
be the map from the cusp form $\phi ^{(12)}(z_{12-1})$ on $\cit^6$ of level $1$ projecting it into the two-dimensional cusp form $f_2(z_{12-1})$ and satisfying the commutative diagram:
 \[ \begin{psmatrix}[colsep=.5cm,rowsep=.3cm]
\phi ^{(12)}(z_{12-1}) && 
f_2(z_{12-1})  \\
& \;\; \raisebox{1cm}{\rotatebox{-45}{$\scriptstyle{\Ms^{-1}_{\Delta \to\phi ^{(12)}}}$}} \hspace{.5cm} \rotatebox{45}{$\scriptstyle{\Ms_{\Delta \to f(12)}}$} &\\
 & \Delta (z)
\psset{arrows=->,nodesep=3pt}
\everypsbox{\scriptstyle}
\ncline{1,1}{1,3}^{CP_{(12)\to 2}}
\ncline{1,1}{3,2}
\ncline{3,2}{1,3}
\end{psmatrix}
\]
Then, {\bbf the two-dimensional cusp form $f_2(z_{12-1})$ of weight $12$ and level $1$ is a universal ``orthogonal'' cusp form corresponding throughout Langlands global correspondences to the sum of the cuspidal representations of six bilinear algebraic semigroups generating three two-dimensional embedded toric bisemisheaves as well as their orthogonal equivalents\/} according to:
\[
 \Pi ^{(12)}(\GL_6(\wt F_{\o\omega }\times_D \wt F_\omega ))=\bigoplus_{i=1}^6
\Pi ^{(2_i)}(\GL_{1_i}(\wt F_{\o \omega }\times_D \wt F_\omega ))
\quad  
 \xrightarrow{CP^{(12)\to 2}_{R\times L\to L}}
\quad
f_2(z_{12-1})
\]
where:
\Bi
\item $\Pi ^{(2i)}(\GL_{1_i}(\wt F_{\o \omega }\times_D \wt F_\omega ))$, being the two-dimensional cuspidal representation of the bilinear algebraic semigroup
$\GL_{1_i}(\wt F_{\o \omega }\times_D \wt F_\omega )$, is in one-to-one correspondence with the two-dimensional cusp form $f_{2_i}(z_{2-1})$ of weight $2$ and level $1$, $1\le i\le 6$;

\item $\Pi ^{(12)}(\GL_6(\wt F_{\o\omega }\times_D \wt F_\omega ))$, being the twelve-dimensional cuspidal representation of\linebreak
$\GL_6(\wt F_{\o\omega }\times_D \wt F_\omega )$, is in one-to-one correspondence with the cusp form
$\phi ^{(12)}(z_{12-1})$.
\Ei

This is summarized in the commutative diagram:
\[ \begin{psmatrix}[colsep=.5cm,rowsep=1cm]
\{\Pi ^{(2i)}(\GL_{1_i}(\wt F_{\o \omega }\times_D \wt F_\omega ))\}^6_{i=1} && &&
\Pi ^{(12)} (\GL_6 (\wt F_{\o\omega }\times_D \wt F_\omega ) )\\
&& &&\phi ^{(12)}(z_{12-1})\\
\{ f_{2_i}(z_{2-1})\}^6_{i=1} && &&f_2(z_{12-1})
\psset{arrows=->,nodesep=3pt}
\everypsbox{\scriptstyle}
\ncline{1,1}{1,5}
\ncline{1,1}{3,1}
\ncline{3,1}{3,5}
\ncline{1,5}{2,5}
\ncline{2,5}{3,5}>{CP_{(12)\to 2}}
\end{psmatrix}
\]
}
\vskip 11pt

\bpr
\Bena
\item Let $\{\phi _R(g^{(2)}_{\o \omega _R}[n])\otimes
\phi _L(g^{(2)}_{\omega _L}[n])\}^{t\le \infty }_{n=1}$ denote the set of two-dimensional bisections of a bisemisheaf
$(M^{(2)}_R(F_{\o\omega })\otimes M^{(2)}_L(F_\omega ))$ of differentiable bifunctions constituting the functional representation space
$\FREPSP(\GL_1(\wt F_{\o\omega }\times \wt F_\omega ))$ of the bilinear algebraic semigroup
$\GL_1(F_{\o\omega }\times F_\omega )$ over sets of complex transcendental extensions $F_{\o\omega }$ and $F_\omega $.

Let $(M^{(2)\perp}_R(F_{\o\omega })\otimes M^{(2)\perp}_L(F_\omega ))$ be the orthogonal complement bisemisheaf of\linebreak
$(M^{(2)}_R(F_{\o\omega })\otimes M^{(2)}_L(F_\omega ))$.

Degenerate singularities on these bisemisheaves can give rise, by versal deformations and blowups of these, to one or two new covering bisemisheaves according to the kind of considered singularities as developed in \cite{Pie9}, \cite{Pie7}.

If there are degenerate singularities of corank $1$ and codimension $3$ on these bisemisheaves, we get, after a process of versal deformations, blowups of these, desingularizations and toroidal compactifications, {\bbf the three shell embedded bisemisheaves\/}:
\begin{multline*}
\Pi ^{(2_1)}(\GL_{1_1}(F_{\o\omega }\times F_\omega ))\oplus\Pi ^{(2_2)}(\GL_{1_2}(F_{\o\omega }\times F_\omega )^\perp)\\[6pt]
\subset
\Pi ^{(2_3)}(\GL_{1_3}(F_{\o\omega_{cov(1)} }\times F_{\omega_{cov(1)}} ))\oplus\Pi ^{(2_4)}(\GL_{1_4}(F_{\o\omega_{cov(1)} }\times F_{\omega_{cov(1)}} )^\perp)\\[6pt]
\subset
\Pi ^{(2_5)}(\GL_{1_5}(F_{\o\omega_{cov(2)} }\times F_{\omega_{cov(2)}} ))\oplus\Pi ^{(2_6)}(\GL_{1_6}(F_{\o\omega_{cov(2)} }\times F_{\omega_{cov(2)}} )^\perp)
\end{multline*}
where:
\Bi
\item $\Pi ^{(2_i)}(\GL_{1_i}(- \times  -))$ is the cuspidal representation of the $i$-th bilinear (algebraic) semigroup
$\GL_{1_i}(- \times  -))$, $1\le i\le 6$ in such a way that the sum of its conjugacy class representatives, which are products, right by left, of two-dimensional semitori according to sections~2.6 to 2.8, is the two-dimensional cusp biform
$f_{2_i}(z_{2-1})\times f_{2_i}(z_{2-1})$ of weight $2$ and level $1$ referring to sections~2.5 and 1.12;

\item $\Pi ^{(2_4)}(\GL_{1_4}(F_{\o\omega_{cov(1)} }\times F_{\omega_{cov(1)}} ))$ is the cuspidal representation of the blinear (algebraic) semigroup constituting the first shell covering of $\Pi ^{(2_2}(\GL_{1_2}(F_{\o\omega }\times F_\omega ))$ with ``$\perp$'' referring to the orthogonal complement cuspidal representation.
\Ei

\item Taking into account the Langlands global functoriality conjecture introduced in \cite{Pie4}, we have that the sum
$\sum\limits_{i=1}^6 \Pi ^{(2_i)}(\GL_{1_i}(F_{\o\omega }\times F_\omega ))$ of the six above cuspidal representations is equal to the $12$-dimensional cuspidal representation $\Pi ^{(12)}(\GL_6(F_{\o\omega }\times_DF_\omega ))$ of the bilinear (abstract) complex semigroup
 $\GL_6(F_{\o\omega }\times_DF_\omega )$.
 
 Indeed, the Langlands global functoriality conjecture states that {\bbf the $12$-\linebreak dimensional cuspidal representation
 $\Pi ^{(12)}(\GL_6(F_{\o\omega }\times_DF_\omega ))$ is nonorthogonally reducible it if decomposes diagonally according to the direct sum of irreducible cuspidal representations of the (algebraic) bilinear semigroups
 $\GL_{1_i}(F_{\o\omega }\times F_\omega )$ and offdiagonally according to the direct sum\/}
 \[
\bigoplus_{i\neq j=1}^6(\Pi ^{(2_i)}(\GL_{1_i}(F_{\o\omega }))\times
\Pi ^{(2_j)}(\GL_{1_j}(F_\omega )))\]
{\bbf of the (tensor) products of irreducible cuspidal representations of cross (algebraic) linear semigroups\/}
\[
\GL_{1_i}(F_{\o\omega }))\times
(\GL_{1_j}(F_\omega )) \equiv T^t_{1_i}(F_{\o\omega })\times T_{1_j}(F_\omega )\;.\]
So,
\begin{multline*}
\Pi ^{(12)}(\GL_6(F_{\o\omega }\times F_\omega ))
= \bigoplus_{i=1}^6 \Pi ^{(2_i)}(\GL_{1_i}(F_{\o\omega }\times F_\omega ))\\[6pt]
 \bigoplus_{i\neq j=1}^6 \Pi ^{(2_i)}(\GL_{1_i}(F_{\o\omega }))\otimes \Pi ^{(2_j)}(\GL_{1_j}( F_\omega ))
 \end{multline*}
 in such a way that, if 
 $\Pi ^{(12)}(\GL_6(F_{\o\omega }\times F_\omega ))$ is orthogonally completely reducible, then the crossed cuspidal representations
 $ \bigoplus_{i\neq j=1}^6 \Pi ^{(2_i)}(\GL_{1_i}(F_{\o\omega }))\otimes \Pi ^{(2_j)}(\GL_{1_j}( F_\omega )$ are equal to zero.
 
 \item Taking into account the sum of the $12$-dimensional cuspidal conjugacy class representatives of
$\Pi ^{(12)}(\GL_6(F_{\o\omega }\times F_\omega ))$, we get {\bbf the $12$-dimensional cusp biform
$\phi ^{(12)}(z^*_{12-1})\times \phi ^{(12)}(z_{12-1})$ of level $1$ which, by the map
$CP_{(12)\to2}$, sends it into the two-dimensional cusp biform
$f_2(z^*_{12-1})\times f_2(z_{12-1})$ of weight $12$ and level $1$\/} referring to propositions~1.14 and 1.17:
\[ \begin{psmatrix}[colsep=.5cm,rowsep=1cm]
\Pi ^{(12)}(\GL_{6}(F_{\o \omega }\times_D  F_\omega )) 
&& &&\phi ^{(12)}(z^*_{12-1})\times _D\phi ^{(12)}(z_{12-1})\\
&&&& f_{2}(z^*_{12-1})\times_D f_2(z_{12-1})\\
&&&& f_2(z_{12-1})
\psset{arrows=->,nodesep=3pt}
\everypsbox{\scriptstyle}
\ncline{1,1}{1,5}
\ncline{1,5}{2,5}>{CP_{(12)\to2}}
\ncline{2,5}{3,5}
\ncline{1,1}{3,5}
\end{psmatrix}
\]
As, by hypothesis, $f_{2}(z^*_{12-1})\times_D f_2(z_{12-1})$ and
$f_{2_i}(z^*_{2-1})\times_D f_{2_i}(z_{2-1})$ are diagonal cusp biforms, they are in one-to-one correspondence respectively with their left equivalents
$f_2(z_{12-1})$ and $f_{2_i}(z_{2-1})$.

Then, {\bbf we get the equality\/}:
\[
f_2(z_{12-1})=\bigoplus_{i=1}^6 f_{2_i}(z_{2-1})\]
resulting from the above commutative diagram and {\bbf corresponding to the announced universal ``orthogonal'' cusp form $f_2(z_{12-1})$ in dimension $2$, weight $12$ and level $1$\/}.

\item Finally, the map
\[
\Ms^{-1}_{\Delta \to\phi ^{(12)}}: \qquad
\phi ^{(12)}(z_{12-1}) \quad \To \quad \Delta (z)\]
directly results from the maps
\[
CP_{(12)\to 2}: \quad \phi ^{(12)}(z_{12-1}) \To f_2(z_{12-1}) \quad \and \quad
\Ms_{\Delta \to f(12)}: \quad \Delta (z) \To f_2(z_{12-1})\;.\]
Then, {\bbf $\Delta (z)$\/}, defined by
\[ \Delta (z) = (\eta (z)^4)^6=(q^{1/6}\prod_n(1-q^n)^4)^6\;, \]
i.e. by the $6$-th power of $(\eta (z))^4$ which is a modular form of weight $2$,
{\bbf is directly related to a $12$-dimensional cusp form $\phi ^{(12)}(z_{12-1})$ of level $1$\/}.\epr
\Ee
\vskip 11pt

\subsection{Universal nonorthogonal bilinear cuspidal representations}

The generalization of proposition~3.22 to universal ``nonorthogonal'' cuspidal representations including crossed cuspidal representations of interaction will now be envisaged.

Referring to the Langlands functoriality conjecture \cite{Pie4}, {\bbf the bilinear cuspidal representation
$\Pi ^{(12)}(\GL_6(F_{\o\omega }\times F_\omega ))$ of the (algebraic) bilinear semigroup
$\GL_6(F_{\o\omega }\times F_\omega )$ decomposes essentially according to\/}:
\begin{multline*}
\qquad \Pi ^{(12)}_{\rm rel}(\GL_6(F_{\o\omega }\times F_\omega ))
= \bigoplus_{i=1}^6 \Pi ^{(2_i)}(\GL_{1_i}(F_{\o\omega }\times_D F_\omega ))\\[6pt]
 \bigoplus_{i\neq j=1}^6 \Pi ^{(2_i)}(\GL_{1_i}(F_{\o\omega }))\otimes_{OD} \Pi ^{(2_j)}(\GL_{1_j}(F_{\omega }))\qquad 
 \end{multline*}
 where the second sum on the right hand side refers to crossed cuspidal representations leading to a nonorthogonally reducible representation of
$\Pi ^{(12)}(\GL_6(F_{\o\omega }\times F_\omega ))$.

According to \cite{Pie9}, {\bbf the six diagonal cuspidal representations
$\Pi ^{(2_i)}(\GL_{1_i}(F_{\o\omega }\times_{D} F_{\omega }))$, $1\le i\le 6$, correspond to an embedded three-shell universal physical structure of space-time (``$ST$''), middle-ground (``$MG$'') and mass (``$M$'')\/}:
\begin{multline*}
\Pi ^{(2_1)}(\GL_{1_1}^{ST}(F_{\o v }\times_D F_v ))_T\oplus\Pi ^{(2_2)}(\GL_{1_1}^{ST}(F_{\o\omega }\times_D F_\omega ))_S\\[6pt]
\subset
\Pi ^{(1_3)}(\GL_{1_3}^{MG}(F_{\o v }\times_D F_v ))_T\oplus\Pi ^{(2_4)}(\GL_{1_4}^{MG}(F_{\o\omega }\times_D F_\omega ))_S\\[6pt]
\subset
\Pi ^{(1_5)}(\GL_{1_5}^{M}(F_{\o v }\times_D F_v ))_T\oplus\Pi ^{(2_6)}(\GL_{1_6}^{M}(F_{\o\omega }\times_D F_\omega ))_S
\end{multline*}
in such a way that:
\Bi
\item $\Pi ^{(1_1)}(\GL_{1_1}^{ST}(F_{\o v }\times_D F_v ))_T\approx
\Pi ^{(2_1)}(\GL_{1_1}(F_{\o \omega  }\times_D F_\omega  ))$,

$\Pi ^{(1_3)}(\GL_{1_3}^{MG}(F_{\o v }\times_D F_v ))_T\approx
\Pi ^{(2_3)}(\GL_{1_3}(F_{\o \omega  }\times_D F_\omega  ))$,

and $\Pi ^{(1_5)}(\GL_{1_5}^{M}(F_{\o v }\times_D F_v ))_T\approx
\Pi ^{(2_5)}(\GL_{1_5}(F_{\o \omega  }\times_D F_\omega  ))$,

are {\bbf the one-dimensional ``time'' diagonal cuspidal representations respectively of the space-time, middle-ground and mass shells\/};

\item $\Pi ^{(2_2)}(\GL_{1_2}^{ST}(F_{\o\omega }\times_D F_\omega ))_S$, $\Pi ^{(2_4)}(\GL_{1_4}^{MG}(F_{\o\omega }\times_D F_\omega ))_S$ and 
$\Pi ^{(2_6)}(\GL_{1_6}^{M}(F_{\o\omega }\times_D F_\omega ))_S$ are similarly {\bbf the two-dimensional ``space'' diagonal cuspidal representations respectively of the space-time, middle-ground and mass shells\/}.
\Ei

In this context, {\bbf the six relevant off diagonal or crossed cuspidal representations
$\Pi ^{(2_i)}(\GL_{1_i}(F_{\o\omega }))\otimes_{OD} \Pi ^{(2_j)}(\GL_{1_j}(F_{\omega }))$, being interaction crossed cuspidal representations
$\Pi ^{(2_i)}(\GL_{2_i}(F_{\o v }))\otimes_{int} \Pi ^{(2_j)}(\GL_{2_j}(F_{v }))$, are for \cite{Pie10}\/}:
\Bi
\item {\bbf $i\neq j= 1$, $2$ and $3$ the crossed cuspidal representations of the interacting fields respectively of the space-time, middle-ground and mass shells\/};

\item {\bbf $i\neq j=4$, $5$ and $6$ the interacting crossed cuspidal representations between respectively the three different right and left semifields
of the ``$ST_R-MG_L$'', ``$MG_R-M_L$'' and ``$ST_R-M_L$'',  mixed shells\/}:

$\Pi ^{(2_4)}(\GL_{2_4}^{ST}(F_{\o v }))_T\otimes_{int}\Pi ^{(2_5)}(\GL_{2_5}^{MG}(F_{v }))_S$,
$\Pi ^{(2_5)}(\GL_{2_5}^{MG}(F_{\o v }))_T\otimes_{int}\Pi ^{(2_6)}(\GL_{2_6}^{M}(F_{v }))_S$ 

and
$\Pi ^{(2_4)}(\GL_{2_4}^{ST}(F_{\o v }))_T\otimes_{int}\Pi ^{(2_6)}(\GL_{2_6}^{M}(F_{v }))_S$.
\Ei
\vskip 11pt

\subsection{Proposition}

{\em
The cross binary product
\[
\mathop{\times}\limits_{\o{i=1}}^6 (f_{2_i}(z^*_{2-1})\times f_{2_i}(z_{2-1}))
=\L( \sum_{i=1}^6 f_{2_i}(z^*_{2-1})\R)\times\L( \sum_{i=1}^6 f_{2_i}(z_{2-1})\R)
\]
between the six cusp biforms
$f_{2_i}(z^*_{2-1})\times f_{2_i}(z_{2-1})$ of dimension $2$, weight $2$ and level $1$ is the image of the injective map:
\[
\Ms_{\Pi ^{(12)}\to\mathop{\u\times}\limits^6}: \qquad
\Pi ^{(12)}_{\rm rel}(\GL_6(F_{\o\omega }\times F_\omega )) \quad \To \quad
\mathop{\times}\limits_{\o{i=1}}^6 (f_{2_i}(z^*_{2-1})\times f_{2_i}(z_{2-1}))
\]
from {\bbf the relevant universal nonorthogonal bilinear cuspidal representation\linebreak
$\Pi ^{(12)}_{\rm rel}(\GL_6(F_{\o\omega }\times F_\omega ))$ decomposing  into six relevant cuspidal representations of interaction\/} by taking into account the section~3.23.
}
\vskip 11pt

\bpr
The cross binary product between cusp biforms was introduced in \cite{Pie3} and recalled in section~1.12.

The six diagonal cusp biforms $\{f_{2_i}(z^*_{2-1})\times_D f_{2_i}(z_{2-1}))\}^6_{i=1}$ are in one-to-one correspondence with the six time and space diagonal cuspidal representations of the embedded three shell universal physical structures $ST$, $MG$ and $M$ mentioned in section~3.23.

The thirty off diagonal bilinear cuspidal representations
$\{f_{2_i}(z^*_{2-1})\times_{OD} f_{2_j}(z_{2-1}))\}^6_{\substack{i\neq j=1\\i>j}}$
reduce by symmetry, $i> j$, to fifteen ones among which the six in one-to-one correspondence with the interacting crossed cuspidal representations taken into account in section~3.23 are relevant.  The other nine remaining off diagonal bilinear cuspidal representations are diagonal crossed cuspidal representations of gravitational fields between mixed shells  $ST-MG$, $ST-M$ and $MG-M$.  Knowing that gravitational fields can be transformed into electro-magnetic fields referring to \cite{Pie11}, we understand why these off diagonal bilinear cuspidal representations are not really relevant.\epr

\vfill

C. Pierre\\
Universit\'e de Louvain\\
Chemin du Cyclotron, 2\\
B-1348 Louvain-la-Neuve,  Belgium\\
pierre.math.be@gmail.com
\eject
}

\end{document}